\documentclass[11pt,reqno]{amsart} 
\makeatletter
\@namedef{subjclassname@2020}{\textup{2020} Mathematics Subject Classification}
\makeatother

\usepackage{amsmath,amssymb,amsthm,amscd}
\usepackage{mathtools,microtype,xcolor}
\usepackage[OT2,T1]{fontenc}
\usepackage{blkarray} \usepackage{mathrsfs,euscript} \usepackage[cal=boondoxo,scr=euler]{mathalfa}
\usepackage{graphicx} \usepackage[all,cmtip]{xy} \usepackage{cancel}
\usepackage[vcentermath]{youngtab} \usepackage{empheq} \usepackage{changepage}
\usepackage{bm} \usepackage{chngcntr}
\usepackage{calligra} 
\usepackage{scalerel,stackengine}  \usepackage{pifont}
\usepackage{booktabs}
\usepackage{adjustbox}  \usepackage{lscape}     \usepackage{longtable}  \usepackage{bbm}
\usepackage{enumitem}
\usepackage[colorlinks=true,urlcolor=blue,citecolor=red,linkcolor=blue,hyperfootnotes=true]{hyperref} \usepackage{array}
\usepackage{booktabs}

\let\mc\mathcal

\let\eu\EuScript
\let\nc\newcommand

\newtheorem{thm}{Theorem}[section]
\newtheorem{cor}[thm]{Corollary}
\newtheorem{lem}[thm]{Lemma}
\newtheorem{prop}[thm]{Proposition}

\theoremstyle{definition}
\newtheorem{defn}[thm]{Definition}
\newtheorem{rem}[thm]{Remark}

\newtheorem{example}[thm]{Example}
\providecommand{\exno}{}

\numberwithin{equation}{section}

\def\beq{\begin{equation}}
\def\eeq{\end{equation}}
\def\be{\begin{equation*}}
\def\ee{\end{equation*}}
\nc{\bea}{\begin{eqnarray*}}
\nc{\eea}{\end{eqnarray*}}

\let\phi\varphi

\def\C{{\mathbb C}}
\def\Z{{\mathbb Z}}
\def\Q{{\mathbb Q}}

\def\R{{\mathbb R}}

\counterwithin{figure}{section}
\counterwithin{table}{section}

\newlanguage\fakelanguage
\newcommand\cyr{\fontencoding{OT2}\fontfamily{wncyr}\selectfont
   \language\fakelanguage}
\DeclareTextFontCommand{\textcyr}{\cyr}

\DeclareMathOperator{\HOM}{\mathscr{H}\text{\kern -3pt {\calligra\large om}}\,}

\makeatletter
\newsavebox{\@brx}
\newcommand{\llangle}[1][]{\savebox{\@brx}{$\m@th{#1\langle}$}  \mathopen{\copy\@brx\kern-0.5\wd\@brx\usebox{\@brx}}}
\newcommand{\rrangle}[1][]{\savebox{\@brx}{$\m@th{#1\rangle}$}  \mathclose{\copy\@brx\kern-0.5\wd\@brx\usebox{\@brx}}}
\makeatother

\newcommand\xqed[1]{  \leavevmode\unskip\penalty9999 \hbox{}\nobreak\hfill
  \quad\hbox{#1}}
\newcommand\qetr{\xqed{$\triangle$}}
\newcommand\qrem{\xqed{$\spadesuit$}}

\stackMath
\newcommand\reallywidehat[1]{\savestack{\tmpbox}{\stretchto{  \scaleto{    \scalerel*[\widthof{\ensuremath{#1}}]{\kern.1pt\mathchar"0362\kern.1pt}    {\rule{0ex}{\textheight}}  }{\textheight}}{2.4ex}}\stackon[-6.9pt]{#1}{\tmpbox}}

\newcommand{\M}{\eu M}
\newcommand{\Mb}{\overline{\eu M}}

\newcommand{\Aut}{\operatorname{Aut}}
\newcommand{\Var}{\operatorname{Var}}
\newcommand{\orb}{\mathrm{orb}}

\allowdisplaybreaks

\usepackage{tikz}
\usepackage{graphicx}
\usepackage{multirow}
\usetikzlibrary{positioning}

\tikzset{
  graph vertex/.style={
    circle,
    draw,
    line width=.65pt,
    minimum size=4mm,
    inner sep=0pt,
    fill=white,
    font=\scriptsize
  },
  graph edge/.style={
    line width=.65pt
  },
  every loop/.style={
    line width=.65pt
  }
}

\newcommand{\cv}[3]{  \node[graph vertex] (#1) at #2 {#3};}

\newcommand{\graphcell}[2]{\begin{tikzpicture}[baseline=-.5ex]
  \path[use as bounding box] (-1.30,-.82) rectangle (1.30,.98);
  \node[font=\scriptsize] at (0,.82) {#1};
  #2
\end{tikzpicture}}

\newcommand{\doubleedge}[2]{  \draw[graph edge,bend left=16] (#1) to (#2);
  \draw[graph edge,bend right=16] (#1) to (#2);
}

\newcommand{\tripleedge}[2]{  \draw[graph edge] (#1) -- (#2);
  \draw[graph edge,bend left=24] (#1) to (#2);
  \draw[graph edge,bend right=24] (#1) to (#2);
}

\newcommand{\quadedge}[2]{  \draw[graph edge,bend left=32] (#1) to (#2);
  \draw[graph edge,bend left=11] (#1) to (#2);
  \draw[graph edge,bend right=11] (#1) to (#2);
  \draw[graph edge,bend right=32] (#1) to (#2);
}

\begin{document}

\title[Heavy--Light Wall Crossings and Motivic Asymptotics
of Moduli Spaces of Curves, I]{Heavy--Light Wall Crossings and Motivic Asymptotics
of Moduli Spaces of Curves\\[2mm]
Part I: Stabilization Graphs and Euler Characteristics}

\author[Giordano Cotti]{Giordano Cotti$\>^{\circ}$}

{\let\thefootnote\relax
\footnotetext{\vskip5pt 
\noindent
$^\circ\>$\textit{ E-mail}:  giordano.cotti@tecnico.ulisboa.pt}}


\maketitle
\begin{center}
\textit{ 
$^\circ\>$Grupo de F\'isica Matem\'atica \\
Departamento de Matem\'atica, Instituto Superior T\'ecnico\\
Av. Rovisco Pais, 1049-001 Lisboa, Portugal\/\\}
\end{center}
\vskip1cm

\begin{abstract}
In this paper, we develop an explicit calculation of ordinary
topological Euler characteristics of the coarse moduli spaces
underlying Hassett spaces in the minimal all-light chamber and
the Deligne--Mumford compactification. For genus $g\ge2$, we
organize the calculation through a decomposition indexed by
stabilization graphs and their decorations. The Burnside--Lefschetz
formula allows us to compute the resulting contributions, taking
into account graph symmetries and their fixed loci.

The calculation in genus two illustrates the method explicitly
and, in the Deligne--Mumford case, recovers the generating function
of G.\,Bini, G.\,Gaiffi, and M.\,Polito~\cite{BGP}, with the
corrections of G.\,Bini and J.\,Harer~\cite{BH}, providing a check
on the construction. We also carry out the calculation in genus
three in an appendix. Together with an elementary genus-one
calculation based on collision partitions and rational necklaces,
these formulas yield complete asymptotic expansions as the number
of labeled markings tends to infinity. The passage from all-light
to ordinary Euler characteristics is governed by the universal
change of variables $\Phi(z)=(2-z)e^z-2$, which we derive by a
direct wall-crossing argument.

For arbitrary fixed genus $g\ge2$, the same framework identifies
the leading contribution of rational trivalent graphs and gives
the ordinary asymptotic scale
$n!(e-2)^{-n}n^{(3g-5)/2}$. The leading constants are the symmetry
weights of these graphs; we derive their generating function,
recurrence, and large-genus behavior. We also determine the
asymptotic ratio of all-light to ordinary Euler characteristics,
answering a question of S.\,Kannan, S.\,Serpente,
and C.\,H.\,Yun~\cite{KSY}.
\end{abstract}

\tableofcontents

\section{Introduction}
This paper is concerned with the following question:
\vskip1,5mm
\noindent{\bf Question:}
{\it How does the stability condition affect the growth of the Euler characteristic of a moduli space of pointed curves?}
\vskip1,5mm
We investigate this question with the genus fixed and the number of labeled markings tending to infinity, comparing the classical Deligne--Mumford compactification with the all-light Hassett compactification.
\vskip1mm
This paper is the first in a series devoted to the asymptotics of
motivic invariants of moduli spaces of curves. Here we concentrate
on the ordinary topological Euler characteristic and on the analytic
structure of its exponential generating functions.
\medskip

\noindent 1.1.\,\,{\bf Weighted compactifications and the all-light chamber.}
Since the pioneering work of B.\,Riemann~\cite{Riemann1857}, the study of moduli spaces of algebraic curves has occupied a central position in mathematics. Their compactifications, developed much later, have revealed a remarkably rich interplay between geometry, topology, and combinatorics.

A fundamental milestone was the construction of the Deligne--Mumford compactification, subsequently extended to pointed curves by F.F.\,Knudsen~\cite{DeligneMumford,Knudsen}. The resulting moduli spaces $\Mb_{g,n}$ parametrize stable nodal curves with distinct marked points, and their boundary stratification is naturally described by stable dual graphs. It was later realized, however, that the Deligne--Mumford construction is not the only possible modular compactification of $\eu M_{g,n}$.

A particularly interesting alternative arises when marked points are allowed to coincide. An early example was introduced by A.\,Losev and Yu.I.\,Manin~\cite{LosevManin}: their genus-zero spaces parametrize chains of projective lines with two distinguished, or \emph{heavy}, markings at the ends, while the remaining \emph{light} markings may coincide. These so-called ``sausage'' compactifications were subsequently incorporated into the general theory of weighted stable curves developed by B.\,Hassett~\cite{Hassett}. His moduli spaces $\Mb_{g,\bm a}$, parametrized by admissible vectors of rational weights $\bm a$, include both the Deligne--Mumford and Losev--Manin spaces as special cases. Varying the weights determines which markings may collide and produces a natural system of reduction morphisms and wall-crossings. The related theory of painted stable curves and its operadic structures was further developed in~\cite{LosevManinExtended,ManinPainted}.

\begin{remark}
Hassett spaces are not the only alternative modular compactifications of $\eu M_{g,n}$. Other stability conditions may allow curves with singularities beyond nodes, as in Smyth's compactifications of pointed genus-one curves~\cite{SmythEllipticI}. A general framework for studying and classifying  modular compactifications, based on the notion of \emph{extremal assignments}, was developed by Smyth in~\cite{SmythClassification}, under suitable hypotheses.
\qrem
\end{remark}

In this paper, we focus on the opposite extreme to the classical
Deligne--Mumford stability condition: the \emph{all-light chamber}\footnote{This chamber is defined for $g\geq1$. In genus zero,
the stability condition $\sum_i a_i>2$ prevents the total weight
of all markings from being arbitrarily small.},
in which all markings have sufficiently small positive weights
and may coincide. Rational tails are excluded, while chains of
rational components carrying light markings may still occur.

For $g\ge2$, forgetting the markings and stabilizing associates
to each all-light curve an unmarked stable curve. We call its
dual graph the \emph{stabilization graph}. The finitely many
possible stabilization graphs, together with decorations
recording the inserted rational chains and the distribution
of the light labels, organize the all-light space into strata.
This gives a concrete framework for computing its ordinary
Euler characteristic: the Burnside--Lefschetz formula accounts
for graph symmetries and their fixed loci, and summing the
resulting contributions yields the generating function.

Our aim is to make this calculation explicit and to use the
resulting formulas to determine complete asymptotic expansions.
The genus-two calculation provides a check on the construction:
after passing to the Deligne--Mumford compactification, it
recovers the generating function of G.\,Bini, G.\,Gaiffi,
and M.\,Polito~\cite{BGP}, with the corrections of G.\,Bini
and J.\,Harer~\cite{BH}. A further calculation in genus three
is presented in Appendix~\ref{app:genus3}.

\medskip

\noindent1.2.\,\,{\bf Euler characteristics and previous results.}  The Euler characteristic is particularly well suited to the study of moduli spaces of curves: its additivity allows one to express the Euler characteristic as a sum of contributions from geometric strata. Throughout this paper, we consider the ordinary topological Euler characteristic of the coarse moduli spaces, with all markings labeled. This invariant must be distinguished from the orbifold Euler characteristic: computing the ordinary Euler characteristic of a finite quotient requires accounting for the fixed loci of group elements, whereas the orbifold Euler characteristic involves stabilizer weights.
\vskip1mm
The classical work of J.\,Harer and D.\,Zagier~\cite{HZ}
provides explicit formulas for the orbifold Euler
characteristics of moduli spaces of smooth pointed curves.
These formulas can subsequently be used to compute the
orbifold Euler characteristics of their Deligne--Mumford
compactifications, by decomposing the latter into strata
indexed by stable graphs.
The organization of such decompositions was developed
in the framework of modular operads by E.\,Getzler
and M.~M.\,Kapranov~\cite{GetzlerKapranov}.
Subsequently, G.\,Bini and J.\,Harer~\cite{BH}
obtained formulas for both the orbifold and the
ordinary Euler characteristics of the
Deligne--Mumford compactifications in arbitrary genus.
More recently, Z.\,Wang and J.\,Zhou~\cite{WangZhou}
developed recursion relations for the orbifold
Euler characteristics of Deligne--Mumford
compactifications in arbitrary genus.
\vskip1mm
In genus zero, S.\,Keel~\cite{Keel} described the intersection theory of $\overline{\eu M}_{0,n}$, while E.\,Getzler~\cite{Getzler} and Yu.~I.\,Manin~\cite{Manin} described the generating-series of its Poincar\'e polynomials as solutions of two functional equations. P.\,Aluffi, M.\,Marcolli, and E.\,Nascimento~\cite{AluffiMarcolliNascimento} later gave explicit formulas for the Grothendieck classes and Poincar\'e polynomials, as well as a Lambert $W$-function expression for the Euler-characteristic series. Further explicit formulas, recursions, and determinantal expressions for Betti numbers, Poincar\'e polynomials, and Euler characteristics were established in~\cite{Cotti}. From the asymptotic viewpoint, Manin~\cite[p.~403]{Manin} recorded a leading asymptotic equivalent for $\chi(\overline{\eu M}_{0,n})$, crediting its derivation to D.\,Zagier. In~\cite{Cotti}, this leading estimate was refined to a complete asymptotic expansion with explicitly computable coefficients.
\vskip1mm
Equivariant methods retain the symmetry information needed in more general graph calculations. E.\,Gorsky~\cite{Gorsky} computed symmetric-group equivariant Euler characteristics of smooth pointed moduli spaces, and A.\,Diaconu~\cite{Diaconu} developed effective formulas for the equivariant Euler characteristics of their stable compactifications using the modular-operad formalism. More recently, S.\,Kannan and T.~D.\,Song~\cite{KannanSong} introduced a graph-enumeration calculus for equivariant motivic invariants. Algorithmic enumeration of stable graphs is also addressed by S.\,Maggiolo and N.\,Pagani~\cite{MaggioloPagani}.
\vskip1mm
Explicit results in low genus provide further points of comparison.
In genus one, Getzler's semi-classical
approximation~\cite{GetzlerSemiClassical} describes the
genus-one part of the modular-operad construction in
terms of the genus-zero and genus-one inputs.
At the level of ordinary Euler characteristics,
Getzler also obtained the leading asymptotic equivalent
and the first relative correction as the number of
markings tends to infinity, in a result explicitly
credited to D.~Zagier
\cite[Corollary~4.2]{GetzlerSemiClassical}.
D.~Petersen~\cite{PetersenSemiClassical} gave a direct
combinatorial proof using necklaces and wreath-product
symmetric functions.
 In genus two, G.\,Bini, G.\,Gaiffi, and M.\,Polito~\cite{BGP} obtained a closed generating function for the ordinary Euler characteristics of $\overline{\eu M}_{2,n}$; corrections to this calculation are discussed in~\cite{BH}.
\vskip1mm
For weighted compactifications, S.\,Kannan, S.\,Serpente, and C.~H.\,Yun~\cite{KSY} expressed the generating functions for equivariant Hodge--Deligne polynomials of heavy/light Hassett spaces in terms of the corresponding series for classical pointed moduli spaces. Their formulas specialize to Euler characteristics. Y.-H.\,Kiem~\cite{Kiem} used wall crossing to reduce cohomological calculations for $\overline{\eu M}_{g,n}$ to those for the minimal-weight spaces $\overline{\eu M}_{g,(0^+)^n}$, and obtained explicit Hodge--Deligne and Poincar\'e polynomial formulas for these all-light spaces in genus one.
\vskip1mm
These results provide structural formulas and effective computational tools. A natural next challenge is to understand the analytic structure of the resulting generating functions and to derive complete asymptotic expansions for their coefficients. The “all-light” geometry described above is crucial in this regard: it provides the key simplification that lets us tackle the problem in a systematic way.

\medskip

\noindent 1.3\,\,\,{\bf Generating functions and their analytic structure.}
For every $g\geq1$, consider the exponential generating functions
\[
A_g^\chi(z):=
\sum_{n\geq0}
\chi\bigl(\Mb_{g,(0^+)^n}\bigr)\frac{z^n}{n!},
\qquad
G_g^\chi(t):=
\sum_{n\geq0}
\chi\bigl(\Mb_{g,n}\bigr)\frac{t^n}{n!},
\]
where $(0^+)^n$ denotes the all-light chamber, with both constant
terms set to zero when $g=1$. These series encode the Euler
characteristics of the two compactifications considered in this
paper.
\vskip1mm
A first remarkable fact is that they are related by a universal
change of variables, independent of the genus:
\begin{equation}\label{eq:intro-composition}
A_g^\chi(z)=G_g^\chi\bigl(\Phi(z)\bigr),
\qquad
G_g^\chi(t)=A_g^\chi\bigl(\Psi(t)\bigr),
\end{equation}
where\footnote{We denote by $f^{\langle-1\rangle}$ the compositional inverse of a function/series $f$.}
\[
\Phi(z)=(2-z)e^z-2,
\qquad
\Psi=\Phi^{\langle-1\rangle}.
\]
The inverse function admits the explicit representation
\[
\Psi(t)=2+W_{-1}\bigl(-(t+2)e^{-2}\bigr),
\]
where $W_{-1}$ denotes the branch of the Lambert $W$-function
satisfying $\Psi(0)=0$.

The identities \eqref{eq:intro-composition} are the Euler
specialization of the Hodge--Deligne transformation obtained by
S.\,Kannan, S.\,Serpente, and C.\,H.\,Yun
in~\cite[Corollary~B]{KSY}. In this paper, we give a direct
derivation from elementary wall-crossing identities, recovering
also the Euler specialization of Kiem's recursion
\cite[Theorem~3.1]{Kiem}.
Their significance for our purposes is that the ordinary
generating functions can be reconstructed from the geometrically
simpler all-light ones by a single composition.
\vskip1mm
Our first main result concerns the analytic structure of these
generating functions. The genus-one case exhibits a logarithmic
singularity, whereas in higher genus the all-light functions
belong to a particularly simple class of exponential-rational
functions. The ordinary functions then inherit a
Laurent--Puiseux singularity through the universal change
of variables.

\begin{theorem}[Structure of the Euler generating functions]
\label{intro:structure}
Let $\rho=e-2$. For every $g\geq1$, the series $G_g^\chi$
has radius of convergence $\rho$, with $\rho$ its unique
dominant singularity.

In genus one, the all-light generating function is
\[
A_1^\chi(z)
=
-\frac12\log(1-z)
+\frac{e^z-1}{3}
+\frac{3(e^{2z}-1)}8
+\frac{e^{3z}-1}{9}
+\frac{e^{4z}-1}{24}
-\frac z{12}.
\]
Moreover,
\[
G_1^\chi(t)
+\frac14\log\left(1-\frac t\rho\right)
\]
is holomorphic in $\sqrt{1-t/\rho}$ near $t=\rho$.

For every $g\geq2$, the function $A_g^\chi$ is a finite sum
of entire exponential polynomials\footnote{An
\emph{exponential polynomial} is a function of the form
$P(z)=\sum_{j=1}^N p_j(z)e^{\lambda_j z}$,
where $p_j\in\mathbb C[z]$ and $\lambda_j\in\mathbb C$.} multiplied by nonnegative
integer powers of $(1-z)^{-1}$. Its only finite singularity
is a pole of order $3g-3$ at $z=1$, with leading term
\[
A_g^\chi(z)
=
\frac{\gamma_g e^{-(2g-2)z}}{(1-z)^{3g-3}}
+O\bigl((1-z)^{-3g+4}\bigr),
\]
where
\[
\gamma_g:=
\sum_{\Gamma\in\mathfrak T_g}
\frac1{|\Aut\Gamma|},
\]
and $\mathfrak T_g$ denotes the set of connected rational
trivalent graphs\footnote{
A \emph{rational trivalent graph} is a connected multigraph
whose vertices all have genus decoration zero and valence three.
Loops and multiple edges are allowed, with each loop contributing
two to the valence of its vertex. The genus of such a graph
is its first Betti number $b_1(\Gamma)=|E(\Gamma)|-|V(\Gamma)|+1$.} of genus $g$.

The ordinary generating function $G_g^\chi$ admits a
convergent Laurent--Puiseux expansion at $t=\rho$, whose
leading term is
\[
G_g^\chi(t)
=
\frac{\gamma_g}
{e^{\frac{g-1}{2}}(2\rho)^{\frac{3g-3}{2}}}
\left(1-\frac t\rho\right)^{-\frac{3g-3}{2}}
\left(
1+O\left(\sqrt{1-\frac t\rho}\right)
\right).
\]
\end{theorem}

The statements collected in Theorem~\ref{intro:structure}
are established in different parts of the paper.
The universal change of variables is proved in
Theorem~\ref{thm:pde}, while the explicit genus-one
formula and its analytic consequences follow from
Theorem~\ref{thm:euler1}.
For $g\geq2$, the structure of the all-light
generating functions is established through
Proposition~\ref{prop:graphcontribution} and
Theorem~\ref{thm:alllightfixedg}.
The corresponding statements for the ordinary
generating functions then follow by composition
with the universal inverse function $\Psi$.

\begin{remark}[The genus-zero case]
Although the all-light chamber does not exist in genus zero,
the same universal inverse function $\Psi$ also appears
in the classical generating function
\[
G_0^\chi(t):=
\sum_{n\geq3}
\chi(\Mb_{0,n})\frac{t^n}{n!}.
\]
More precisely, one has the identity
\begin{equation}\label{eq:intro-genuszero}
1+t+(G_0^\chi)'(t)=e^{\Psi(t)}.
\end{equation}
This relation, established in
Proposition~\ref{prop:Psi-genuszero},
provides a further connection between the universal
change of variables governing weighted compactifications
and the classical genus-zero generating series
studied in~\cite{Manin,AluffiMarcolliNascimento,Cotti}.
\qrem
\end{remark}

\noindent{\bf Stabilization graphs and explicit computations.}
The structural results above follow from a finite decomposition
indexed by graphs. For $g\geq2$, forgetting the light markings
and stabilizing gives an unmarked stable curve. We call its
dual graph the \emph{stabilization graph}; only finitely many
such graphs occur in each genus. The original curve is
reconstructed by inserting chains of rational components
at the nodes of the stabilized curve and distributing the
light markings over the resulting components. We describe
these reconstructions in terms of decorated expansions.
\vskip1mm
The contribution of each graph can be computed from the
Euler characteristics of suitable moduli spaces of smooth
pointed curves, together with the action of its automorphism
group. Since we work with ordinary Euler characteristics
of coarse moduli spaces, the fixed loci of nontrivial
automorphisms must be taken into account through the
Burnside--Lefschetz formula of Lemma~\ref{lem:burnside}. The resulting expressions
are finite combinations of exponential polynomials
and powers of $(1-z)^{-1}$, as established in
Proposition~\ref{prop:graphcontribution}.
These contributions determine the exact generating functions,
including the terms that do not affect the leading singularity.

This construction also explains the leading singularity
in Theorem~\ref{intro:structure}. Each edge of a stable
graph contributes at most one factor $(1-z)^{-1}$, and
a graph of genus $g$ has at most $3g-3$ edges.
The maximal number of edges is attained precisely
by rational trivalent graphs, whose contributions
determine the coefficient $\gamma_g$.
Thus the dominant singular behavior is governed
by a finite collection of graphs, while the complete
generating functions are reconstructed by including
the contributions of the remaining graphs and
their automorphisms.
\vskip1mm
The procedure is entirely effective. We illustrate it explicitly
by computing the contributions of the seven stabilization graphs
of genus two. This gives closed expressions for both
$A_2^\chi$ and $G_2^\chi$. In the Deligne--Mumford case, we
recover the generating function of G.\,Bini, G.\,Gaiffi,
and M.\,Polito~\cite{BGP}, with the corrections of G.\,Bini
and J.\,Harer~\cite{BH}, providing a check on the construction.
We also carry out the calculation for the forty-two stabilization
graphs of genus three in Appendix~\ref{app:genus3}.
These calculations illustrate how the same finite
combinatorial procedure applies in arbitrary fixed genus.

In genus one, a related decomposition into smooth elliptic
curves and rational necklaces provides a direct geometric
explanation for the logarithmic term in $A_1^\chi$.

\medskip

\noindent 1.4\,\,\,{\bf Asymptotic expansions of Euler characteristics.}
The analytic structure described in
Theorem~\ref{intro:structure} has direct consequences
for the asymptotic growth of Euler characteristics.
The inverse function $\Psi$ has a square-root
singularity at $\rho=e-2$, so that a pole of order
$3g-3$ in the all-light variable becomes a
Laurent--Puiseux singularity of order
$(3g-3)/2$ in the ordinary variable.
In genus one, the corresponding singularity
is logarithmic.
\vskip1mm
By applying singularity-analysis techniques
\cite{FO,FS}, we obtain complete asymptotic
expansions for the ordinary Euler characteristics,
together with a remarkably simple asymptotic
description in the all-light chamber.

\begin{theorem}[Asymptotics of Euler characteristics]
\label{intro:asymptotics}
Let $\rho=e-2$. For every fixed $g\geq1$,
the sequence
\[
\frac{\chi(\Mb_{g,n})}{n!\rho^{-n}}
\]
admits a complete asymptotic expansion in
the powers
\(
n^{\frac{3g-5-k}{2}},\) with
\( k\geq0.
\)
Its leading term is given, in genus one, by
\[
\chi(\Mb_{1,n})
\sim
\frac{n!\rho^{-n}}{4n},
\]
and, for every $g\geq2$, by
\[
\chi(\Mb_{g,n})
\sim
\frac{
\gamma_g\,n!\rho^{-n}n^{\frac{3g-5}{2}}
}{
\Gamma\!\left(\frac{3g-3}{2}\right)
e^{\frac{g-1}{2}}
(2\rho)^{\frac{3g-3}{2}}
}.
\]

In the all-light chamber, the genus-one Euler
characteristics admit the explicit formula
\[
\chi(\Mb_{1,(0^+)^n})
=
\frac{(n-1)!}{2}
+\frac13+\frac38\,2^n
+\frac19\,3^n+\frac1{24}\,4^n
-\frac1{12}\delta_{n,1},
\qquad n\geq1,
\]
where $\delta_{n,1}$ denotes the Kronecker delta.
In particular,
\[
\chi(\Mb_{1,(0^+)^n})
=
\frac{(n-1)!}{2}+O(4^n)
\sim\frac{(n-1)!}{2}.
\]For every fixed $g\geq2$,
there exists a unique polynomial $P_g(n)\in\R[n]$ of degree
$3g-4$ such that
\[
\chi(\Mb_{g,(0^+)^n})
=
n!P_g(n)+O(b_g^n)
\]
for some $b_g>0$, with
\[
P_g(n)
=
\frac{\gamma_g e^{-2g+2}}{\Gamma(3g-3)}
n^{3g-4}
+O(n^{3g-5}).
\]
\end{theorem}

The genus-one expansion is established in
Theorem~\ref{thm:full1}. For $g\geq2$,
the ordinary leading equivalent is
Theorem~\ref{thm:ordinaryfixedg},
while the complete expansion follows from the
Laurent--Puiseux structure, as explained in
Remark~\ref{rem:fullfixedg}.
The all-light expansion is
Proposition~\ref{prop:alllightfullfixedg}.
\vskip1mm
The distinction between the two regimes is
particularly striking. For $g\geq2$, the principal
part of the Laurent expansion of $A_g^\chi(z)$
at $z=1$ determines precisely the polynomial
$P_g(n)$ appearing in the preceding theorem.
The remaining, entire part of $A_g^\chi(z)$
contributes only an exponentially bounded
term to the Euler characteristics.
Consequently, the all-light asymptotic expansion
at factorial scale contains only finitely many
terms.

By contrast, in the ordinary chamber, the
square-root singularity of the universal
inverse function $\Psi$ gives rise to a
Laurent--Puiseux expansion of $G_g^\chi$,
and hence to a generally infinite asymptotic
expansion in half-integer powers of $n^{-1}$.
\vskip1mm
A particularly interesting consequence concerns
the comparison of the two compactifications.
Kannan, Serpente, and Yun asked whether the
ratio of the all-light and ordinary Euler
characteristics tends to zero in every genus
$g\geq2$~\cite[Question~5.3]{KSY}.
Our results provide a precise asymptotic
equivalent, answering their question.

\begin{corollary}[Answer to Kannan--Serpente--Yun]
\label{intro:KSY}
For every fixed $g\geq2$, one has
\[
\frac{
\chi(\Mb_{g,(0^+)^n})
}{
\chi(\Mb_{g,n})
}
\sim
\left(\frac{2\rho}{e}\right)^{\frac{3g-3}{2}}
\frac{
\Gamma\!\left(\frac{3g-3}{2}\right)
}{
\Gamma(3g-3)
}
\rho^n n^{\frac{3g-3}{2}},
\]
where $\rho=e-2$. In particular, the ratio
tends to zero as $n\to\infty$.
\end{corollary}

This is Corollary~\ref{cor:KSYquestion}
in the main text.
Remarkably, the trivalent constant $\gamma_g$
cancels from the quotient: the relative
growth is governed by a universal exponential
factor $\rho^n$, together with an explicit
genus-dependent polynomial correction.
In genus one, the corresponding equivalent
$2\rho^n$ was already obtained
in~\cite[Corollary~5.2]{KSY}.

\medskip
\noindent{\bf Trivalent graphs and asymptotic constants.}
The coefficients $\gamma_g$ appearing in the preceding
formulas are themselves of independent
combinatorial interest. By definition, they are weighted
enumerations of connected trivalent graphs,
with each graph contributing the reciprocal
of the order of its automorphism group.
The same symmetry weights arise in the
enumeration of connected cubic Feynman
diagrams~\cite{CLP,BenderMilton,Borinsky}
and in the asymptotic enumeration of
connected graphs of fixed excess~\cite{FSS}.

We obtain the following generating function,
recurrence, and large-genus asymptotic equivalent.

\begin{theorem}[Trivalent constants]
\label{intro:trivalent}
As formal power series,
\[
\sum_{g\geq2}\gamma_gs^{g-1}
=
\log\left(
\sum_{m\geq0}
\frac{(6m-1)!!}{(2m)!6^{2m}}s^m
\right).
\]
In particular,
\[
\gamma_2=\frac{5}{24},
\qquad
\gamma_3=\frac{5}{16},
\]
and, for every $g\geq3$,
\[
\gamma_g=
\frac32(g-2)\gamma_{g-1}
+
\frac{3}{2(g-1)}
\sum_{h=2}^{g-2}
(h-1)(g-1-h)\gamma_h\gamma_{g-h}.
\]
Moreover,
\[
\gamma_g\sim
\frac1{2\pi}
\left(\frac32\right)^{g-1}(g-2)!
\qquad(g\to\infty).
\]
\end{theorem}

The generating-function identity, recurrence,
and large-genus equivalent are established
in Theorem~\ref{thm:cubicpartition},
Proposition~\ref{prop:gammagrecurrence},
and Theorem~\ref{thm:gammaglargegenus},
respectively. The limit $g\to\infty$
considered here is independent of the
fixed-genus regime $n\to\infty$ studied
above; no uniformity in the two parameters
is asserted.

\medskip
\noindent 1.5\,\,\,{\bf Organization of the paper.}
The paper is organized as follows.

Section~\ref{sec:euler-heavy-light} establishes the
basic framework for Euler characteristics of weighted
moduli spaces. After fixing the conventions concerning
coarse moduli spaces and ordinary Euler characteristics,
we review the relevant stability conditions and derive
an elementary Euler-characteristic identity for a single
wall crossing. Iterating this identity leads to the
heavy/light transformation of Theorem~\ref{thm:pde},
relating the ordinary and all-light generating functions
through the universal change of variables $\Phi$ and
its inverse $\Psi$. We also explain the genus-zero
interpretation of $\Psi$, connecting the transformation
with the classical generating function of
$\Mb_{0,n}$.

Section~\ref{sec:g1} is devoted to genus one.
We describe the geometry of all-light curves in terms
of smooth elliptic curves with possibly coincident
markings and rational necklaces. The contributions of
these two types of curves are computed separately,
leading to explicit closed formulas for the corresponding
Euler generating functions. Their analytic properties
are then used to derive the complete asymptotic expansion
of the ordinary Euler characteristics in genus one.

Section~\ref{sec:graph-formalism} develops the main
geometric and combinatorial construction for higher
genera. We introduce the stabilization-graph decomposition
and describe the strata lying over each unmarked
stable curve in terms of decorated graph expansions.
The contributions of these strata are expressed
through vertex and edge generating functions, while
the effects of graph automorphisms are handled
by fixed-locus calculations and the
Burnside--Lefschetz formula.
This yields an effective finite procedure for
computing the all-light Euler generating functions
and determining their analytic structure.

In Section~\ref{sec:g2}, we apply this procedure
to the seven stabilization graphs of genus two.
We compute their individual contributions and obtain
closed formulas for both the all-light and ordinary
Euler generating functions, recovering the 
formula of Bini--Gaiffi--Polito with its known
correction. We then derive explicit asymptotic
expansions in both stability chambers, illustrating
the general analytic mechanism in a concrete case.

Section~\ref{sec:fixed-genus} extends the asymptotic
analysis to arbitrary fixed genus $g\geq2$.
We show that the leading pole of the all-light
generating function is determined by rational
trivalent graphs, while nontrivial automorphism
sectors contribute only to lower pole orders.
The universal change of variables then provides
the corresponding ordinary asymptotics and an
explicit equivalent for the ratio of the Euler
characteristics in the two chambers.

Section~\ref{sec:gammag} investigates the trivalent
constants $\gamma_g$ governing these leading terms.
Using Gaussian pairings and the combinatorics
of connected cubic graphs, we derive their
generating function and recurrence, and determine
their asymptotic behavior as $g\to\infty$.
This provides a complementary large-genus
perspective on the constants arising from the
fixed-genus analysis.

Finally, Appendix~\ref{app:genus3} contains the
explicit genus-three calculation. We enumerate
the forty-two stabilization graphs, compute their exact
contributions, and obtain closed generating
functions and asymptotic expansions for both
compactifications, together with independent
consistency checks.

\vskip1,5mm
\noindent\textbf{Notation and conventions.}
We write $[n]=\{1,\ldots,n\}$ for $n\ge1$ and $[0]=\varnothing$,
and denote by $\mathfrak S_n$ the symmetric group on $[n]$.
For a formal power series $F(z)$, the notation $[z^n]F(z)$
denotes the coefficient of $z^n$.
The symbols $\square$, $\triangle$, and $\spadesuit$ mark the end
of a proof, an example, and a remark, respectively.

\vskip1,5mm
\noindent\textbf{Acknowledgements.}
This research was supported by Portuguese national funds
through the Fundação para a Ciência e a Tecnologia (FCT),
under project UID/00208/2025.

\noindent\textbf{Use of AI tools.}
AI tools accessed through IAedu, the platform provided by FCT through its digital services unit FCCN, were used for an extensive literature search, English-language revision, and editing parts of the \LaTeX{} source, including figure and table commands. These tools included OpenAI's ChatGPT~5.6 and Anthropic's Claude. Computations supporting this work were carried out in
\textit{Mathematica}, with assistance from ChatGPT in developing
the code used for the calculations reported in the appendix. The author takes full responsibility for the mathematical content, references, and final text.

\section{Euler characteristics and the heavy/light transformation}
\label{sec:euler-heavy-light}

\subsection{Coarse spaces and Euler characteristics}
\label{sec:euler-conventions}
All varieties and stacks are defined over $\mathbb C$. We use Betti
cohomology of the associated complex analytic spaces, with rational
coefficients. In particular, $H_c^i(X;\mathbb Q)=H^i\bigl(R\Gamma_c(X^{\mathrm{an}},\mathbb Q)\bigr)$ denotes cohomology
with compact supports. For a complex algebraic variety $X$, we write
\[
\chi(X)=\sum_i(-1)^i\dim_{\mathbb Q}H_c^i(X;\mathbb Q).
\]
\begin{remark}It would be more appropriate to denote 
by $\chi_c(X)$ the sum on the r.h.s.. It turns out, however, that for complex algebraic varieties $\chi_c(X)$ equals 
the ordinary topological Euler characteristic
$\chi(X)=\sum_i(-1)^i\dim_{\mathbb Q}H^i(X;\mathbb Q)$. See \cite[Theorem~2.1]{Virk}, applied to the structure morphism
$X\to\mathrm{pt}$ and the constant sheaf.\footnote{For $X$ smooth of pure complex dimension $d$, the equality also follows
from Poincar\'e duality: $H_c^i(X;\mathbb Q)$ is dual to
$H^{2d-i}(X;\mathbb Q)$.} We therefore denote both by \(\chi\). This equality does not identify the corresponding cohomology groups.
\qrem
\end{remark}

The compact-support definition gives additivity for finite locally
closed algebraic decompositions, by the cohomology localization
sequence, and multiplicativity for products, by the K\"unneth formula.
In particular,
\[
\chi(\mathbb A^m)=1,\qquad \chi(\mathbb P^{m-1})=m,
\qquad \chi(\mathbb G_m)=0.
\]
For a separated Deligne--Mumford stack $\mathcal X$ of finite type
with coarse moduli space $X$, we adopt the convention
$\chi(\mathcal X):=\chi(X)$. This convention is consistent with the
rational cohomology of the associated analytic stack: the coarse
moduli morphism induces an isomorphism\footnote{This is a comparison theorem,
not a definition of stack cohomology. Locally, the stack is presented
as a quotient by a finite group $G$, and the comparison rests on the
vanishing $H^q(BG;\mathbb Q)=0$ for $q>0$.}
\[
H^\bullet(X^{\mathrm{an}};\mathbb Q)
\simeq H^\bullet(\mathcal X^{\mathrm{an}};\mathbb Q);
\]
see \cite[Proposition~36]{Behrend}. 

The invariants $\chi(\mc X)$ should be distinguished from the orbifold Euler
characteristic $\chi_{\mathrm{orb}}(\mc X)$, which weights stabilizers.
For instance, for a finite group $G$, our convention gives
$\chi(BG)=1$, whereas the stabilizer-weighted invariant gives
$\chi_{\mathrm{orb}}(BG)=1/|G|$.
We will use the orbifold Euler characteristic $\chi_\orb$ only when recalling the 
Harer--Zagier formula for smooth component moduli spaces.

\subsection{Hassett spaces and heavy/light stability}\label{sec:hassett}

We retain $\Mb$ for the moduli stack and write $\overline M$ for
its coarse moduli space. Euler characteristics always follow the coarse-space
convention above.

Let $g,n\geq 0$, and
\[
 \bm a=(a_1,\ldots,a_n)\in\Q^n,\qquad
 0<a_i\le1,
 \qquad
 2g-2+\sum_{i=1}^n a_i>0.
\]
An $\bm a$-stable curve is a connected projective nodal curve $C$ of
arithmetic genus $g$, with labeled smooth markings
$p_1,\ldots,p_n$, such that marked points may coincide only when the
sum of their weights is at most $1$, and
\(
 \omega_C\Bigl(\sum_i a_ip_i\Bigr)
\)
is ample. Equivalently, on an irreducible component whose normalization
has genus $h$ and on which the normalization has $d$ preimages of
nodes (including self-nodes), stability is
\[
 2h-2+d+\sum_{p_i\text{ on the component}}a_i>0.
\]

Hassett constructs a smooth proper Deligne--Mumford stack
$\Mb_{g,\bm a}$ parametrizing such curves; its coarse moduli space is
projective of dimension $3g-3+n$, see \cite[Theorem~2.1]{Hassett}. When all
weights are equal to $1$, one recovers the usual moduli space
$\Mb_{g,n}$ of stable pointed curves \cite{DeligneMumford,Knudsen}.

As the weights vary, these spaces are related by reduction morphisms.
If $\bm b=(b_1,\ldots,b_n)\in\Q^n$ satisfies $b_i\le a_i$, then lowering the
weights permits further collisions and may destabilize some rational
components. Hassett's morphism
\(
 \rho_{\bm a,\bm b}:
 \Mb_{g,\bm a}\longrightarrow\Mb_{g,\bm b}
\)
contracts precisely the components that become unstable for the
weights $\bm b$; see \cite[Sections~4--5]{Hassett}.

For fixed $g$ and $n$, the admissible weight domain
is divided into finitely many chambers by the walls
$\sum_{i\in I}a_i=1$, for subsets $I$ of cardinality
at least two.
Inside a chamber the collision and stability inequalities are
unchanged, hence the corresponding Hassett space is independent, up to
canonical isomorphism, of the particular weights chosen there.
Crossing a wall changes the stability condition and is realized by a
reduction morphism.

We shall mainly use heavy/light chambers. We write
\(
 \Mb_{g,r|s}
\)
for the Hassett space with $r$ \emph{heavy} markings of weight $1$ and
$s$ \emph{light} markings of a common sufficiently small positive
weight $\epsilon$. Heavy markings cannot collide with any other marked
point, whereas light markings may collide provided that their total
weight is at most $1$.

For fixed $s>0$, we choose $\epsilon$ in the minimal positive-weight
chamber, for instance
\(
 0<\epsilon<\frac1s.
\)
Any two such choices give the same Hassett space. We therefore also
write the corresponding weight vector as
\(
 (1^r,(0^+)^s),
\)
and in the all-light case simply $(0^+)^s$.\footnote{The symbol $0^+$
denotes the chamber of sufficiently small positive weights, not the
weight $0$. The representative $\epsilon$ is chosen separately for
each fixed $s$.}

In this minimal chamber every subset of light markings may coincide,
while the positivity of $\epsilon$ still affects stability. 
The admissible heavy/light configurations are determined by
the condition $2g-2+r+s\epsilon>0$.
In the minimal chamber, this requires $r\geq3$, or
$r=2$ and $s\geq1$, when $g=0$; at least one marking
when $g=1$; and imposes no restriction when $g\geq2$.

\begin{lemma}[Stability of components in the all-light chamber]
\label{lem:all-light-stability}
Let $C$ be an all-light stable curve of genus $g$, with
$n\geq 0$ markings of weight $\epsilon>0$, where $n\epsilon<1$.
Let $C_v$ be an irreducible component of $C$, and denote by
$h(v)$ its geometric genus, by $d(v)$ the number of nodal
branches\footnote{Here a self-node contributes two nodal branches, so that
$d(v)$ is the valence of the corresponding vertex in the
dual graph, with loops counted twice.} on its normalization, and by $m(v)$ the number of
light markings lying on it.

Then the following statements hold.
\begin{enumerate}
\item If $h(v)=0$, stability is equivalent to \[
d(v)\geq 3,
\qquad\text{or}\qquad
d(v)=2\ \text{and}\ m(v)\geq 1.
\]
\item If $h(v)=1$, stability is equivalent to
\[
d(v)\geq1,
\qquad\text{or}\qquad
d(v)=0\ \text{and}\ m(v)\geq1.
\]

\item If $h(v)\geq 2$, the component is stable regardless
of the number of markings or nodal branches.
\end{enumerate}

Moreover, $C$ contains no rational tails, namely no
connected proper subcurves of arithmetic genus zero
meeting the closure of their complements at exactly one node.
\end{lemma}

\begin{proof}
The stability condition on $C_v$ is
\(
2h(v)-2+d(v)+\epsilon m(v)>0.
\)
Since $0\leq m(v)\leq n$ and $n\epsilon<1$, the assertions
follow by distinguishing the three possible ranges of
$h(v)$.

Indeed, for $h(v)=0$, if $d(v)\leq1$ then
\(
-2+d(v)+\epsilon m(v)
\leq -1+n\epsilon<0.
\)
If $d(v)=2$, the condition reduces to
$\epsilon m(v)>0$, while for $d(v)\geq3$ it
holds automatically.

For $h(v)=1$, the stability inequality becomes
$d(v)+\epsilon m(v)>0$, which is equivalent to
$d(v)\geq1$ or $m(v)\geq1$.
Finally, for $h(v)\geq2$ the inequality is
always satisfied.

For the last assertion, suppose that $T\subsetneq C$ is a connected
subcurve of arithmetic genus zero meeting the closure
of its complement at exactly one node.
The degree of the weighted dualizing sheaf on $T$ is
\[
\deg\left(
\omega_C\Bigl(\epsilon\sum_i p_i\Bigr)\big|_T
\right)
=-1+\epsilon m(T)<0,
\]
contradicting ampleness.
\end{proof}

\subsection{A single-wall Euler identity}
\label{sec:single-wall}

\begin{defn}
Let $\bm a\ge\bm b$ be admissible weight vectors. We say that they are
\emph{separated by a single wall indexed by $I$} if $I$ is the unique
subset of the marking labels satisfying
\[
\sum_{i\in I}b_i\le1<\sum_{i\in I}a_i.
\]
\end{defn}

\begin{lemma}\label{lem:motivic-wall}
Let $\bm a\ge\bm b$ be admissible weight vectors separated by a single
collision wall, indexed by $I$, with $k=|I|\ge2$. Let
$Z_I\subset\overline M_{g,\bm b}$ be the locus where the markings
indexed by $I$ coincide. Then
\begin{equation}\label{eq:euler-single-wall}
\chi(\overline M_{g,\bm a})-\chi(\overline M_{g,\bm b})
=(k-2)\chi(Z_I).
\end{equation}
\end{lemma}

\begin{proof}
We use Hassett's description of reduction morphisms
\cite[Theorem~4.1, Proposition~4.5, Remark~4.6 and
Corollary~4.7]{Hassett}, applied to coarse spaces. Put
$J=[n]\setminus I$ and $b_I=\sum_{i\in I}b_i$.
The single-wall assumption gives
\[
\sum_{i\in S}a_i\le1\quad(S\subsetneq I),\qquad
b_I+b_j>1\quad(j\in J).
\]
Otherwise the wall indexed by $S$ or $I\cup\{j\}$ would also be
crossed. Merging the markings in $I$ identifies
\[
Z_I\simeq\overline M_{g,((b_j)_{j\in J},b_I)}
\simeq\overline M_{g,((b_j)_{j\in J},1)}
\simeq\overline M_{g,((a_j)_{j\in J},1)}.
\]
Indeed, the merged marking cannot collide with any other marking.
A rational tail carrying it alone is unstable for either weight, and
a tail carrying it together with any other marking has total weight
greater than one. No wall involving only labels in $J$ is crossed.

For $k\ge3$, the reduction is an isomorphism away from $Z_I$.
Its inverse image $D_I$ is the divisor with a rational tail carrying
precisely the labels in $I$, and
\[
D_I\simeq\overline M_{0,((a_i)_{i\in I},1)}\times Z_I
\simeq\mathbb P^{k-2}\times Z_I.
\]
For completeness, the tail is irreducible: any terminal component
not containing its attaching point would carry a proper subset of
$I$ of weight greater than one. Put the attaching point at infinity.
The remaining moduli are
\[
(\mathbb A^k\setminus\Delta)/\operatorname{Aff}_1
\simeq\mathbb P^{k-2},\qquad \Delta=\{x\in\mathbb A^k\colon x_1=x_2=\dots=x_k\}.
\]
Translation sets the last coordinate to zero, and dilation
projectivizes the nonzero vector of differences. The attaching point
and labeled markings leave no nontrivial tail automorphisms;
automorphisms of the main component act trivially on this factor.
Thus the product description also holds on coarse spaces.
Additivity, and the isomorphism of the complements, give
\[
\chi(\overline M_{g,\bm a})-\chi(\overline M_{g,\bm b})
=\bigl(\chi(\mathbb P^{k-2})-1\bigr)\chi(Z_I)
=(k-2)\chi(Z_I).
\]
For $k=2$ the reduction is an isomorphism, and both sides vanish.
\end{proof}

\subsection{Recursions and generating functions}
\label{sec:wall}
Fix $g\ge1$ and write
\[
\chi_{r,s}:=\chi(\Mb_{g,r|s}).
\]
We set $\chi_{0,0}=0$ when $g=1$, while for $g\ge2$ it is
$\chi(\Mb_g)$. Define the formal exponential generating series
\begin{align*}
H^\chi_{g,r}(z)&=\sum_{s\ge0}\chi_{r,s}\frac{z^s}{s!},&
H_g^\chi(y,z)&=\sum_{r\ge0}H^\chi_{g,r}(z)\frac{y^r}{r!},\\
A_g^\chi(z)&=H_g^\chi(0,z),&G_g^\chi(t)&=H_g^\chi(t,0).
\end{align*}
Thus $A_g^\chi$ encodes the all-light chamber and $G_g^\chi$ the
ordinary Deligne--Mumford compactification.

\begin{proposition}[One-step Euler recursion]\label{prop:motivic-onestep}
For $g\ge1$, $r\ge1$ and $s\ge0$,
\begin{equation}\label{eq:euler-onestep}
\chi_{r,s}-\chi_{r-1,s+1}
=\sum_{m=2}^s\binom{s}{m}(m-1)\chi_{r,s-m}.
\end{equation}
Empty sums vanish.
\end{proposition}

\begin{proof}
The case $s=0$ is immediate, so assume $s\geq1$.
Choose positive rational weights
$\epsilon_1,\ldots,\epsilon_s$ with total less than one
and pairwise distinct subset sums. Choose also a rational
number $\eta>0$ such that
$\eta+\sum_i\epsilon_i<1$, and consider
\[
(1^{r-1},t,\epsilon_1,\ldots,\epsilon_s),
\qquad \eta\leq t\leq1.
\]
The endpoint chambers correspond to
$\Mb_{g,r|s}$ and $\Mb_{g,r-1|s+1}$, respectively.

Enumerate the subsets of $[s]$ as
$S_1,\ldots,S_{2^s}$ in increasing order of their
weight sums, with $S_1=\varnothing$, and set
\[
t_k=1-\sum_{i\in S_k}\epsilon_i,
\qquad k=1,\ldots,2^s.
\]
Then
$1=t_1>t_2>\cdots>t_{2^s}>\eta$.
As $t$ decreases, the collision walls are crossed
one at a time, precisely at $t=t_k$ for $k\geq2$.

Fix $k\geq2$, and write $S=S_k$, $m=|S|$.
Choose a rational number $\delta_k>0$ such that
\[
0<\delta_k<
\min\{t_{k-1}-t_k,\ t_k-t_{k+1}\},
\]
where $t_{2^s+1}:=\eta$.
Consider the two weight vectors
\(
\bm a_k^\pm
=(1^{r-1},t_k\pm\delta_k,
\epsilon_1,\ldots,\epsilon_s).
\)
They belong to the chambers immediately above
and below the $k$-th wall, respectively.

In the upper chamber, the distinguished marking
cannot collide with all markings indexed by $S$,
since their combined weight is $1+\delta_k>1$.
In the lower chamber, their combined weight is
$1-\delta_k<1$, so this collision becomes allowed.

The corresponding center
$Z_S\subset\overline M_{g,\bm a_k^-}$
is the locus where these $m+1$ markings coincide.
By the choice of $\delta_k$, for every $j\notin S$
we have
\(
\delta_k<t_k-t_{k+1}\leq\epsilon_j,
\)
since the wall corresponding to $S\cup\{j\}$
occurs at $t_k-\epsilon_j\leq t_{k+1}$.
Thus the merged marking, whose weight is $1-\delta_k$,
cannot collide with any remaining marking.
Replacing its weight by one does not change stability,
and consequently
\(
Z_S\simeq\overline M_{g,r|s-m}.
\)
By Lemma~\ref{lem:motivic-wall}, the contribution of
this wall crossing is therefore
\[
\chi(\overline M_{g,\bm a_k^+})
-\chi(\overline M_{g,\bm a_k^-})
=(m-1)\chi_{r,s-m}.
\]
Notice that for $2\leq k<2^s$, the weight vectors $\bm a_k^-$
and $\bm a_{k+1}^+$ lie in the same chamber, so the
corresponding moduli spaces are isomorphic.
Consequently, as the walls are crossed successively, we have
\[
\chi_{r,s}-\chi_{r-1,s+1}
=
\sum_{k=2}^{2^s}
\left[
\chi(\overline M_{g,\bm a_k^+})
-\chi(\overline M_{g,\bm a_k^-})
\right]
=
\sum_{k=2}^{2^s}
(|S_k|-1)\chi_{r,s-|S_k|}.
\]

For each fixed $m$, there are exactly $\binom{s}{m}$
subsets $S\subseteq[s]$ of cardinality $m$, each
appearing exactly once in the sum.
Their contributions are equal, since the corresponding
centers are isomorphic to the same heavy/light moduli
space, up to relabeling.
Grouping the terms by cardinality and observing that
the contributions for $m=1$ vanish gives
\eqref{eq:euler-onestep}.
\end{proof}

\begin{corollary}[Euler form of Kiem's recursion]\label{cor:motivic-kiem}
For $N=r+s>0$,
\[
\chi_{r,s}=\chi_{0,N}
+\sum_{j=1}^{\min(r,N-2)}\sum_{k=0}^{N-j-2}
\binom{N-j}{k}(N-j-k-1)\chi_{j,k}.
\]
\end{corollary}
\begin{proof}
Sum \eqref{eq:euler-onestep} over $j=1,\ldots,r$, with
$(r,s)=(j,N-j)$, and set $k=N-j-m$. This is the Euler specialization of Kiem's recursion
\cite[Theorem~3.1]{Kiem}.
\end{proof}

Put
\[
Q_1(z)=\sum_{m\ge2}(m-1)\frac{z^m}{m!}
=(z-1)e^z+1,
\]
and define
\begin{equation}\label{eq:phi-euler}
\Phi(z)=z-\sum_{m\ge2}(m-1)\frac{z^{m+1}}{(m+1)!}
=(2-z)e^z-2.
\end{equation}
Then $\Phi'=1-Q_1=(1-z)e^z$ and $\Phi=z+O(z^3)$.
Its unique formal compositional inverse is denoted by $\Psi$.

\begin{theorem}[Euler heavy/light transformation]\label{thm:pde}
For every $g\ge1$,
\begin{equation}\label{eq:euler-biv}
H_g^\chi(y,z)=G_g^\chi\bigl(y+\Phi(z)\bigr).
\end{equation}
Equivalently, it is the unique formal solution of
\[
\partial_zH_g^\chi=\Phi'(z)\partial_yH_g^\chi,
\qquad H_g^\chi(y,0)=G_g^\chi(y).
\]
In particular,
\begin{equation}\label{eq:composition-euler}
A_g^\chi(z)=G_g^\chi(\Phi(z)),\qquad
G_g^\chi(t)=A_g^\chi(\Psi(t)).
\end{equation}
\end{theorem}
\begin{proof}
Multiply \eqref{eq:euler-onestep} by $z^s/s!$ and sum. Binomial
convolution gives
\[
\Phi'(z)H^\chi_{g,r}(z)=(H^\chi_{g,r-1})'(z).
\]
Summation against $y^{r-1}/(r-1)!$ yields the PDE. With
$x=\Phi(z)$ and $J(y,x)=H_g^\chi(y,\Psi(x))$, the formal chain
rule gives $\partial_xJ=\partial_yJ$. Writing
$J=\sum_{n\ge0}J_n(y)x^n$ gives
$(n+1)J_{n+1}=J_n'$ and $J_0=G_g^\chi$; thus the unique solution is
$J(y,x)=G_g^\chi(y+x)$. This proves the composition identities.
\end{proof}

\begin{remark}[Relation to Kannan--Serpente--Yun]
Theorem~\ref{thm:pde} is the Euler specialization of the scalar
Hodge--Deligne transformation in \cite[Corollary~B]{KSY}.
The direct wall-crossing proof above makes the Euler argument
self-contained. The extension to coarse Grothendieck classes and
general motivic realizations belongs to the companion study.
\qrem\end{remark}

\subsection{The genus-zero comparison}
\label{sec:g0comparison}

There is no all-light space $\Mb_{0,(0^+)^n}$ in the minimal-weight
chamber, since stability of a genus-zero curve requires total weight
strictly larger than $2$. Thus genus zero does not admit a direct
analogue of the all-light generating functions $A_g^\chi$ introduced above.

For ordinary markings, set
\[
G_0^\chi(t)
:=
\sum_{n\ge3}
\chi(\Mb_{0,n})\frac{t^n}{n!}.
\]
The generating function used in \cite{Cotti} is
\[
Y(t)
:=
1+t+(G_0^\chi)'(t)
=
1+t+
\sum_{n\ge3}
\chi(\Mb_{0,n})\frac{t^{n-1}}{(n-1)!}.
\]
The genus-zero series gives an intrinsic interpretation of the
inverse heavy/light coordinate.

\begin{proposition}[Genus-zero interpretation of the inverse coordinate]
\label{prop:Psi-genuszero}
\[
Y(t)=e^{\Psi(t)}.
\]
\end{proposition}
\begin{proof}
The genus-zero recursion gives
\[
Y'(t)=\frac{Y(t)}{2+t-Y(t)},\qquad Y(0)=1;
\]
see \cite{Manin,AluffiMarcolliNascimento,Cotti}.
Writing $W(t)=\log Y(t)$ gives $W'=1/(2+t-e^W)$ and $W(0)=0$.
On the other hand,
\[
\Phi'(w)=2+\Phi(w)-e^w.
\]
Formal inversion shows that $\Psi$ satisfies the same differential
equation and initial condition as $W$. Uniqueness implies $W=\Psi$.
\end{proof}

The analytic description of $\Psi$ in
Section~\ref{sec:sing1} shows that the unique dominant singularity of
$Y$ is
\[
\rho=e-2.
\]

The complete asymptotic expansion of
$\chi(\Mb_{0,n})$ was obtained in \cite[Theorem~2.18]{Cotti}. We shall
only need its leading term:
\begin{equation}
\chi(\Mb_{0,n})
\sim
C_0\,n!\rho^{-n}n^{-5/2},
\qquad
C_0
=
\sqrt{\frac{e}{2\pi}}\,\rho^{3/2}.
\label{eq:g0leading}
\end{equation}
Thus genus zero already exhibits the exponential scale $\rho^{-n}$
which will reappear in every fixed genus. 
\section{Genus one: collision strata, rational necklaces, and Euler asymptotics}
\label{sec:g1}

\subsection{Smooth elliptic curves and collisions}
\label{sec:g1-collisions}

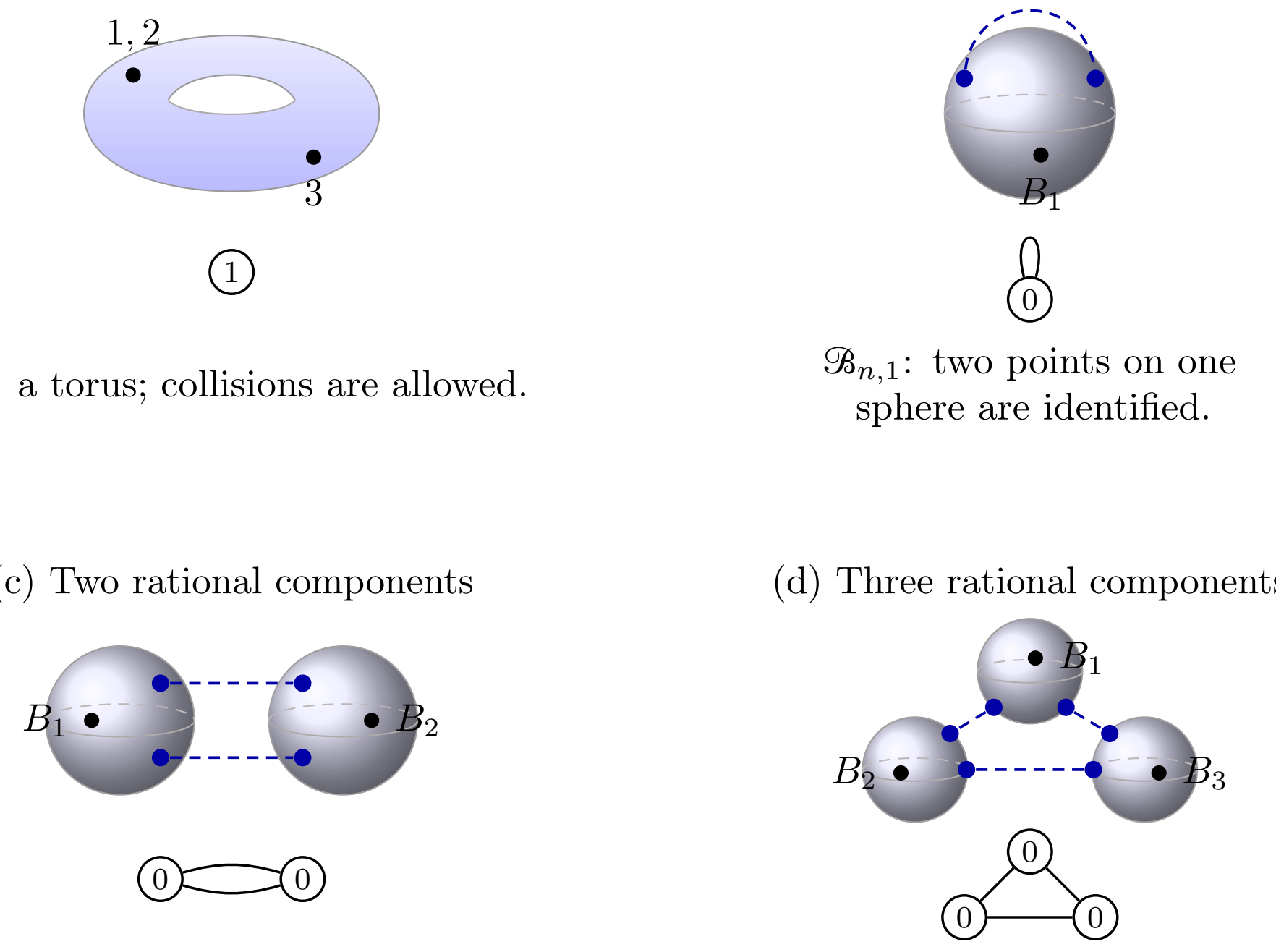
\begin{figure}[ht]
\centering
\begin{tikzpicture}[
  every node/.style={font=\small},
  identification/.style={densely dashed,draw=blue!65!black,line width=.7pt},
  node point/.style={circle,fill=blue!65!black,inner sep=1.6pt},
  light point/.style={circle,fill=black,inner sep=1.4pt}]
\newcommand{\topsphere}[2]{  \begin{scope}[shift={#1}]
  \shade[ball color=blue!9] (0,0) circle (#2);
  \draw[gray!75] (0,0) circle (#2);
  \draw[gray!55,densely dashed] (-#2,0)
    arc[start angle=180,end angle=0,x radius=#2,y radius=.22*#2];
  \draw[gray!65] (-#2,0)
    arc[start angle=180,end angle=360,x radius=#2,y radius=.22*#2];
  \end{scope}}
\begin{scope}[shift={(-3.65,3.1)}]
\node at (0,1.25) {(a) Smooth elliptic curve};
\shade[top color=blue!5,bottom color=blue!30]
 (-1.35,0) .. controls (-1.35,.95) and (1.35,.95) .. (1.35,0)
 .. controls (1.35,-.95) and (-1.35,-.95) .. (-1.35,0);
\draw[gray!80] (-1.35,0)
 .. controls (-1.35,.95) and (1.35,.95) .. (1.35,0)
 .. controls (1.35,-.95) and (-1.35,-.95) .. (-1.35,0);
\fill[white] (-.58,.12)
 .. controls (-.42,.43) and (.42,.43) .. (.58,.12)
 .. controls (.42,-.05) and (-.42,-.05) .. (-.58,.12);
\draw[gray!80] (-.58,.12)
 .. controls (-.42,.43) and (.42,.43) .. (.58,.12);
\draw[gray!65] (-.58,.12)
 .. controls (-.42,-.05) and (.42,-.05) .. (.58,.12);
\node[light point,label=above:{$1,2$}] at (-.9,.35) {};
\node[light point,label=below:{$3$}] at (.75,-.4) {};
\cv{topell}{(0,-1.45)}{1}
\node[text width=6.1cm,align=center] at (0,-2.5)
 {$\mathcal I_n$: a torus; collisions are allowed.};
\end{scope}
\begin{scope}[shift={(3.65,3.1)}]
\node at (0,1.25) {(b) One self-nodal rational component};
\topsphere{(0,0)}{.78}
\node[node point] (selfa) at (-.6,.32) {};
\node[node point] (selfb) at (.6,.32) {};
\draw[identification] (selfa)
 .. controls (-.5,1.12) and (.5,1.12) .. (selfb);
\node[light point,label=below:{$B_1$}] at (.1,-.38) {};
\cv{topone}{(0,-1.7)}{0}
\draw[graph edge] (topone) edge[loop above] (topone);
\node[text width=6.1cm,align=center] at (0,-2.5)
 {$\mathcal B_{n,1}$: two points on one sphere are identified.};
\end{scope}
\begin{scope}[shift={(-3.65,-2.45)}]
\node at (0,1.25) {(c) Two rational components};
\topsphere{(-1.02,0)}{.68} \topsphere{(1.02,0)}{.68}
\node[node point] (twoleftup) at (-.65,.34) {};
\node[node point] (tworightup) at (.65,.34) {};
\node[node point] (twoleftdown) at (-.65,-.34) {};
\node[node point] (tworightdown) at (.65,-.34) {};
\draw[identification] (twoleftup)--(tworightup);
\draw[identification] (twoleftdown)--(tworightdown);
\node[light point,label=left:{$B_1$}] at (-1.28,0) {};
\node[light point,label=right:{$B_2$}] at (1.28,0) {};
\cv{toptwoa}{(-.65,-1.45)}{0} \cv{toptwob}{(.65,-1.45)}{0}
\doubleedge{toptwoa}{toptwob}
\node[text width=6.1cm,align=center] at (0,-2.5)
 {$\mathcal B_{n,2}$: two spheres joined at two nodes.};
\end{scope}
\begin{scope}[shift={(3.65,-2.45)}]
\node at (0,1.25) {(d) Three rational components};
\topsphere{(0,.45)}{.48}
\topsphere{(-1.05,-.45)}{.48} \topsphere{(1.05,-.45)}{.48}
\node[node point] (triupperleft) at (-.33,.12) {};
\node[node point] (trileftupper) at (-.73,-.12) {};
\node[node point] (triupperright) at (.33,.12) {};
\node[node point] (trirightupper) at (.73,-.12) {};
\node[node point] (trileftinner) at (-.58,-.45) {};
\node[node point] (trirightinner) at (.58,-.45) {};
\draw[identification] (triupperleft)--(trileftupper);
\draw[identification] (triupperright)--(trirightupper);
\draw[identification] (trileftinner)--(trirightinner);
\node[light point,label=right:{$B_1$}] at (.05,.57) {};
\node[light point,label=left:{$B_2$}] at (-1.18,-.48) {};
\node[light point,label=right:{$B_3$}] at (1.18,-.48) {};
\cv{topthreea}{(0,-1.2)}{0}
\cv{topthreeb}{(-.6,-1.8)}{0} \cv{topthreec}{(.6,-1.8)}{0}
\draw[graph edge] (topthreea)--(topthreeb)--(topthreec)--(topthreea);
\node[text width=6.1cm,align=center] at (0,-2.5)
 {$\mathcal B_{n,3}$: three spheres joined in a cycle.};
\end{scope}
\end{tikzpicture}
\caption{Topological models of the genus-one loci.
A smooth elliptic curve is topologically a torus, whereas
the normalization of each rational component is a sphere.
Blue points connected by dashed lines are identified
to form nodes; the dashed lines indicate identifications,
not additional segments of the curve.
Black points indicate light markings.
The small diagrams record the incidence relations among
irreducible components, whose systematic graph-theoretic
description is postponed to
Section~\ref{sec:graph-formalism}. The numbers inside the diagram vertices indicate
the geometric genera of the corresponding components.
Each $B_i$ denotes a nonempty set of labels carried by
a rational component, without requiring the corresponding
markings to coincide.}
\label{fig:g1-loci}
\end{figure}

Fix $n\geq1$ and choose a common light weight $\epsilon$ with
$n\epsilon<1$. We first describe the curves parametrized by
$\Mb_{1,(0^+)^n}$.
\vskip1 mm
By Lemma~\ref{lem:all-light-stability}, no rational tails
can occur. If a curve of arithmetic genus\footnote{For a connected
nodal curve $C$ with irreducible components
$C_1,\ldots,C_\nu$ and $\delta$ nodes, the arithmetic genus is
\(
p_a(C)=\sum_{i=1}^{\nu}g(\widetilde C_i)
+\delta-\nu+1,
\)
where $\widetilde C_i$ denotes the normalization of $C_i$.
} one has a component of geometric genus one, all other
components, if any, belong to rational trees attached to it.
Such trees would be rational tails, hence the curve must be
smooth and irreducible.
\vskip1mm
Otherwise, every irreducible component is rational.
The arithmetic genus condition implies that the components
are arranged in a single cycle: any branches departing
from this cycle would form rational tails. Such curves
are called \emph{rational necklaces}.
Each rational component has exactly two nodal branches
and must therefore carry at least one light marking,
again by Lemma~\ref{lem:all-light-stability}.
\vskip1mm
Let $\mathcal I_n$ denote the locus of smooth elliptic curves
with $n$ labeled light markings, where collisions are allowed,
and let $\mathcal B_n$ denote the locus of rational necklaces.
The former is open and the latter is its closed complement.
Hence
\[
\overline M_{1,(0^+)^n}
=\mathcal I_n\sqcup\mathcal B_n,
\qquad A_1^\chi=I^\chi+B^\chi,
\]
where
\[
I^\chi(z)=\sum_{n\geq1}\chi(\mathcal I_n)\frac{z^n}{n!},
\qquad
B^\chi(z)=\sum_{n\geq1}\chi(\mathcal B_n)\frac{z^n}{n!}
\]
are their Euler exponential generating series.
Both have zero constant term.
The two types of curves are illustrated in Figure~\ref{fig:g1-loci}.

On $\mathcal I_n$, the collision patterns of the light
markings correspond to partitions of $[n]$ into
nonempty blocks, each block representing a collection
of coincident markings. For example, when $n=3$,
the partition $\{\{1,3\},\{2\}\}$ defines a stratum
isomorphic to $M_{1,2}$: the two distinct marked
points carry the blocks $\{1,3\}$ and $\{2\}$.
Since the blocks are distinguished by their labels,
no additional symmetric-group quotient is required.

\begin{lemma}[The elliptic collision contribution]\label{lem:elliptic-collisions}
\begin{equation}\label{eq:euler-interior-integrated}
I^\chi(z)=\frac{e^z-1}{3}+\frac{e^{2z}-1}{4}
+\frac{e^{3z}-1}{9}-\frac{e^{4z}-1}{48}-\frac z{12}.
\end{equation}
\end{lemma}
\begin{proof}
A partition of the $n$ labels into $k$ nonempty collision blocks
defines a stratum isomorphic to $M_{1,k}$. 
Thus
\[
\chi(\mathcal I_n)=\sum_{k=1}^n
\genfrac{\{}{\}}{0pt}{}{n}{k}\chi(M_{1,k}),
\qquad
I^\chi(z)=\sum_{k\ge1}\chi(M_{1,k})\frac{(e^z-1)^k}{k!}.
\]
Here the braces denote Stirling numbers of the second kind.
The ordinary Euler characteristics of the smooth elliptic spaces are
\[
\chi(M_{1,1})=\chi(M_{1,2})=1,\qquad
\chi(M_{1,3})=\chi(M_{1,4})=0,
\]
and $\chi(M_{1,k})=(-1)^k(k-1)!/12$ for $k\ge5$;
see \cite[Section~4]{BGP}. The last formula is the genus-one
Harer--Zagier formula \cite{HZ}; with five distinct markings no
nontrivial automorphism remains.
For
$F(x)=\sum_{k\ge1}\chi(M_{1,k})x^k/k!$, correcting the first four
terms of the formal Harer--Zagier series gives
\[
F(x)=-\frac1{12}\log(1+x)
+\frac{13}{12}x+\frac{11}{24}x^2+\frac1{36}x^3-\frac1{48}x^4.
\]
Substitute $x=e^z-1$. The logarithm becomes $z$, and expanding the
polynomial gives \eqref{eq:euler-interior-integrated}.
\end{proof}

\subsection{Rational necklaces and the short-cycle quotients}
\label{sec:g1-necklaces}
Write $\mathcal B_n=\bigsqcup_{j=1}^n\mathcal B_{n,j}$ according
to the number of components. A rational component with two
distinguished nodal branches and $m\ge1$ light labels has moduli space
$(\mathbb G_m)^{m-1}$: put the nodes at $0,\infty$, normalize one
light coordinate to one, and allow the remaining coordinates to
coincide freely in $\mathbb G_m$.

The following is a special case of the Burnside–Lefschetz formula, which we shall prove in greater generality in Lemma \ref{lem:burnside}.

\begin{lemma}[Euler characteristic of the inversion quotient]
\label{lem:euler-inversion-quotient}
For $m\ge1$, let $Q_m=(\mathbb G_m)^m/\langle\iota_m\rangle$,
where $\iota_m$ inverts every coordinate and let $Q_0$ be a point. Then
\[
\chi(Q_m)=\frac12\bigl(\mathbf1_{m=0}+2^m\bigr),\qquad m\geq 0.
\]
\end{lemma}
\begin{proof}
Let $\iota$ be an algebraic involution of a complex variety $X$,
and put $F=X^\iota$ and $U=X\setminus F$.
The action on $U$ is free, so $U\to U/\langle\iota\rangle$ is a
two-sheeted covering. Euler characteristic with compact supports
is multiplicative by the degree of a finite covering, and therefore
\(
\chi_c(U)=2\chi_c(U/\langle\iota\rangle).
\)
The image of $F$ in the quotient is isomorphic to $F$. Additivity gives
\[
\chi_c(X)=\chi_c(F)+\chi_c(U),\qquad
\chi_c(X/\langle\iota\rangle)
=\chi_c(F)+\chi_c(U/\langle\iota\rangle).
\]
Combining these identities, and using
Section~\ref{sec:euler-conventions}, yields
\[
\chi(X/\langle\iota\rangle)
=\frac{\chi(X)+\chi(X^\iota)}2.
\]
For simultaneous inversion on $(\mathbb G_m)^m$, the fixed locus is
$\{\pm1\}^m$, whereas
$\chi((\mathbb G_m)^m)=\mathbf1_{m=0}$.
The claim follows; it also covers the trivial action for
$m=0$.
\end{proof}

\begin{example}For example, $Q_1$ is isomorphic to $\mathbb A^1$, with quotient
coordinate $t+t^{-1}$. Its Euler characteristic is
$(0+2)/2=1$: the two fixed points supply the contribution that would
be lost by merely dividing $\chi(\mathbb G_m)$ by two.
\qetr
\end{example}

\begin{lemma}[The rational necklace contribution]\label{lem:euler-necklaces}
\begin{equation}\label{eq:euler-boundary}
B^\chi(z)=-\frac12\log(1-z)
+\frac{e^{2z}-1}{8}+\frac{e^{4z}-1}{16}.
\end{equation}
\end{lemma}
\begin{proof}
For one component, exchanging the two branches of the node gives
the simultaneous inversion quotient $Q_{n-1}$. Consequently,
\[
\chi(\mathcal B_{n,1})
=\tfrac12(\mathbf1_{n=1}+2^{n-1}),\qquad
B_1^\chi(z)=\frac z2+\frac{e^{2z}-1}{4}.
\]
For two components, the label sets form an unordered bipartition,
with $2^{n-1}-1$ choices. Each nonempty label set fixes its component;
the remaining graph symmetry exchanges the two connecting nodes and
acts by simultaneous inversion on $(\mathbb G_m)^{n-2}$.
Hence, for $n\ge2$,
\[
\chi(\mathcal B_{n,2})
=\frac{2^{n-1}-1}{2}(\mathbf1_{n=2}+2^{n-2}),
\qquad
B_2^\chi(z)=\frac{z^2}{4}+\frac{(e^{2z}-1)^2}{16}.
\]

For $j\ge3$, there are
$\genfrac{\{}{\}}{0pt}{}{n}{j}$
ways to partition the $n$ labels into $j$ nonempty sets.
Once such a partition is fixed, its $j$ blocks can be arranged
in $j!$ linear orders. Since cyclic rotations describe the
same arrangement along the necklace, this gives
$(j-1)!$ cyclic orders. Moreover, the cycle has no
preferred orientation, so each cyclic order must be
identified with its reversal.\footnote{For example, assume $j=4$. If the
blocks are $A,B,C,D$, the cyclic arrangements
$(A,B,C,D)$ and $(B,C,D,A)$ differ only by a rotation,
while $(A,B,C,D)$ and $(A,D,C,B)$ differ by a reversal of cycle orientation.
All three therefore describe the same unoriented cycle.}
Hence the number of
strata is
\[
\genfrac{\{}{\}}{0pt}{}{n}{j}\frac{(j-1)!}{2}.
\]

Since every component carries at least one labeled marking,
no nontrivial graph automorphism preserves all the labels.
Each stratum is $(\mathbb G_m)^{n-j}$, whose Euler
characteristic vanishes unless $n=j$.
Thus, the contribution to the exponential generating series
for fixed $j$ is
\[
\frac{(j-1)!}{2}\frac{z^j}{j!}
=\frac{z^j}{2j}.
\]
Summing over all cycle lengths, including the cases
$j=1,2$ treated above, gives
\[
B^\chi(z)=\frac z2+\frac{e^{2z}-1}{4}
+\frac{z^2}{4}+\frac{(e^{2z}-1)^2}{16}
+\sum_{j\ge3}\frac{z^j}{2j},
\]
which simplifies to \eqref{eq:euler-boundary}.
\end{proof}

\subsection{Closed Euler generating functions}
\label{sec:euler1}
The collision and necklace calculations give the all-light series
directly. The ordinary series is then obtained from
Theorem~\ref{thm:pde}.
\begin{theorem}[Closed genus-one Euler formulas]
\label{thm:euler1}
We have
\begin{align}
A_1^{\chi}(z)
={}&
-\frac12\log(1-z)
+\frac{e^z-1}{3}
+\frac{3(e^{2z}-1)}8
+\frac{e^{3z}-1}{9}
+\frac{e^{4z}-1}{24}
-\frac z{12},
\label{eq:A1}\\
G_1^{\chi}(t)
={}&
A_1^{\chi}\!
\left(
2+W_{-1}\bigl(-(t+2)e^{-2}\bigr)
\right),
\label{eq:Lambert}\\
\chi(\Mb_{1,n})
={}&
(n-1)![z^{n-1}]
(A_1^\chi)'(z)
\left(
\frac{z}{(2-z)e^z-2}
\right)^n .
\label{eq:LB}
\end{align}
Here $W_{-1}$ denotes the branch of the Lambert $W$-function satisfying
\(
W_{-1}(-2e^{-2})=-2.
\)
Moreover, for every $n\geq1$,
\begin{equation}
\chi(\Mb_{1,(0^+)^n})
=
\frac{(n-1)!}{2}
+\frac13
+\frac38\,2^n
+\frac19\,3^n
+\frac1{24}\,4^n
-\frac{\mathbf 1_{n=1}}{12}.
\label{eq:euler-light-coeff}
\end{equation}
\end{theorem}

\begin{proof}
Adding \eqref{eq:euler-interior-integrated} and
\eqref{eq:euler-boundary} proves \eqref{eq:A1}.

For the ordinary spaces, the universal heavy/light transformation
\eqref{eq:composition-euler} gives
\[
G_1^\chi(t)
=
A_1^\chi(\Psi(t)),
\]
where $\Psi=\Phi^{\langle-1\rangle}$ and
\[
\Phi(z)=(2-z)e^z-2.
\]
Solving $t=\Phi(z)$ for $z$ gives
\[
\Psi(t)
=
2+W\bigl(-(t+2)e^{-2}\bigr).
\]
The condition $\Psi(0)=0$ selects the branch $W_{-1}$, proving
\eqref{eq:Lambert}; for the standard properties of the Lambert
function and its branches see \cite{LambertRef}.

Finally, Lagrange--B\"urmann inversion applied to
$t=\Phi(z)$ yields
\[
[t^n]A_1^\chi(\Psi(t))
=
\frac1n[z^{n-1}]
(A_1^\chi)'(z)
\left(\frac z{\Phi(z)}\right)^n.
\]
Multiplication by $n!$ gives \eqref{eq:LB}. Formula
\eqref{eq:euler-light-coeff} follows directly by extracting
coefficients from \eqref{eq:A1}.
\end{proof}

\begin{remark}[Relation with Bini--Gaiffi--Polito]
\label{rem:BGP-genus-one}
Put
\[
D=e^z-1,\qquad E=z.
\]
The entire part of \eqref{eq:A1} can equivalently be written as
\[
\frac{19}{12}D
+\frac{23}{24}D^2
+\frac5{18}D^3
+\frac1{24}D^4
-\frac E{12}.
\]
In these variables, \eqref{eq:Lambert} agrees with the genus-one
generating function of Bini--Gaiffi--Polito
\cite[Theorem~1.1]{BGP}. The derivation above computes the all-light series directly from
collision and necklace strata; composition with the Euler heavy/light
coordinate then recovers their ordinary genus-one expression.
\qrem\end{remark}

\subsection{Full Euler-characteristic asymptotic expansion}
\label{sec:sing1}

We now extract the complete asymptotic expansion of
$\chi(\Mb_{1,n})$ from the closed formula for $G_1^{\chi}$.

Write
\[
 A_1^{\chi}(z)=-\frac12\log(1-z)+R(z),
\]
where
\[
 R(z)
 =
 \frac{e^z-1}{3}
 +\frac{3(e^{2z}-1)}8
 +\frac{e^{3z}-1}{9}
 +\frac{e^{4z}-1}{24}
 -\frac z{12}
\]
is entire, and set
\[
 \rho:=e-2,\qquad
 \zeta:=1-\frac{t}{\rho},\qquad
 a:=\sqrt{\frac{2\rho}{e}}.
\]

The coefficients entering the expansion below can be computed
effectively as follows. Let
\[
 \mathfrak d(u)=u+\xi_2u^2+\xi_3u^3+\cdots
\]
be the analytic solution at the origin of
\begin{equation}\label{eq:def-d}
 2\bigl(1-(1+\mathfrak d(u))e^{-\mathfrak d(u)}\bigr)=u^2.
\end{equation}
Its coefficients are obtained successively by comparing powers of $u$
in \eqref{eq:def-d}; equivalently, Lagrange inversion gives
\begin{equation}\label{eq:dcoeff}
 [u^m]\mathfrak d(u)
 =
 \frac1m[x^{m-1}]
 \left(
 \frac{x^2}{2(1-(1+x)e^{-x})}
 \right)^{m/2}.
\end{equation}
Thus, for example,
\[
 \mathfrak d(u)
 =
 u+\frac{u^2}{3}+\frac{11u^3}{72}
 +\frac{43u^4}{540}+\cdots.
\]

Define $b_j$ by the convergent expansion
\begin{equation}\label{eq:def-bj}
 -\frac12\log a
 -\frac12\log
 \left(\frac{\mathfrak d(as)}{as}\right)
 +R\bigl(1-\mathfrak d(as)\bigr)
 =
 \sum_{j\ge0}b_js^j.
\end{equation}
Equations \eqref{eq:dcoeff} and \eqref{eq:def-bj} give a finite
algorithm for computing $b_0,\ldots,b_N$ to any prescribed order.

We define the coefficients $C_r(\alpha)$ by the asymptotic expansion
\beq\label{eq:Cr-gamma}
 \frac{\Gamma(n-\alpha)}{\Gamma(n+1)}
 \sim
 n^{-\alpha-1}
 \sum_{r\ge0}\frac{C_r(\alpha)}{n^r},
 \qquad C_0(\alpha)=1.
\eeq
For more details and properties of the polynomials $C_r(\alpha)$ see \cite[App.\,C]{Cotti}.

\begin{theorem}[Full genus-one expansion]\label{thm:full1}
The radius of convergence of $G_1^{\chi}$ is
\(
 \rho=e-2,
\)
and $\rho$ is its unique dominant singularity. In a dented
neighborhood of $\rho$ one has the convergent Puiseux expansion
\begin{equation}\label{eq:puiseux1}
 G_1^{\chi}(t)
 =
 -\frac14\log\!\left(1-\frac{t}{\rho}\right)
 +
 \sum_{j\ge0}
 b_j\left(1-\frac{t}{\rho}\right)^{j/2}.
\end{equation}
In particular, no logarithmic terms other than the leading one occur.

Consequently,
\begin{equation}\label{eq:asym1}
 \chi(\Mb_{1,n})
 \sim
 n!\rho^{-n}
 \left\{
 \frac1{4n}
 +
 \sum_{m\ge0}\frac{d_m}{n^{m+3/2}}
 \right\},\end{equation}
where
\begin{equation}\label{eq:dm}
 d_m
 =
 \sum_{k=0}^m
 \frac{
 b_{2k+1}\,
 C_{m-k}\!\left(k+\frac12\right)
 }{
 \Gamma\!\left(-k-\frac12\right)
 }.
\end{equation}
\end{theorem}

Hence all coefficients in the full asymptotic expansion are
effectively computable from \eqref{eq:dcoeff}--\eqref{eq:dm}. By setting $\kappa:=\sqrt{\frac{e-2}{2\pi e}}$, the first coefficients are
\[
d_0
=
\kappa\,
\frac{2e^4+4e^3+9e^2+4e+1}{12},
\]
\[
d_1
=
-\kappa\,
\frac{
1470e^5-2070e^4-2685e^3-1785e^2+43e-101
}{
2160\,e
},
\]
and
\[
d_2
=
-\kappa\,
\frac{
20146e^6+9856e^5-90853e^4-40586e^3
+610e^2+962e-2439
}{
24192\,e^2
}.
\]
\begin{example}
We have
\[
\chi(\Mb_{1,50})
=
8326206901893890747516006661713788337956557083023720618082983829568912.
\]
Using only the terms through $n^{-7/2}$ in
\eqref{eq:asym1} gives
\[
\chi(\Mb_{1,50})
\approx
8.32628017946\times10^{69},
\]
with relative error approximately $8.8\times10^{-6}$.
For comparison, retaining terms only through $n^{-3/2}$ and
$n^{-5/2}$ gives relative errors approximately
$4.5\times10^{-3}$ and $2.9\times10^{-4}$, respectively.
\qetr\end{example}
\begin{proof}[Proof of Theorem \ref{thm:full1}]
Recall that
\[
 G_1^{\chi}(t)=A_1^{\chi}(\Psi(t)),
 \qquad
 \Psi(t)
 =
 2+W_{-1}\bigl(-(t+2)e^{-2}\bigr).
\]
The branch $W_{-1}$ is chosen so that $\Psi(0)=0$.

The closest singularity of this branch of the Lambert function occurs
when its argument reaches the branch point $-e^{-1}$. Since
\[
 -(t+2)e^{-2}=-e^{-1}
 \qquad\Longleftrightarrow\qquad
 t=e-2=\rho,
\]
the radius of convergence of $\Psi$ is $\rho$. The other singularity
of the branch lies at the origin of the Lambert variable and corresponds
to $t=-2$, which has larger modulus. Hence $\rho$ is the unique
singularity of $\Psi$ on its circle of convergence. The standard
analytic continuation of $W_{-1}$ away from its branch cut also gives
a dented neighborhood of $\rho$; see \cite{LambertRef}.

To determine the local form of the singularity, put
\[
 \delta:=1-\Psi(t).
\]
Since $t=\Phi(1-\delta)$ and
\[
 \Phi(z)=(2-z)e^z-2,
\]
we obtain
\[
 t+2=e(1+\delta)e^{-\delta}.
\]
Using $t=\rho(1-\zeta)$ and $\rho+2=e$, this becomes
\[
 \rho\zeta
 =
 e\bigl(1-(1+\delta)e^{-\delta}\bigr).
\]
With $a^2=2\rho/e$, this is equivalent to
\[
 2\bigl(1-(1+\delta)e^{-\delta}\bigr)
 =
 (a\sqrt\zeta)^2.
\]
By the definition of $\mathfrak d$,
\begin{equation}\label{eq:delta-d}
 \delta=\mathfrak d(a\sqrt\zeta).
\end{equation}
Since $\mathfrak d(u)=u+O(u^2)$, the quotient
$\mathfrak d(u)/u$ is analytic and nonzero at the origin. Therefore
\[
\begin{aligned}
 G_1^{\chi}(t)
 &=
 -\frac12\log\delta+R(1-\delta)\\
 &=
 -\frac14\log\zeta
 -\frac12\log a
 -\frac12
 \log\left(
 \frac{\mathfrak d(a\sqrt\zeta)}
      {a\sqrt\zeta}
 \right)
 +R\bigl(1-\mathfrak d(a\sqrt\zeta)\bigr).
\end{aligned}
\]
The last three terms are analytic functions of $\sqrt\zeta$.
Equation \eqref{eq:def-bj} therefore gives precisely
\eqref{eq:puiseux1}. It also shows directly that no further logarithms
can occur.

It remains to extract coefficients. The logarithmic term gives exactly
\[
 [t^n]\left(-\frac14\log(1-t/\rho)\right)
 =
 \frac{\rho^{-n}}{4n}.
\]
For $\alpha\notin\Z_{\ge0}$, the standard transfer formula
\cite[Chapter~VI]{FS} gives
\[
 [t^n](1-t/\rho)^\alpha
 =
 \rho^{-n}
 \frac{\Gamma(n-\alpha)}
      {\Gamma(-\alpha)\Gamma(n+1)}.
\]
The even powers $\zeta^k$ are locally analytic and therefore do not
contribute to the algebraic singular expansion. Only the odd powers
\(
 \zeta^{k+1/2}
\)
remain. Applying the preceding transfer formula to each of them and
expanding the gamma ratio in inverse powers of $n$ gives
\eqref{eq:asym1}--\eqref{eq:dm}.

Since $G_1^{\chi}$ is analytic in a dented disk at its unique dominant
singularity $\rho$, the transfer may be performed to any prescribed
finite order, which proves the complete asymptotic expansion.
\end{proof}

\section{Stabilization graphs, decorated expansions, and exact Euler contributions}
\label{sec:graph-formalism}

This section develops the constructions and computational tools
for determining, in higher genera $g\geq 2$, the all-light Euler series 
$A_g^\chi(z)$ and studying its asymptotics. Our strategy is to use the additivity
of Euler characteristic to decompose the all-light space into
geometrically manageable pieces whose contributions have a
tractable analytic form.
\vskip1mm

For $g\ge2$, forgetting all markings and stabilizing gives the
commutative diagram
\[
\begin{tikzpicture}[baseline=(C.base)]
\node (A) at (0,1.5) {$\Mb_{g,n}$};
\node (B) at (5,1.5) {$\Mb_{g,(0^+)^n}$};
\node (C) at (2.5,0) {$\Mb_g$};
\draw[->] (A) -- node[above] {$\rho_n$} (B);
\draw[->] (A) -- node[below left] {$\pi_n$} (C);
\draw[->] (B) -- node[below right] {$f_n$} (C);
\end{tikzpicture}
\]
where $\rho_n$ is the reduction morphism. The dual graph of the
resulting unmarked curve is its \emph{stabilization graph}. Writing
$\mathcal Z_\Gamma\subseteq\Mb_g$ for the stratum with dual graph $\Gamma$,
the inverse images $f_n^{-1}(\mathcal Z_\Gamma)$ stratify the all-light
space. Additivity therefore gives
\beq
\label{eq:decomposition_A_g_1}
A_g^\chi(z)=\sum_{[\Gamma]}A_{g,\Gamma}^\chi(z),
\qquad
A_{g,\Gamma}^\chi(z)
:=
\sum_{n\ge0}
\chi\!\left(f_n^{-1}(\mathcal Z_\Gamma)\right)\frac{z^n}{n!},
\eeq
where the outer sum runs over the finitely many isomorphism
classes of stabilization graphs of genus $g$, and Euler characteristics
are taken on coarse moduli spaces.
\vskip1mm
The advantage of this decomposition lies in the geometry above
each graph. In the ordinary space, forgetting the markings may
contract rational tails and branching trees. In the all-light
chamber, these are excluded (Lemma \ref{lem:all-light-stability}): 
the only components contracted are
rational chains along the edges of the graph. 
After fixing these
chains and assigning the light labels to the components, the
component data are parametrized by a product of moduli spaces.
Gluing and quotienting by the remaining graph symmetries then
give the corresponding piece of $f_n^{-1}(\mathcal Z_\Gamma)$.
\vskip1mm
Thus each graph contribution can be assembled from vertex and
edge data, with finite symmetries accounted for. The purpose
of this section is to make this construction precise and turn
it into a procedure for computing $A_{g,\Gamma}^\chi(z)$ and
determining its analytic structure.

\subsection{Graphs, stable graphs, and boundary strata}
\label{sec:stablegraphs}

A finite graph with legs
is a quadruple
\(
\Gamma=(V,H,\iota,\partial),
\)
where $V$ and $H$ are finite sets, $\iota:H\to H$ is an involution,
and $\partial:H\to V$ is the incidence map. The elements of $H$ will
be called \emph{half-edges}. The two-element orbits of $\iota$ are the
\emph{edges} of $\Gamma$, while the fixed points of $\iota$ are its
\emph{legs}. Thus every edge consists of two half-edges. A loop is an
edge whose two half-edges are incident to the same vertex. Multiple
edges and loops are allowed.
\vskip1mm
An automorphism of $\Gamma$ is a pair of bijections of $V$ and $H$
commuting with $\iota$ and $\partial$. In particular, an automorphism
which preserves an edge may either fix its two half-edges individually
or exchange them; in the latter case we shall say that it
\emph{reverses} the edge.

\vskip1mm

Let $\Gamma$ be a connected finite graph with legs. We shall use the
following notation throughout this section:
\begin{itemize}
\item $\operatorname{Vert}(\Gamma)$, $\operatorname{Edge}(\Gamma)$ and
$\operatorname{Leg}(\Gamma)$ denote respectively the sets of vertices,
edges and legs;

\item
\[
V(\Gamma):=\#\operatorname{Vert}(\Gamma),
\qquad
E(\Gamma):=\#\operatorname{Edge}(\Gamma);
\]

\item for a vertex $v$, $d_v$ denotes the number of half-edges
belonging to edges and incident to $v$, with a loop counted twice;

\item $m_v$ denotes the number of legs incident to $v$, and
\[
n_v:=d_v+m_v
\]
is the total number of half-edges and legs incident to $v$;

\item the integer
\[
b_1(\Gamma):=E(\Gamma)-V(\Gamma)+1
\]
is the first Betti number of $\Gamma$.
\end{itemize}

A \emph{stable graph of type $(g,n)$} is such a connected graph
$\Gamma$, together with a bijection
\(
\operatorname{Leg}(\Gamma)\xrightarrow{\sim}[n]
\)
and a genus decoration
\[
v\longmapsto h_v\in\mathbb Z_{\ge0},
\]
satisfying the following two conditions. First, setting
\[
h(\Gamma):=\sum_{v\in\operatorname{Vert}(\Gamma)}h_v,
\qquad
g(\Gamma):=h(\Gamma)+b_1(\Gamma),
\]
one requires
\[
g(\Gamma)=g.
\]
Second, the stability condition
\[
2h_v-2+n_v>0
\]
must hold at every vertex $v$.

The automorphism group $\Aut\Gamma$ consists of the graph
automorphisms preserving the genus decoration and the labels of the
legs.

\vskip1mm

For the ordinary moduli space $\Mb_{g,n}$, stable graphs describe the
boundary stratification; see, for instance,
\cite{GetzlerKapranov,ManinFrob,Hassett}. Fix a stable graph $\Gamma$ and
temporarily distinguish its vertices and half-edges. The corresponding
curves are parametrized by
\[
\prod_{v\in\operatorname{Vert}(\Gamma)}
\M_{h_v,n_v},
\]
the gluing being performed at the marked points indexed by the two
half-edges of each edge. Forgetting these auxiliary distinctions
amounts to quotienting by graph automorphisms; thus the corresponding
stack stratum is
\[
\left[
\prod_{v\in\operatorname{Vert}(\Gamma)}
\M_{h_v,n_v}
/\Aut\Gamma
\right].
\]
For a weighted pointed curve the dual graph has the same underlying
structure: the markings give labeled legs and the nodes give edges.
The graph alone, however, no longer records all the marking data,
since several weighted markings may coincide. One must additionally
record, at each vertex, which light labels determine the same marked
point, or equivalently the corresponding partition into collision
blocks. We shall not include this collision data in the graph
notation. In the minimal all-light chamber it will instead be absorbed
into the component moduli spaces introduced below.

\begin{figure}[ht]
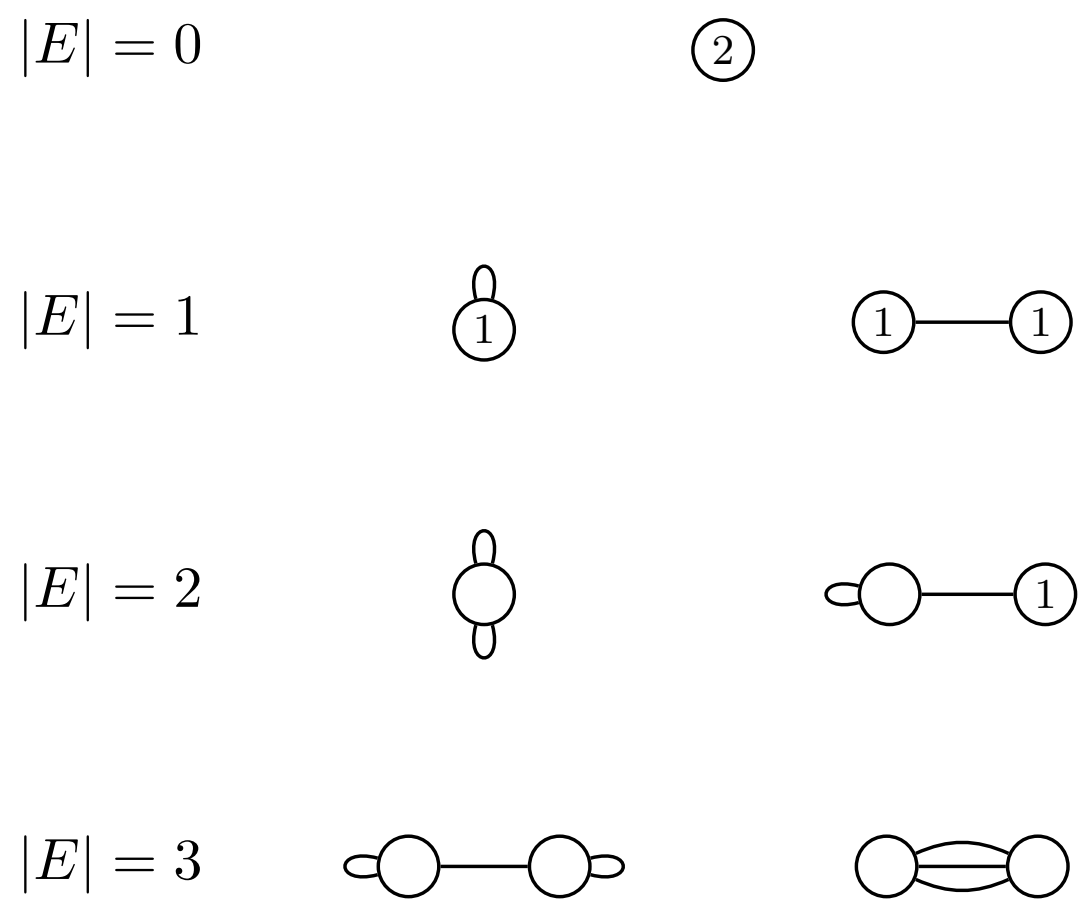

\centering
\setlength{\tabcolsep}{8pt}
\renewcommand{\arraystretch}{1.25}

\begin{tabular}{c cc}

$|E|=0$
&
\multicolumn{2}{c}{
\graphcell{}{
  \cv{a}{(0,0)}{$2$}
}}
\\[3mm]

$|E|=1$
&
\graphcell{}{
  \cv{a}{(0,-.05)}{$1$}
  \draw[graph edge] (a) edge[loop above,looseness=7] (a);
}
&
\graphcell{}{
  \cv{a}{(-.52,0)}{$1$}
  \cv{b}{(.52,0)}{$1$}
  \draw[graph edge] (a)--(b);
}
\\[3mm]

$|E|=2$
&
\graphcell{}{
  \cv{a}{(0,0)}{}
  \draw[graph edge] (a) edge[loop above,looseness=7] (a);
  \draw[graph edge] (a) edge[loop below,looseness=7] (a);
}
&
\graphcell{}{
  \cv{a}{(-.48,0)}{}
  \cv{b}{(.55,0)}{$1$}
  \draw[graph edge] (a)--(b);
  \draw[graph edge] (a) edge[loop left,looseness=7] (a);
}
\\[3mm]

$|E|=3$
&
\graphcell{}{
  \cv{a}{(-.5,0)}{}
  \cv{b}{(.5,0)}{}
  \draw[graph edge] (a)--(b);
  \draw[graph edge] (a) edge[loop left,looseness=7] (a);
  \draw[graph edge] (b) edge[loop right,looseness=7] (b);
}
&
\graphcell{}{
  \cv{a}{(-.5,0)}{}
  \cv{b}{(.5,0)}{}
  \tripleedge{a}{b}
}
\end{tabular}

\caption{The seven stabilization graphs of genus $2$, grouped according to the
number of edges (equivalently, the codimension of the corresponding
stratum of $\Mb_2$). An unlabeled vertex has genus $0$; a numerical
label denotes its genus weight.}
\label{fig:stable-graphs-g2}
\end{figure}

\subsection{Stabilization graphs}
\label{sec:stablegraphs}
Assume $g\ge2$. By Lemma~\ref{lem:all-light-stability},
forgetting all light markings destabilizes precisely the
two-valent rational components. Suppressing the corresponding
maximal chains in the dual graph gives the \emph{stabilization graph}:
a connected genus-decorated graph $\Gamma$ without legs whose
vertices satisfy
\begin{equation}
2h_v-2+d_v>0.
\label{eq:graphstability}
\end{equation}
Its automorphism group $\Aut\Gamma$ consists of graph
automorphisms preserving the genus decoration.

\vskip1mm

Conversely, at the level of the underlying graph, every weighted dual
graph having stabilization graph $\Gamma$ is obtained by subdividing the edges
(including loops)
of $\Gamma$ by chains of two-valent genus-zero vertices and then
distributing the light markings among the resulting vertices, with
every inserted two-valent vertex receiving at least one light
marking. The precise discrete data involved in this reconstruction
will be formalized in the next subsection.

\begin{remark}The suppressed vertices form paths: a cycle consisting
entirely of two-valent rational vertices would be the
whole connected graph and would have genus one.
The two ends of a path may meet the same surviving
vertex, in which case suppression produces a loop.
\qrem\end{remark}

\vskip1mm

With the notation introduced above,
\[
g
=
h(\Gamma)+E(\Gamma)-V(\Gamma)+1,
\]
and, since $\sum_vd_v=2E(\Gamma)$,
\begin{equation}
\sum_{v\in\operatorname{Vert}(\Gamma)}
\bigl(2h_v-2+d_v\bigr)
=
2g-2.
\label{eq:graphidentity}
\end{equation}

\begin{lemma}[The maximum number of graph edges]
\label{lem:maxedges}
For every stabilization graph of genus $g\ge2$,
\[
E(\Gamma)\le3g-3.
\]
Equality holds if and only if every vertex has genus zero and valence
three. Equivalently, equality holds precisely for connected rational
trivalent graphs of first Betti number $g$.
\end{lemma}

\begin{proof}
By \eqref{eq:graphstability}, every summand in
\eqref{eq:graphidentity} is a positive integer. Hence
\[
V(\Gamma)\le2g-2.
\]
Moreover,
\[
E(\Gamma)
=
g-1+V(\Gamma)-h(\Gamma)
\le
g-1+V(\Gamma)
\le
3g-3.
\]
If equality holds, then $V(\Gamma)=2g-2$ and $h(\Gamma)=0$.
Every summand in \eqref{eq:graphidentity} must therefore equal $1$.
Thus $h_v=0$ and $d_v=3$ for every vertex. The converse is immediate.
\end{proof}

\begin{remark}
A stabilization graph of genus $g$ is precisely a stable graph of
type $(g,0)$. Its isomorphism class therefore indexes a
dual-graph stratum of $\Mb_g$, of dimension
$3g-3-E(\Gamma)$. In particular, rational trivalent graphs
index the zero-dimensional strata.
\qrem\end{remark}

\begin{remark}
The number of rational trivalent stabilization graphs of genus $g$ is the
number of connected trivalent graphs on $2g-2$ vertices, with loops
and multiple edges allowed. These numbers are recorded in OEIS
sequence \href{https://oeis.org/A005967}{A005967}, with index $g-1$. Thus, for
$g=2,3,4,5,6,\ldots$, their numbers are
\(
2,\ 5,\ 17,\ 71,\ 388,\ldots.
\)
\qrem\end{remark}

Figures~\ref{fig:stable-graphs-g2} and~\ref{fig:stable-graphs-g3}
display all stabilization graphs in genera $2$ and $3$, respectively.

\begin{figure}[ht]
\centering
\setlength{\tabcolsep}{3pt}
\renewcommand{\arraystretch}{1.1}

\resizebox{\textwidth}{!}{\begin{tabular}{@{}c *{6}{c}@{}}

$|E|=0$
&
\graphcell{}{
  \cv{a}{(0,0)}{$3$}
}
& & & & &
\\[3mm]

$|E|=1$
&
\graphcell{}{
  \cv{a}{(0,-.05)}{$2$}
  \draw[graph edge] (a) edge[loop above,looseness=7] (a);
}
&
\graphcell{}{
  \cv{a}{(-.52,0)}{$1$}
  \cv{b}{(.52,0)}{$2$}
  \draw[graph edge] (a)--(b);
}
& & & &
\\[3mm]

$|E|=2$
&
\graphcell{}{
  \cv{a}{(0,0)}{$1$}
  \draw[graph edge] (a) edge[loop above,looseness=7] (a);
  \draw[graph edge] (a) edge[loop below,looseness=7] (a);
}
&
\graphcell{}{
  \cv{a}{(-.5,0)}{}
  \cv{b}{(.55,0)}{$2$}
  \draw[graph edge] (a)--(b);
  \draw[graph edge] (a) edge[loop left,looseness=7] (a);
}
&
\graphcell{}{
  \cv{a}{(-.5,0)}{$1$}
  \cv{b}{(.55,0)}{$1$}
  \draw[graph edge] (a)--(b);
  \draw[graph edge] (a) edge[loop left,looseness=7] (a);
}
&
\graphcell{}{
  \cv{a}{(-.5,0)}{$1$}
  \cv{b}{(.5,0)}{$1$}
  \doubleedge{a}{b}
}
&
\graphcell{}{
  \cv{a}{(-.72,0)}{$1$}
  \cv{b}{(0,0)}{$1$}
  \cv{c}{(.72,0)}{$1$}
  \draw[graph edge] (a)--(b)--(c);
}
&
\\[3mm]

\multirow{2}{*}{$|E|=3$}
&
\graphcell{}{
  \cv{a}{(0,0)}{}
  \draw[graph edge] (a) edge[loop above,looseness=7] (a);
  \draw[graph edge] (a) edge[loop left,looseness=7] (a);
  \draw[graph edge] (a) edge[loop right,looseness=7] (a);
}
&
\graphcell{}{
  \cv{a}{(-.48,0)}{}
  \cv{b}{(.58,0)}{$1$}
  \draw[graph edge] (a)--(b);
  \draw[graph edge] (a) edge[loop above,looseness=7] (a);
  \draw[graph edge] (a) edge[loop left,looseness=7] (a);
}
&
\graphcell{}{
  \cv{a}{(-.48,0)}{}
  \cv{b}{(.55,0)}{$1$}
  \doubleedge{a}{b}
  \draw[graph edge] (a) edge[loop left,looseness=7] (a);
}
&
\graphcell{}{
  \cv{a}{(-.5,0)}{}
  \cv{b}{(.5,0)}{$1$}
  \draw[graph edge] (a)--(b);
  \draw[graph edge] (a) edge[loop left,looseness=7] (a);
  \draw[graph edge] (b) edge[loop right,looseness=7] (b);
}
&
\graphcell{}{
  \cv{a}{(-.5,0)}{}
  \cv{b}{(.5,0)}{$1$}
  \tripleedge{a}{b}
}
&
\\[-1mm]

&
\graphcell{}{
  \cv{a}{(-.72,0)}{$1$}
  \cv{b}{(0,0)}{}
  \cv{c}{(.72,0)}{$1$}
  \draw[graph edge] (a)--(b)--(c);
  \draw[graph edge] (b) edge[loop above,looseness=7] (b);
}
&
\graphcell{}{
  \cv{a}{(-.72,0)}{}
  \cv{b}{(0,0)}{$1$}
  \cv{c}{(.72,0)}{$1$}
  \draw[graph edge] (a)--(b)--(c);
  \draw[graph edge] (a) edge[loop left,looseness=7] (a);
}
&
\graphcell{}{
  \cv{a}{(0,0)}{}
  \cv{b}{(-.72,0)}{$1$}
  \cv{c}{(.72,0)}{$1$}
  \doubleedge{a}{b}
  \draw[graph edge] (a)--(c);
}
&
\graphcell{}{
  \cv{a}{(0,-.08)}{}
  \cv{b}{(-.62,.38)}{$1$}
  \cv{c}{(.62,.38)}{$1$}
  \cv{d}{(0,-.62)}{$1$}
  \draw[graph edge] (a)--(b);
  \draw[graph edge] (a)--(c);
  \draw[graph edge] (a)--(d);
}
&
&
\\[4mm]

\multirow{2}{*}{$|E|=4$}
&
\graphcell{}{
  \cv{a}{(-.48,0)}{}
  \cv{b}{(.48,0)}{}
  \draw[graph edge] (a)--(b);
  \draw[graph edge] (a) edge[loop above,looseness=7] (a);
  \draw[graph edge] (a) edge[loop left,looseness=7] (a);
  \draw[graph edge] (b) edge[loop right,looseness=7] (b);
}
&
\graphcell{}{
  \cv{a}{(-.5,0)}{}
  \cv{b}{(.5,0)}{}
  \tripleedge{a}{b}
  \draw[graph edge] (a) edge[loop left,looseness=7] (a);
}
&
\graphcell{}{
  \cv{a}{(-.5,0)}{}
  \cv{b}{(.5,0)}{}
  \doubleedge{a}{b}
  \draw[graph edge] (a) edge[loop left,looseness=7] (a);
  \draw[graph edge] (b) edge[loop right,looseness=7] (b);
}
&
\graphcell{}{
  \cv{a}{(-.5,0)}{}
  \cv{b}{(.5,0)}{}
  \quadedge{a}{b}
}
&
\graphcell{}{
  \cv{a}{(-.7,0)}{}
  \cv{b}{(0,0)}{}
  \cv{c}{(.72,0)}{$1$}
  \doubleedge{a}{b}
  \draw[graph edge] (b)--(c);
  \draw[graph edge] (a) edge[loop left,looseness=7] (a);
}
&
\graphcell{}{
  \cv{a}{(0,0)}{}
  \cv{b}{(-.67,0)}{}
  \cv{c}{(.7,0)}{$1$}
  \draw[graph edge] (a)--(b);
  \draw[graph edge] (a)--(c);
  \draw[graph edge] (a) edge[loop above,looseness=7] (a);
  \draw[graph edge] (b) edge[loop left,looseness=7] (b);
}
\\[-1mm]

&
\graphcell{}{
  \cv{a}{(-.72,0)}{}
  \cv{b}{(0,0)}{}
  \cv{c}{(.72,0)}{$1$}
  \draw[graph edge] (a)--(b);
  \doubleedge{b}{c}
  \draw[graph edge] (a) edge[loop left,looseness=7] (a);
}
&
\graphcell{}{
  \cv{a}{(-.68,0)}{}
  \cv{b}{(.68,0)}{}
  \cv{c}{(0,0)}{$1$}
  \draw[graph edge] (a)--(c)--(b);
  \draw[graph edge] (a) edge[loop left,looseness=7] (a);
  \draw[graph edge] (b) edge[loop right,looseness=7] (b);
}
&
\graphcell{}{
  \cv{a}{(0,0)}{}
  \cv{b}{(-.68,0)}{}
  \cv{c}{(.72,0)}{$1$}
  \tripleedge{a}{b}
  \draw[graph edge] (a)--(c);
}
&
\graphcell{}{
  \cv{a}{(-.52,-.22)}{}
  \cv{b}{(.52,-.22)}{}
  \cv{c}{(0,.43)}{$1$}
  \doubleedge{a}{b}
  \draw[graph edge] (a)--(c)--(b);
}
&
\graphcell{}{
  \cv{a}{(-.68,.18)}{}
  \cv{b}{(0,0)}{}
  \cv{c}{(.65,.32)}{$1$}
  \cv{d}{(.05,-.58)}{$1$}
  \draw[graph edge] (a)--(b);
  \draw[graph edge] (b)--(c);
  \draw[graph edge] (b)--(d);
  \draw[graph edge] (a) edge[loop left,looseness=7] (a);
}
&
\graphcell{}{
  \cv{a}{(-.25,0)}{}
  \cv{b}{(.25,0)}{}
  \cv{c}{(-.82,0)}{$1$}
  \cv{d}{(.82,0)}{$1$}
  \doubleedge{a}{b}
  \draw[graph edge] (c)--(a);
  \draw[graph edge] (b)--(d);
}
\\[4mm]

\multirow{2}{*}{$|E|=5$}
&
\graphcell{}{
  \cv{a}{(-.7,0)}{}
  \cv{b}{(0,0)}{}
  \cv{c}{(.7,0)}{}
  \doubleedge{a}{b}
  \draw[graph edge] (b)--(c);
  \draw[graph edge] (a) edge[loop left,looseness=7] (a);
  \draw[graph edge] (c) edge[loop right,looseness=7] (c);
}
&
\graphcell{}{
  \cv{a}{(0,0)}{}
  \cv{b}{(-.62,-.12)}{}
  \cv{c}{(.62,-.12)}{}
  \draw[graph edge] (a)--(b);
  \draw[graph edge] (a)--(c);
  \draw[graph edge] (a) edge[loop above,looseness=7] (a);
  \draw[graph edge] (b) edge[loop left,looseness=7] (b);
  \draw[graph edge] (c) edge[loop right,looseness=7] (c);
}
&
\graphcell{}{
  \cv{a}{(-.75,0)}{}
  \cv{b}{(-.15,0)}{}
  \cv{c}{(.68,0)}{}
  \draw[graph edge] (a)--(b);
  \tripleedge{b}{c}
  \draw[graph edge] (a) edge[loop left,looseness=7] (a);
}
&
\graphcell{}{
  \cv{a}{(-.5,.22)}{}
  \cv{b}{(.5,.22)}{}
  \cv{c}{(0,-.22)}{}
  \doubleedge{a}{b}
  \draw[graph edge] (a)--(c);
  \draw[graph edge] (b)--(c);
  \draw[graph edge] (c) edge[loop below,looseness=7] (c);
}
&
&
\\[-1mm]

&
\graphcell{}{
  \cv{a}{(0,.4)}{}
  \cv{b}{(-.53,-.28)}{}
  \cv{c}{(.53,-.28)}{}
  \doubleedge{a}{b}
  \doubleedge{a}{c}
  \draw[graph edge] (b)--(c);
}
&
\graphcell{}{
  \cv{a}{(-.82,0)}{}
  \cv{b}{(-.25,0)}{}
  \cv{c}{(.32,0)}{}
  \cv{d}{(.88,0)}{$1$}
  \draw[graph edge] (a)--(b);
  \doubleedge{b}{c}
  \draw[graph edge] (c)--(d);
  \draw[graph edge] (a) edge[loop left,looseness=7] (a);
}
&
\graphcell{}{
  \cv{a}{(0,0)}{}
  \cv{b}{(-.62,.28)}{}
  \cv{c}{(.62,.28)}{}
  \cv{d}{(0,-.55)}{$1$}
  \draw[graph edge] (a)--(b);
  \draw[graph edge] (a)--(c);
  \draw[graph edge] (a)--(d);
  \draw[graph edge] (b) edge[loop above,looseness=7] (b);
  \draw[graph edge] (c) edge[loop above,looseness=7] (c);
}
&
\graphcell{}{
  \cv{a}{(-.45,.2)}{}
  \cv{b}{(.45,.2)}{}
  \cv{c}{(0,-.15)}{}
  \cv{d}{(0,-.65)}{$1$}
  \doubleedge{a}{b}
  \draw[graph edge] (a)--(c);
  \draw[graph edge] (b)--(c);
  \draw[graph edge] (c)--(d);
}
&
&
\\[4mm]

$|E|=6$
&
\graphcell{}{
  \cv{a}{(-.82,0)}{}
  \cv{b}{(-.28,0)}{}
  \cv{c}{(.28,0)}{}
  \cv{d}{(.82,0)}{}
  \draw[graph edge] (a)--(b);
  \doubleedge{b}{c}
  \draw[graph edge] (c)--(d);
  \draw[graph edge] (a) edge[loop left,looseness=7] (a);
  \draw[graph edge] (d) edge[loop right,looseness=7] (d);
}
&
\graphcell{}{
  \cv{a}{(0,0)}{}
  \cv{b}{(-.62,.38)}{}
  \cv{c}{(.62,.38)}{}
  \cv{d}{(0,-.55)}{}
  \draw[graph edge] (a)--(b);
  \draw[graph edge] (a)--(c);
  \draw[graph edge] (a)--(d);
  \draw[graph edge] (b) edge[loop left,looseness=7] (b);
  \draw[graph edge] (c) edge[loop right,looseness=7] (c);
  \draw[graph edge] (d) edge[loop below,looseness=7] (d);
}
&
\graphcell{}{
  \cv{a}{(-.72,.22)}{}
  \cv{b}{(-.1,0)}{}
  \cv{c}{(.58,.35)}{}
  \cv{d}{(.58,-.4)}{}
  \draw[graph edge] (a)--(b);
  \draw[graph edge] (b)--(c);
  \draw[graph edge] (b)--(d);
  \doubleedge{c}{d}
  \draw[graph edge] (a) edge[loop left,looseness=7] (a);
}
&
\graphcell{}{
  \cv{a}{(-.48,.35)}{}
  \cv{b}{(.48,.35)}{}
  \cv{c}{(-.48,-.35)}{}
  \cv{d}{(.48,-.35)}{}
  \doubleedge{a}{b}
  \draw[graph edge] (a)--(c);
  \draw[graph edge] (b)--(d);
  \doubleedge{c}{d}
}
&
\graphcell{}{
  \cv{a}{(-.48,.38)}{}
  \cv{b}{(.48,.38)}{}
  \cv{c}{(-.48,-.38)}{}
  \cv{d}{(.48,-.38)}{}
  \draw[graph edge] (a)--(b);
  \draw[graph edge] (a)--(c);
  \draw[graph edge] (a)--(d);
  \draw[graph edge] (b)--(c);
  \draw[graph edge] (b)--(d);
  \draw[graph edge] (c)--(d);
}
&
\\

\end{tabular}}

\caption{The forty-two stabilization graphs of genus $3$, grouped according to the
number of edges (equivalently, the codimension of the corresponding
stratum of $\Mb_3$). An unlabeled vertex has genus $0$; a numerical
label denotes its genus weight.}
\label{fig:stable-graphs-g3}

\end{figure}

\subsection{Decorated expansions of a stabilization graph}
\label{sec:decorated-expansions}

Fix a stabilization graph $\Gamma$ and an integer $n\ge0$. We now make precise
the discrete data needed to recover the weighted dual graphs lying
above $\Gamma$.

By a \emph{discrete decoration} of $\Gamma$ we mean a triple
\[
\delta=(\widetilde\Gamma,c,\ell)
\]
with the following properties:
\begin{enumerate}
\item $\widetilde\Gamma$ is obtained by subdividing
the edges of $\Gamma$ with two-valent genus-zero
vertices. 
\item Let $\widetilde\Gamma^{\mathrm{st}}$ denote
the stabilization of $\widetilde\Gamma$, obtained
by suppressing its two-valent genus-zero vertices,
and let
\[
c:\widetilde\Gamma^{\mathrm{st}}\xrightarrow{\sim}\Gamma
\]
be a genus-preserving graph isomorphism.

\item the function
\[
\ell:[n]\longrightarrow\operatorname{Vert}(\widetilde\Gamma)
\]
assigns each light label to a vertex of $\widetilde\Gamma$, and every
inserted two-valent vertex receives at least one label.
\end{enumerate}

Graph vertices are allowed to receive no light labels.
We use $c$ to identify the vertices of $\Gamma$
with the surviving vertices of $\widetilde\Gamma$.
Likewise, each half-edge of $\Gamma$ indexes the
corresponding end half-edge of its subdividing chain,
also for loops.

An isomorphism over $\Gamma$ between decorations
$\delta=(\widetilde\Gamma,c,\ell)$ and
$\delta'=(\widetilde\Gamma',c',\ell')$
is a genus-preserving graph isomorphism
$\varphi:\widetilde\Gamma\to\widetilde\Gamma'$
such that the following diagrams commute:
\[
\begin{CD}
\widetilde\Gamma^{\mathrm{st}}
@>{\varphi^{\mathrm{st}}}>>
\widetilde\Gamma'^{\mathrm{st}}\\
@V{c}VV @VV{c'}V\\
\Gamma @>{\mathrm{id}}>> \Gamma
\end{CD}
\qquad
\begin{CD}
\operatorname{Vert}(\widetilde\Gamma)
@>{\varphi_V}>>
\operatorname{Vert}(\widetilde\Gamma')\\
@A{\ell}AA @AA{\ell'}A\\
[n] @>{\mathrm{id}}>> [n].
\end{CD}
\]Here $\varphi^{\mathrm{st}}$ is the graph isomorphism
induced by suppression.
\vskip1mm
Let $\mathscr D_n(\Gamma)$ be the set of decorations
up to isomorphism over $\Gamma$. This set is finite,
since every inserted vertex carries a light label,
so at most $n$ vertices can be inserted.

\begin{figure}[ht]
\centering

\begin{tikzpicture}[
  dec graph/.style={
  circle,draw,line width=.65pt,fill=white,
  minimum size=4mm,inner sep=0pt,font=\scriptsize
},
  dec inserted/.style={
    circle,fill=black,
    minimum size=2.4mm,inner sep=0pt
  },
  dec edge/.style={line width=.65pt},
  dec marking/.style={densely dashed,line width=.55pt},
  contraction/.style={->,>=stealth,line width=.6pt},
  font=\small
]

\begin{scope}[xshift=0cm]

\node at (0,1.8) {$\widetilde\Gamma$};

\node[dec graph] (v1) at (-2,0) {$1$};
\node[dec graph] (v2) at ( 2,0) {$1$};

\node[below=4pt] at (v1) {$v_1$};
\node[below=4pt] at (v2) {$v_2$};

\node[dec inserted] (u1) at (-.7,0) {};
\node[dec inserted] (u2) at ( .7,0) {};

\draw[dec edge] (v1)--(u1)--(u2)--(v2);
\node at (0,.23) {$e_1$};

\draw[dec edge]
  (v1) .. controls (-1.1,1.2) and (1.1,1.2) .. (v2);
\node at (0,1.15) {$e_2$};

\draw[dec marking] (v1)--(-2.65,.65)
  node[above left,inner sep=1pt] {$1$};
\draw[dec marking] (v1)--(-2.65,-.55)
  node[below left,inner sep=1pt] {$5$};

\node[above=3pt] at (u1) {$u_1$};
\node[above=3pt] at (u2) {$u_2$};

\draw[dec marking] (u1)--(-.7,-.42)
  node[below,inner sep=1pt] {$2$};

\draw[dec marking] (u2)--(.48,-.38)
  node[below left,inner sep=1pt] {$3$};
\draw[dec marking] (u2)--(.92,-.38)
  node[below right,inner sep=1pt] {$4$};
\end{scope}

\draw[contraction] (2.65,0)--(4.15,0)
  node[midway,above=3pt] {suppress, then $c$};

\begin{scope}[xshift=6.8cm]

\node at (0,1.8) {$\Gamma$};

\node[dec graph] (w1) at (-1.55,0) {$1$};
\node[dec graph] (w2) at ( 1.55,0) {$1$};

\node[below=4pt] at (w1) {$v_1$};
\node[below=4pt] at (w2) {$v_2$};

\draw[dec edge] (w1)--(w2);
\node at (0,.23) {$e_1$};

\draw[dec edge]
  (w1) .. controls (-.8,1.1) and (.8,1.1) .. (w2);
\node at (0,1.05) {$e_2$};

\end{scope}

\end{tikzpicture}

\caption{A decorated expansion of a stabilization graph. The open circles are
graph vertices; their numerical labels denote the genus decoration.
The filled circles are inserted two-valent genus-zero vertices, and
the dashed segments indicate the light labels assigned by $\ell$.
Here
\[
B_{v_1}=\{1,5\},\qquad
B_{u_1}=\{2\},\qquad
B_{u_2}=\{3,4\},\qquad
B_{v_2}=\varnothing.
\]
The edge $e_1$ is subdivided by $u_1,u_2$, while $e_2$ is left
unchanged. Suppressing the inserted vertices gives
$\widetilde\Gamma^{\mathrm{st}}$, which $c$ identifies
with the fixed graph $\Gamma$.}
\label{fig:decorated-expansion}

\end{figure}

The decoration records the expanded graph and the
component carrying each light label. Positions and
collisions of light markings are included in the
component moduli spaces introduced below.

\vskip1mm
The group $\Aut\Gamma$ acts on $\mathscr D_n(\Gamma)$
by
\[
\sigma\cdot
(\widetilde\Gamma,c,\ell)
=
(\widetilde\Gamma,\sigma\circ c,\ell).
\]For $\delta\in\mathscr D_n(\Gamma)$ we write
\[
G_\delta
:=
\operatorname{Stab}_{\Aut\Gamma}(\delta).
\]

\begin{remark}
Suppression induces an injective homomorphism
\[
\Aut(\widetilde\Gamma,\ell)
\longrightarrow \Aut\Gamma,
\qquad
\varphi\longmapsto
c\circ\varphi^{\mathrm{st}}\circ c^{-1},
\]
whose image is $G_\delta$.
Indeed, the image consists precisely of the graph
symmetries preserving the decoration, and an
automorphism in the kernel fixes both ends of every
subdividing chain separately, hence is the identity.
\qrem
\end{remark}

The group $G_\delta$ is a combinatorial stabilizer,
not the automorphism group of a particular pointed curve.

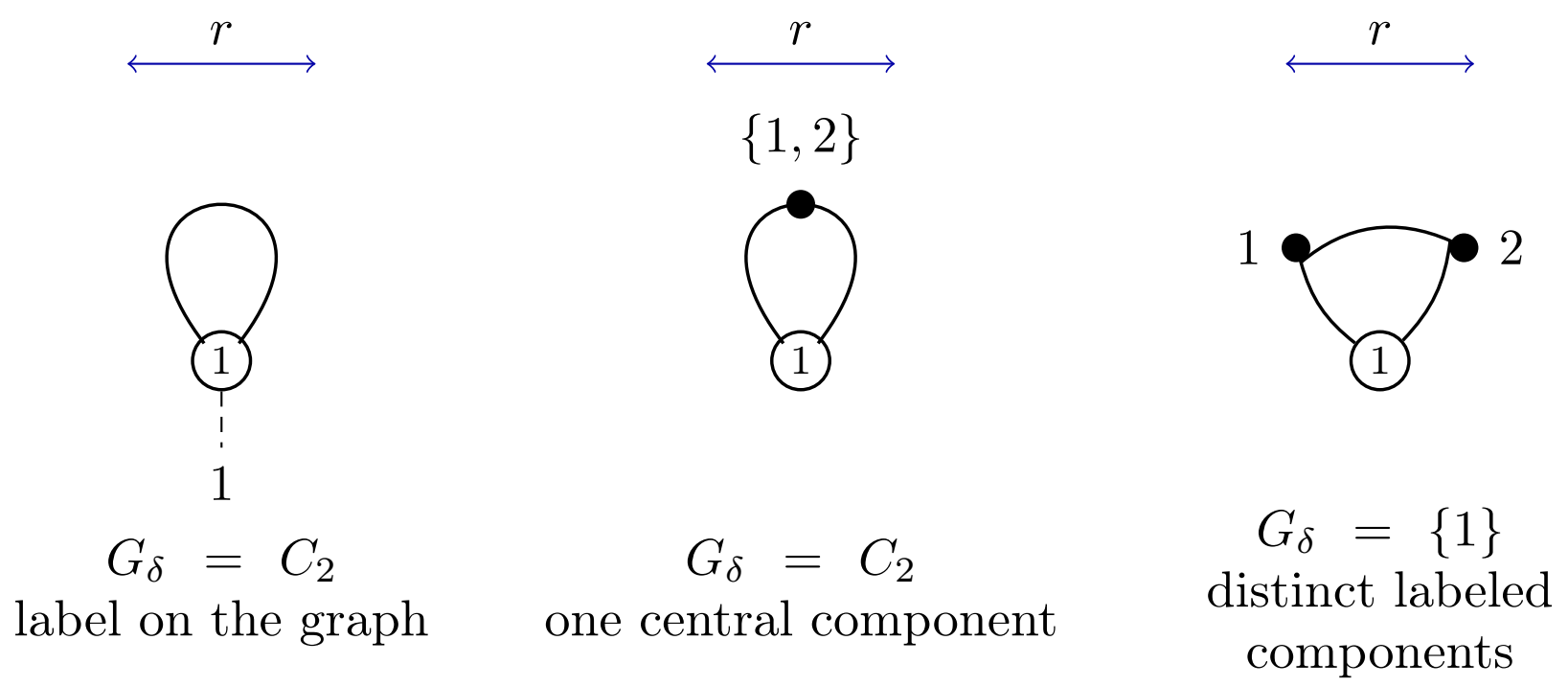
\begin{figure}[ht]
\centering
\begin{tikzpicture}[every node/.style={font=\small}]
\foreach \x/\name in {-4/a,0/b}{
 \cv{\name}{(\x,0)}{1}
 \draw[graph edge] (\x-.12,.12) .. controls (\x-1.1,1.4)
 and (\x+1.1,1.4) .. (\x+.12,.12);
}
\draw[densely dashed] (a)--(-4,-.6) node[below] {$1$};
\node[circle,fill=black,inner sep=2pt] (mid) at (0,1.08) {};
\node[above=4pt] at (mid) {$\{1,2\}$};
\cv{c}{(4,0)}{1}
\node[circle,fill=black,inner sep=2pt] (left) at (3.42,.78) {};
\node[circle,fill=black,inner sep=2pt] (right) at (4.58,.78) {};
\draw[graph edge] (c) to[bend left=18] (left)
 to[bend left=32] (right) to[bend left=18] (c);
\node[left=3pt] at (left) {$1$}; \node[right=3pt] at (right) {$2$};
\foreach \x in {-4,0,4}{
\draw[<->,blue!65!black] (\x-.65,2.05)--(\x+.65,2.05);
\node[above] at (\x,2.05) {$r$};}
\node[text width=3.7cm,align=center] at (-4,-1.6)
 {$G_\delta=C_2$\\label on the graph};
\node[text width=3.7cm,align=center] at (0,-1.6)
 {$G_\delta=C_2$\\one central component};
\node[text width=3.7cm,align=center] at (4,-1.6)
 {$G_\delta=\{1\}$\\distinct labeled components};
\end{tikzpicture}
\caption{Reversal of a loop on a genus-one graph vertex. Filled dots
are inserted rational vertices. In the right panel, reversal exchanges
the components carrying labels $1$ and $2$, so it does not preserve the
decoration. The double arrows indicate the reversal symmetry.}
\label{fig:loop-stabilizers}
\end{figure}

\begin{example}[A loop reversal and its decoration stabilizers]
\label{ex:loop-decoration-stabilizers}
Let $\Gamma$ have one genus-one vertex $v$ and a single loop. It is a
stabilization graph of total genus two, with $\Aut\Gamma=\langle r\rangle
\simeq\mathbb Z/2\mathbb Z$, where $r$ exchanges the half-edges of
the loop. Figure~\ref{fig:loop-stabilizers} shows three decorations.
If label $1$ lies on $v$ and the loop is not subdivided, then
$G_\delta=\Aut\Gamma$. The same holds if the loop is subdivided
by one vertex $u$ carrying both labels $1,2$: reversal fixes $u$.
If instead two successive inserted vertices carry labels $1$ and $2$
separately, reversal interchanges their order. It produces a different
decoration, since labels must be preserved individually, and
$G_\delta$ is trivial.
\qetr
\end{example}

We next describe the geometric stratum associated with a decoration.
\vskip1mm

For a finite set $B$ and a finite set $S$ of cardinality $d$ satisfying
$2h-2+d>0$, let
\(
\mathscr V_{h,S,B}
\)
denote the moduli stack\footnote{Equivalently, $\mathscr V_{h,S,B}$ is the smooth-curve
locus of $\overline{\eu M}_{h,(1^S,\epsilon^B)}$, with
$\epsilon>0$ sufficiently small and $|B|\epsilon<1$.}
of smooth genus-$h$ curves endowed with
pairwise distinct heavy markings $(x_s)_{s\in S}$ and light markings
$(p_b)_{b\in B}$, where the light markings may collide with one
another but are disjoint from all heavy markings. Let
$\mathcal V_{h,S,B}$ be its coarse moduli space. Up to relabeling,
these spaces depend only on $d=\#S$ and $n=\#B$; when only the
cardinalities matter we shall therefore write
\[
\mathscr V_{h,d,n},
\qquad
\mathcal V_{h,d,n}.
\]

Let now $\delta=(\widetilde\Gamma,c,\ell)$ and, for every
$w\in\operatorname{Vert}(\widetilde\Gamma)$, put
\(
B_w:=\ell^{-1}(w).
\)
For a graph vertex $v\in\operatorname{Vert}(\Gamma)$, let $S_v$ be the
set of half-edges of $\Gamma$ incident to $v$. These half-edges index
the branches which will be glued to the adjacent components.

If $u$ is an inserted two-valent vertex, then $B_u\neq\varnothing$.
The corresponding rational component has two distinguished nodal
points and light markings indexed by $B_u$. After normalizing the two
nodal points to $0$ and $\infty$, its moduli space is
\[
\eu M_{0,2|B_u}
\cong
(\mathbb G_m)^{|B_u|-1}.
\]

The \emph{parameter stack} associated with $\delta$ is therefore
the product stack
\begin{equation}
\mathscr X_\delta
:=
\prod_{v\in\operatorname{Vert}(\Gamma)}
\mathscr V_{h_v,S_v,B_v}
\times
\prod_{u\in
\operatorname{Vert}(\widetilde\Gamma)\setminus
\operatorname{Vert}(\Gamma)}
\eu M_{0,2|B_u}.
\label{eq:decoration-parameter-stack}
\end{equation}
The product is indexed by a fixed representative of $\Gamma$:
vertices distinguish the factors and half-edges distinguish
their heavy markings, including the two branches of a loop.
These indices and all light labels are preserved within each
factor; graph symmetries are imposed by the action of
$G_\delta$ below.
Let $X_\delta$ denote the coarse moduli space of
$\mathscr X_\delta$; concretely,
\[
X_\delta
\cong
\prod_{v\in\operatorname{Vert}(\Gamma)}
\mathcal V_{h_v,S_v,B_v}
\times
\prod_{u\in
\operatorname{Vert}(\widetilde\Gamma)\setminus
\operatorname{Vert}(\Gamma)}
(\mathbb G_m)^{|B_u|-1}.
\]

Via $G_\delta\cong\Aut(\widetilde\Gamma,\ell)$,
graph symmetries act on $\mathscr X_\delta$ by
permuting the component factors and relabeling
their attaching points, while fixing every light
label. The coarse moduli map is equivariant,
and
\[
\operatorname{coarse}
\bigl([\mathscr X_\delta/G_\delta]\bigr)
\cong X_\delta/G_\delta;
\]
see \cite[Theorem~3.2 and Corollary~3.3]{AOV08}.
Internal component automorphisms are retained
in $\mathscr X_\delta$; graph symmetries are
accounted for by $G_\delta$.
\vskip1mm
The following examples illustrate the $G_\delta$-action on inserted
components and on graph-vertex moduli.

\begin{example} Consider the middle decoration in
Example~\ref{ex:loop-decoration-stabilizers}. The graph
consists of one genus-one vertex with a loop. The loop
is subdivided by a single rational vertex carrying the
light labels $1,2$, while the elliptic vertex carries
no light markings. Loop reversal preserves this decoration,
so
\(
G_\delta=\Aut\Gamma=\langle r\rangle\cong C_2.
\)
The parameter stack and its coarse moduli space are
\[
\mathscr X_\delta
\cong \eu M_{1,2}\times\mathbb G_m,
\qquad
X_\delta
\cong M_{1,2}\times\mathbb G_m.
\]
The first factor parametrizes a smooth elliptic curve
with two ordered attaching points. The second parametrizes
the inserted rational component: put its attaching points
at $0,\infty$, its first light marking at $1$, and write
$t\in\mathbb G_m$ for its second light marking.
The value $t=1$ is allowed, since the light markings
may coincide.

The generator $r$ exchanges the two attaching points
on each component. On the rational component, the change
of coordinate $x\mapsto x^{-1}$ restores the chosen
normalization, sending $t$ to $t^{-1}$. Thus its action
on the coarse parameter space is
\[
r\cdot\bigl([E;p,q],t\bigr)
=
\bigl([E;q,p],t^{-1}\bigr).
\]
Here the light labels are not exchanged: $t^{-1}$ is
the new coordinate of the marking still labeled $2$.

In this example, $[E;p,q]=[E;q,p]$ in $M_{1,2}$:
the elliptic involution $x\mapsto p+q-x$ exchanges
$p$ and $q$. Hence $r$ acts trivially on the coarse
elliptic factor, and the $r$-fixed sublocus in $X_\delta$ is
\(
X_\delta^r=M_{1,2}\times\{1,-1\}.
\)
Thus $G_\delta$ preserves the discrete decoration,
but fixes only those parameter points for which
$t=\pm1$.
\qetr
\end{example}

\begin{example}For a nontrivial action on a vertex factor, consider
the graph consisting of one rational vertex with two loops,
without subdivisions or light markings. Since the decoration
adds no further data, every graph automorphism preserves it:
\[
G_\delta=\Aut\Gamma
\cong(C_2\times C_2)\rtimes C_2,
\qquad |G_\delta|=8.
\]
The first two factors reverse the loops independently,
while the last exchanges them. The parameter stack is
\(
\mathscr X_\delta=\mathscr V_{0,4,0}
\cong\eu M_{0,4}.
\)
Write its four attaching points as
$(p_1,p_2,p_3,p_4)=(0,\infty,1,\lambda)$, with
$\lambda\in\mathbb C\setminus\{0,1\}$, pairing
$p_1,p_2$ to form the first loop and $p_3,p_4$
to form the second.

Let $r\in G_\delta$ reverse only the first loop.
It exchanges $p_1,p_2$ and fixes $p_3,p_4$.
Renormalizing by $x\mapsto x^{-1}$ gives
\(
r\colon\lambda\longmapsto\lambda^{-1}.
\)
Thus the action on the vertex moduli stack is nontrivial;
the fixed locus of $r$ is the single configuration
$\lambda=-1$.
\qetr
\end{example}

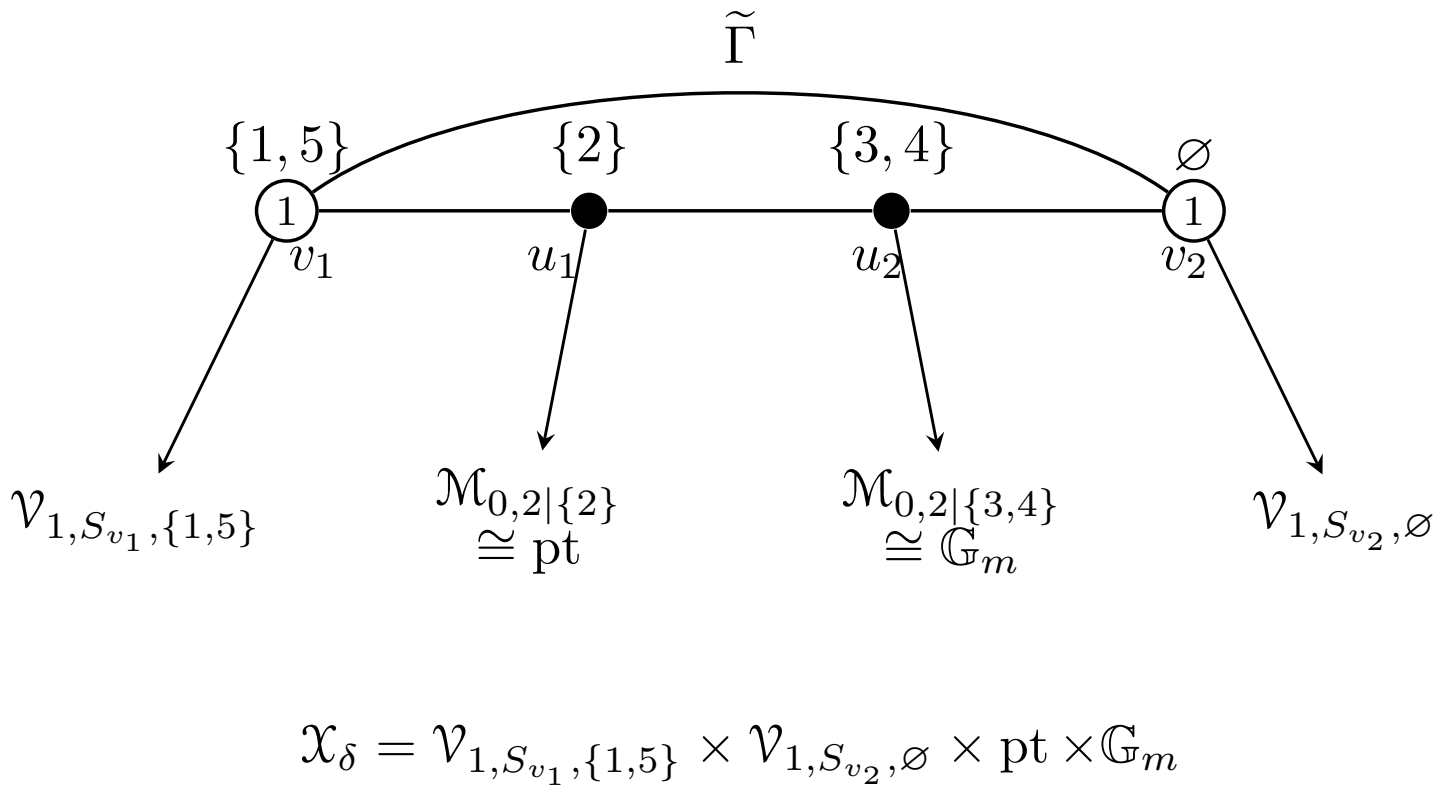
\begin{figure}[ht]
\centering

\begin{tikzpicture}[
  par graph/.style={
  circle,draw,line width=.65pt,fill=white,
  minimum size=4mm,inner sep=0pt,font=\scriptsize
},
  par inserted/.style={
    circle,fill=black,
    minimum size=2.4mm,inner sep=0pt
  },
  par edge/.style={line width=.65pt},
  par arrow/.style={->,>=stealth,line width=.55pt},
  factor/.style={
    align=center,
    inner sep=3pt,
    font=\small
  },
  font=\small
]

\node at (0,2.15) {$\widetilde\Gamma$};

\node[par graph] (v1) at (-3,1) {$1$};
\node[par inserted] (u1) at (-1,1) {};
\node[par inserted] (u2) at ( 1,1) {};
\node[par graph] (v2) at ( 3,1) {$1$};

\draw[par edge] (v1)--(u1)--(u2)--(v2);
\draw[par edge]
  (v1) .. controls (-1.6,2) and (1.6,2) .. (v2);

\node[above=3pt] at (v1) {$\{1,5\}$};
\node[above=3pt] at (u1) {$\{2\}$};
\node[above=3pt] at (u2) {$\{3,4\}$};
\node[above=3pt] at (v2) {$\varnothing$};

\node[below=3pt] at (v1) {$\quad v_1$};
\node[below=3pt] at (u1) {$u_1\quad\,\,$};
\node[below=3pt] at (u2) {$u_2\,\,\,$};
\node[below=3pt] at (v2) {$v_2\,\,$};

\node[factor] (f1) at (-4,-1.05)
  {$\mathscr V_{1,S_{v_1},\{1,5\}}$};

\node[factor] (f2) at (-1.4,-1.05)
  {$\eu M_{0,2|\{2\}}$\\[-1pt]
   $\cong \operatorname{pt}$};

\node[factor] (f3) at (1.4,-1.05)
  {$\eu M_{0,2|\{3,4\}}$\\[-1pt]
   $\cong \mathbb G_m$};

\node[factor] (f4) at (4,-1.05)
  {$\mathscr V_{1,S_{v_2},\varnothing}$};

\draw[par arrow] (v1) -- (f1);
\draw[par arrow] (u1) -- (f2);
\draw[par arrow] (u2) -- (f3);
\draw[par arrow] (v2) -- (f4);

\node[align=center] at (0,-2.6) {
$\displaystyle
\mathscr X_\delta
=
\mathscr V_{1,S_{v_1},\{1,5\}}
\times
\mathscr V_{1,S_{v_2},\varnothing}
\times
\operatorname{pt}
\times
\mathbb G_m
$};

\end{tikzpicture}

\caption{The parameter stack associated with the decoration in
Figure~\ref{fig:decorated-expansion}. Each graph vertex contributes
the corresponding vertex moduli space
$\mathscr V_{h_v,S_v,B_v}$, whereas an inserted two-valent vertex
$u$ contributes
$\eu M_{0,2|B_u}\cong(\mathbb G_m)^{|B_u|-1}$.
The parameter stack $\mathscr X_\delta$ is the product of these
component moduli stacks before quotienting by graph automorphisms.}
\label{fig:decoration-parameter-space}

\end{figure}

Thus, for a fixed stabilization graph, the geometry separates into two levels:
the finite combinatorics of the decorations
$\delta\in\mathscr D_n(\Gamma)$ and the continuous component moduli
encoded by $\mathscr X_\delta$. The Euler-characteristic analysis
below will handle the finite quotients by the Burnside--Lefschetz
formula. 
\subsection{The contribution of a fixed graph}
\label{sec:graphcontribution}

We shall use the following Euler-characteristic analogue of Burnside's
lemma.

\begin{lemma}[Burnside--Lefschetz formula]
\label{lem:burnside}
Let a finite group $G$ act algebraically on a complex algebraic variety
$X$. Then
\begin{equation}
\chi(X/G)
=
\frac1{|G|}
\sum_{\sigma\in G}\chi(X^\sigma),
\label{eq:burnside}
\end{equation}
where $X^\sigma$ denotes the fixed-point locus of $\sigma$.
\end{lemma}

This result goes back at least to Conner
\cite[Theorem~2]{Conner}, in the topological and differentiable setting of finite group
actions. The present algebraic situation reduces to that framework
upon passing to the underlying complex analytic spaces. For
completeness, we include a short proof.

\begin{proof}
Let $Y:=X/G$ and let $\pi:X\to Y$ be the quotient map. Since $\pi$
is finite, hence proper,
\(
H_c^\bullet(X;\mathbb Q)
\simeq
H_c^\bullet(Y;\pi_*\mathbb Q_X).
\)
Moreover,
\(
(\pi_*\mathbb Q_X)^G\simeq\mathbb Q_Y.
\)
Indeed, the stalk of $\pi_*\mathbb Q_X$ at $y\in Y$ is the space of
$\mathbb Q$-valued functions on the finite orbit $\pi^{-1}(y)$, whose
$G$-invariant elements are precisely the constant functions. Since
taking $G$-invariants is exact over $\mathbb Q$, it commutes with
cohomology. Therefore
\(
H_c^\bullet(X/G;\mathbb Q)
\simeq
H_c^\bullet(X;\mathbb Q)^G.
\)

For each $i$, consider the Reynolds operator
\[
\Pi_G
:=
\frac1{|G|}
\sum_{\sigma\in G}\sigma
\in
\operatorname{End}_{\mathbb Q}\bigl(H_c^i(X;\mathbb Q)\bigr).
\]
It satisfies $\Pi_G^2=\Pi_G$ and
\(
\operatorname{Im}\Pi_G=H_c^i(X;\mathbb Q)^G,
\)
so it is the projector onto the invariant subspace. Hence
\[
\dim H_c^i(X;\mathbb Q)^G
=
\operatorname{Tr}(\Pi_G)
=
\frac1{|G|}
\sum_{\sigma\in G}
\operatorname{Tr}
\bigl(\sigma\mid H_c^i(X;\mathbb Q)\bigr).
\]

Taking alternating sums yields
\[
\chi(X/G)
=
\frac1{|G|}
\sum_{\sigma\in G}L_c(\sigma;X),\qquad \text{where }\qquad L_c(\sigma;X)
:=
\sum_i(-1)^i
\operatorname{Tr}
\bigl(\sigma\mid H_c^i(X;\mathbb Q)\bigr).
\] For a finite-order automorphism, the compactly supported
Lefschetz fixed-point formula gives
\(
L_c(\sigma;X)=\chi_c(X^\sigma).
\)
This remains valid when $X$ is singular: choose a finite
$\langle\sigma\rangle$-invariant stratification of $X$ into
smooth locally closed subvarieties, apply the formula on
each stratum, and use the additivity of both sides.
Since $\chi_c(X^\sigma)=\chi(X^\sigma)$ for complex
algebraic varieties, the claim follows.

For background on finite group actions and fixed-point
formulas, see \cite[Chapter~III]{Bredon}; for the classical
Lefschetz fixed-point theorem on manifolds, see also
\cite[\S15.3]{DFN}.
\end{proof}

\begin{remark}
The factor $1/|G|$ in \eqref{eq:burnside} is a consequence of the
averaging projector and the Lefschetz trace formula. In particular,
no analogous identity is being assumed for classes in
$K_0(\Var_{\C})$ or for an arbitrary motivic measure.
\qrem\end{remark}

We now fix a stabilization graph $\Gamma$ of genus $g\ge2$. For each $n\ge0$,
consider the finite disjoint union
\[
X_{\Gamma,n}
:=
\coprod_{\delta\in\mathscr D_n(\Gamma)}X_\delta .
\]
The group $\Aut\Gamma$ acts naturally on $X_{\Gamma,n}$, and the
quotient decomposes as
\[
X_{\Gamma,n}/\Aut\Gamma
\cong
\coprod_{[\delta]\in\mathscr D_n(\Gamma)/\Aut\Gamma}
X_\delta/G_\delta.
\]

Recall that $f_n$ forgets all light markings and stabilizes
the curve, while $\mathcal Z_\Gamma\subseteq\Mb_g$ is the stratum
of unmarked curves with dual graph $\Gamma$.
Thus $f_n^{-1}(\mathcal Z_\Gamma)$ parametrizes all-light curves
whose stabilization graph is isomorphic to $\Gamma$.

\begin{lemma}
\label{lem:graph-stratum-euler}
For every stabilization graph $\Gamma$ of genus $g\ge2$ and every
$n\ge0$,
\[
\chi\!\left(f_n^{-1}(\mathcal Z_\Gamma)\right)
=
\chi\!\left(X_{\Gamma,n}/\Aut\Gamma\right)
=
\sum_{[\delta]\in\mathscr D_n(\Gamma)/\Aut\Gamma}
\chi(X_\delta/G_\delta).
\]
\end{lemma}

\begin{proof}
Let $\mathscr T_{[\delta]}\subseteq f_n^{-1}(\mathcal Z_\Gamma)$
be the reduced locally closed substack parametrizing
curves with decoration orbit $[\delta]$, allowing
collisions among light markings.
These substacks give a finite decomposition\footnote{We use $\bigsqcup$ for a decomposition into locally
closed pieces, and $\coprod$ for their disjoint union
as algebraic spaces or stacks. In the latter, each
piece is open and closed; in the former, different
pieces may be related by specialization.}
\[
f_n^{-1}(\mathcal Z_\Gamma)
=
\bigsqcup_{[\delta]\in
\mathscr D_n(\Gamma)/\Aut\Gamma}
\mathscr T_{[\delta]}.
\]

Gluing the paired attaching sections induces
\[
[\mathscr X_\delta/G_\delta]
\xrightarrow{\sim}
\mathscr T_{[\delta]}.
\]
Indeed, etale-locally a family in
$\mathscr T_{[\delta]}$ admits an identification
of its decorated dual graph with
$\widetilde\Gamma$. Separating its nodes gives
the inverse construction, and changes of
identification are exactly
$\Aut(\widetilde\Gamma,\ell)\cong G_\delta$.
Both constructions commute with base change.

Taking coarse moduli spaces and using additivity,
\[
\chi\bigl(f_n^{-1}(\mathcal Z_\Gamma)\bigr)
=
\sum_{[\delta]\in
\mathscr D_n(\Gamma)/\Aut\Gamma}
\chi(X_\delta/G_\delta)
=
\chi(X_{\Gamma,n}/\Aut\Gamma).\tag*{\qedhere}
\]
\end{proof}

\begin{remark}
The disjoint union $X_{\Gamma,n}$ does not retain
the specialization relations among the decorated
pieces inside $f_n^{-1}(\mathcal Z_\Gamma)$.
Lemma~\ref{lem:graph-stratum-euler} identifies their
Euler characteristics, not the two moduli spaces.
\qrem
\end{remark}

By \eqref{eq:decomposition_A_g_1} and
Lemma~\ref{lem:graph-stratum-euler}, we have
\[
A^\chi_{g,\Gamma}(z)
=
\sum_{n\ge0}
\chi(X_{\Gamma,n}/\Aut\Gamma)\frac{z^n}{n!}.
\]

By the Burnside--Lefschetz formula,
\begin{equation}
A^\chi_{g,\Gamma}(z)
=
\frac1{|\Aut\Gamma|}
\sum_{\sigma\in\Aut\Gamma}
F_{\Gamma,\sigma}(z),
\qquad
F_{\Gamma,\sigma}(z)
:=
\sum_{n\ge0}
\chi(X_{\Gamma,n}^{\sigma})\frac{z^n}{n!}.
\label{eq:graph-burnside}
\end{equation}
We analyze separately the identity and nonidentity terms.

\medskip
\noindent{\bf The identity contribution.}
For a stable pair $(h,d)$, with $2h-2+d>0$, define the vertex series
associated with the coarse component spaces introduced in
Section~\ref{sec:decorated-expansions} by
\begin{equation}
V^\chi_{h,d}(z)
:=
\sum_{n\ge0}
\chi(\mathcal V_{h,d,n})\frac{z^n}{n!}.
\label{eq:eulervertexseries}
\end{equation}

\begin{lemma}
\label{lem:eulervertex}
For every stable pair $(h,d)$,
\begin{equation}
V^\chi_{h,d}(z)
=
\sum_{k\ge0}
\chi(\M_{h,d+k})
\frac{(e^z-1)^k}{k!}.
\label{eq:eulervertex-stirling}
\end{equation}
Moreover, $V^\chi_{h,d}(z)$ is an exponential polynomial\,\footnote{An
exponential polynomial in $z$ is a finite sum
$\sum_{j=1}^r p_j(z)e^{\lambda_j z}$, with $p_j$ polynomials and
$\lambda_j\in\mathbb C$. Such a function is entire. The vertex and
fixed-sector series occurring here have constant $p_j$, so they are
finite linear combinations of exponentials.}, and hence
an entire function. In genus zero,
\begin{equation}
V^\chi_{0,d}(z)
=
(-1)^{d-3}(d-3)!\,e^{-(d-2)z},
\qquad d\ge3.
\label{eq:eulervertex-rational}
\end{equation}
In particular,
\begin{equation}
V^\chi_{0,3}(z)=e^{-z}.
\label{eq:eulervertex-trivalent}
\end{equation}
\end{lemma}

\begin{proof}
Stratify $\mathcal V_{h,d,n}$ according to the partition of the $n$
light labels determined by collisions. If the light markings occupy
exactly $k$ distinct points, then for each set partition of $[n]$
into $k$ nonempty blocks the corresponding stratum is isomorphic to
$M_{h,d+k}$. Hence
\[
\chi(\mathcal V_{h,d,n})
=
\sum_{k=0}^n
\genfrac{\{}{\}}{0pt}{}{n}{k}
\chi(\M_{h,d+k}),
\]
where
$\genfrac{\{}{\}}{0pt}{}{n}{k}$
denotes a Stirling number of the second kind. Using
\[
\sum_{n\ge k}
\genfrac{\{}{\}}{0pt}{}{n}{k}
\frac{z^n}{n!}
=
\frac{(e^z-1)^k}{k!}
\]
gives \eqref{eq:eulervertex-stirling}.

For $h=0$,
\[
\chi(\M_{0,N})
=
(-1)^{N-3}(N-3)!,
\qquad N\ge3,
\]
and substitution into \eqref{eq:eulervertex-stirling} gives
\[
V^\chi_{0,d}(z)
=
(-1)^{d-3}(d-3)!
\sum_{k\ge0}
(-1)^k
\binom{d+k-3}{k}
(e^z-1)^k,
\]
which is \eqref{eq:eulervertex-rational}.

Let now $h\ge1$. The Harer--Zagier orbifold Euler characteristic is
\[
\chi_{\mathrm{orb}}(\M_{h,N})
=
(-1)^N
\frac{{\sf B}_{2h}}{2h(2h-2)!}
(2h+N-3)!,\qquad \text{where }
{\sf B}_{2h}=(-1)^{h-1}\frac{2(2h)!}{(2\pi)^{2h}}\zeta(2h);
\]
see \cite{HZ}. Moreover, the ordinary and orbifold Euler
characteristics coincide once $N\ge2h+3$, since an $N$-pointed smooth
genus-$h$ curve is then automorphism-free; see
\cite[Remark~4.4]{BH}. Thus the ordinary sequence
$\chi(\M_{h,d+k})$ differs from the Harer--Zagier sequence for only
finitely many values of $k$. Setting $y=e^z-1$, the infinite
Harer--Zagier contribution in \eqref{eq:eulervertex-stirling} sums to
\[
(-1)^d
\frac{{\sf B}_{2h}}{2h(2h-2)!}
(2h+d-3)!\,
(1+y)^{-(2h+d-2)}.
\]
Consequently
\[
V^\chi_{h,d}(z)
=
Q_{h,d}(e^z-1)
+
(-1)^d
\frac{{\sf B}_{2h}}{2h(2h-2)!}
(2h+d-3)!\,
e^{-(2h+d-2)z}
\]
for a polynomial $Q_{h,d}$. Hence $V^\chi_{h,d}$ is an exponential
polynomial and, in particular, entire.
\end{proof}

For an inserted two-valent component carrying $m\ge1$ light labels,
Section~\ref{sec:decorated-expansions} gives the factor
$(\mathbb G_m)^{m-1}$. Since
\[
\chi\bigl((\mathbb G_m)^{m-1}\bigr)
=
\begin{cases}
1,&m=1,\\
0,&m\ge2,
\end{cases}
\]
only components carrying a single label contribute.
For a chain of length $k$, the $k$ labels admit
$k!$ orders along the chain. Hence its exponential
generating series is
\begin{equation}
\sum_{k\ge0}k!\frac{z^k}{k!}
=1+z+z^2+z^3+\dots=\frac1{1-z}.
\label{eq:eulerpropagator}
\end{equation}
We shall refer to this universal edge factor as the
\emph{Euler propagator}\footnote{We borrow this terminology from physics,
where propagators are associated with the edges of Feynman diagrams
in diagrammatic expansions.}.

Before taking the quotient by graph symmetries, each vertex and each
edge is a separate slot for component data. Distributing the light
labels among these slots gives a product of exponential generating
series. The vertex data and the chains on distinct edges are therefore
independent, and the
identity term in \eqref{eq:graph-burnside} is
\begin{equation}
F_{\Gamma,1}(z)
=
\left(
\prod_{v\in\operatorname{Vert}(\Gamma)}
V^\chi_{h_v,d_v}(z)
\right)
(1-z)^{-E(\Gamma)}.
\label{eq:identitysector}
\end{equation}

\medskip
\noindent{\bf The nonidentity contributions.}
For $\sigma\in\Aut\Gamma$, let $r_\Gamma(\sigma)$ denote the number of
edges of $\Gamma$ whose two half-edges are fixed individually by $\sigma$.

\begin{lemma}[Fixed sectors]
\label{lem:euler-fixed-sectors}
For every $\sigma\in\Aut\Gamma$ there exists an entire exponential
polynomial $E_{\Gamma,\sigma}(z)$ such that
\begin{equation}
F_{\Gamma,\sigma}(z)
=
E_{\Gamma,\sigma}(z)
(1-z)^{-r_\Gamma(\sigma)}.
\label{eq:euler-fixed-sector}
\end{equation}
If $\sigma\ne1$, then
\[
r_\Gamma(\sigma)\le E(\Gamma)-1.
\]
\end{lemma}

\begin{proof}
We examine separately the graph vertices and the chains inserted on
graph edges. Throughout, $\sigma$ acts on the graph but fixes each of
the $n$ light labels individually. Thus a fixed point must lie in a
component $X_\delta$ for which the decoration $\delta$ itself is fixed.
In particular, a vertex carrying a light label cannot be transported
to a different vertex.

\smallskip
\noindent\emph{1. Vertices moved by $\sigma$.}
Let $\mathcal O=\{v_0,\ldots,v_{\ell-1}\}$ be an orbit
of vertices of length $\ell>1$, with
$\sigma(v_j)=v_{j+1}$ (indices modulo $\ell$).

First, no vertex in $\mathcal O$ can carry a light
label, since
$\sigma$ fixes each light label individually. 
Moreover, all vertices in $\mathcal O$ have the same genus
$h$ and valence $d$. Their component space is
a product of $\ell$ copies of $M_{h,d}$, on which
$\sigma$ acts by cyclically permuting the factors,
together with the corresponding reindexing of
the heavy markings.

Let
$\pi_{\mathcal O}\in\mathfrak S_d$ be the permutation of
the heavy branches at $v_0$ induced by $\sigma^\ell$.
A fixed tuple is determined by its component at $v_0$,
which must be fixed by $\pi_{\mathcal O}$. Thus the fixed
locus of the vertex orbit is isomorphic to
\(
M_{h,d}^{\pi_{\mathcal O}},
\)
where the fixed locus is taken in the coarse moduli space.
Its points represent curves admitting an automorphism
realizing $\pi_{\mathcal O}$, as in Step~3 below.

The orbit therefore contributes the constant
\(
\chi\bigl(M_{h,d}^{\pi_{\mathcal O}}\bigr).
\)
If $\pi_{\mathcal O}=1$, the fixed locus is a diagonal
copy of $M_{h,d}$, after the corresponding reindexing
of the heavy markings, and its contribution is
$\chi(M_{h,d})$.

\smallskip
\noindent\emph{2. {Fixed vertices with individually fixed branches.}}
At such a vertex there is no additional condition on the component
or its light markings. Its contribution is exactly
$V^\chi_{h_v,d_v}(z)$, which is entire by
Lemma~\ref{lem:eulervertex}.

\smallskip
\noindent\emph{3. Fixed vertices with permuted branches.}
Let $\pi\ne1$ be the permutation induced by $\sigma$ on
the $d$ heavy branches at a fixed vertex of genus $h$.
We regard the attaching points as heavy markings and
consider the coarse moduli space $B:=M_{h,d}$ of smooth
stable $d$-pointed curves of genus $h$.

\smallskip
\noindent\emph{Automorphism strata.}
For $b=[(C;p_1,\ldots,p_d)]\in B$, let
\[
A_C:=\Aut(C,\{p_1,\ldots,p_d\})
\]
be the group of automorphisms preserving the set of
heavy markings. Its action on these markings defines
a homomorphism $\rho_C:A_C\to\mathfrak S_d$.
We set
\[
H_C:=\ker\rho_C,\qquad
R_C:=\rho_C^{-1}(\pi),
\]
and denote by $B^\pi\subseteq B$ the locus where
$R_C\ne\varnothing$.

We wish to stratify $B^\pi$ so that the following
data are constant on each stratum, up to relabeling:

\begin{enumerate}
\item The group $A_C$ together with its homomorphism
$\rho_C:A_C\to\mathfrak S_d$, and hence the subgroup
$H_C$ and the subset $R_C$.

\item The fixed-point incidence numbers
\[
q_{J,\beta}(C):=
\#\left(
\operatorname{Fix}(\beta)\cap
\bigcap_{\alpha\in J}\operatorname{Fix}(\alpha)
\cap C^\circ
\right),
\]
for every $\varnothing\ne J\subseteq R_C$ and
$\beta\in H_C$, where
$C^\circ:=C\setminus\{p_1,\ldots,p_d\}$.
\end{enumerate}

Let us  justify the existence of such a finite
stratification.

First, by stability, $2h-2+d>0$, the orders of the
groups $A_C$ are uniformly bounded. Indeed:
\begin{itemize}
\item If $h=0$, then $d\ge3$. An automorphism of
$\mathbb P^1$ is determined by the images of three
distinct points, so $\rho_C$ is injective and
$|A_C|\le d!$.
\item If $h=1$, then $d\ge1$. The subgroup fixing
one heavy point has order at most $6$, while its
index in $A_C$ is at most $d$. Thus $|A_C|\le6d$.
\item If $h\ge2$, the Hurwitz bound gives
$|A_C|\le84(h-1)$.
\end{itemize}
Consequently, only finitely many group-theoretic
configurations $(A_C,\rho_C)$ can occur.

Second, automorphisms and their fixed points vary
algebraically in families. Working étale-locally on
the moduli stack, the automorphism groups become
finite étale after a suitable stratification.
Since $\pi\ne1$, every $\alpha\in R_C$ is
nonidentity and has finitely many fixed points.
The conditions
\[
\alpha(x)=x,\qquad
\beta(x)=x,\qquad
x=p_j,
\]
for $\alpha\in R_C$, $\beta\in H_C$, and
$1\le j\le d$, are algebraic. Their simultaneous
imposition describes precisely the intersections
entering the definition of $q_{J,\beta}(C)$,
including the exclusion of the heavy markings.

By constructibility, a further finite refinement
makes all these incidence numbers constant.
We therefore obtain a finite stratification into
locally closed subsets
\[
B^\pi=\bigsqcup_{i=1}^N S_i
\]
on which both the automorphism data
$(A_C,\rho_C)$ and the full collection of
integers $q_{J,\beta}(C)$ are constant up to
relabeling.

\smallskip
\noindent\emph{Contribution of a single stratum.}
Fix a stratum $S_i$ and a representative
$b=[(C;p_1,\ldots,p_d)]\in S_i$.
For simplicity, write $H=H_C$ and $R=R_C$.

A coarse point with light markings
$(x_1,\ldots,x_n)$ is fixed by $\pi$ precisely
when some $\alpha\in R$ fixes every $x_j$.

For each $\alpha\in R$, let
\(
Q_\alpha:=\operatorname{Fix}(\alpha)\cap C^\circ.
\)
This is a finite set, since $\alpha$ is
nonidentity. Writing $q_\alpha=\#Q_\alpha$,
the generating series for assignments fixed
by this particular automorphism is
\[
\sum_{n\ge0}q_\alpha^n\frac{z^n}{n!}
=e^{q_\alpha z}.
\]
Indeed, collisions are allowed, so each light
label independently chooses a point of $Q_\alpha$,
with the convention $0^0=1$.

To account for all realizations of $\pi$ and
for isomorphic assignments, put
\[
Y_n:=\bigcup_{\alpha\in R}Q_\alpha^{\,n}.
\]
The union avoids counting an assignment more
than once, while the quotient by $H$ identifies
isomorphic assignments. Thus the coarse fiber
is $Y_n/H$.

By Burnside's formula and inclusion--exclusion,
its exponential generating series is
\[
f_i(z)
:=\sum_{n\ge0}\#(Y_n/H)\frac{z^n}{n!}
=\frac1{|H|}
\sum_{\beta\in H}
\sum_{\varnothing\ne J\subseteq R}
(-1)^{|J|+1}e^{q_{J,\beta}z}.
\]
By construction of the stratification, the
groups $H$, $R$ and all the integers
$q_{J,\beta}$ are constant up to relabeling.
Hence $f_i(z)$ is independent of the chosen
$C\in S_i$.

In particular, $f_i$ is a finite linear
combination of exponentials, and therefore entire.

\smallskip
\noindent\emph{Summation over the strata.}
Write
\(
f_i(z)=\sum_{n\ge0}a_{i,n}\frac{z^n}{n!},
\)
and let $Z_{n,i}$ be the locus of fixed
configurations with $n$ light markings lying
over $S_i$. The forgetful map
\(
p_{n,i}:Z_{n,i}\longrightarrow S_i
\)
has finite coarse fibers of constant cardinality
$a_{i,n}$.

Although monodromy may permute the available
fixed points as $C$ varies, this does not affect
the Euler characteristic. Indeed, integration
of constructible functions with respect to Euler
characteristic gives
\[
\chi(Z_{n,i})
=\int_{S_i}
\chi\bigl(p_{n,i}^{-1}(b)\bigr)\,d\chi(b)
=a_{i,n}\chi(S_i).
\]
Equivalently, finite coverings have Euler
characteristic equal to their degree times that
of the base, regardless of their monodromy.

Summing over $n$ and over the finitely many
strata, we obtain
\[
\sum_{n\ge0}\chi(Z_n^\pi)\frac{z^n}{n!}
=\sum_{i=1}^N\chi(S_i)
  \sum_{n\ge0}a_{i,n}\frac{z^n}{n!}
=\sum_{i=1}^N\chi(S_i)f_i(z).
\]
Since each $f_i(z)$ is an entire exponential
polynomial, the total vertex contribution is as well.

For an example of computation, see Example \ref{ex:fixed-point-monodromy} below.

\smallskip
\noindent\emph{4. Chains on the edges.}
There are three possibilities, illustrated in
Figure~\ref{fig:euler-fixed-edge-types}.
If both half-edges of an edge are fixed individually, any inserted
chain is allowed and its series is $(1-z)^{-1}$.
If the edge is moved to a different edge, its chain must be empty:
a nonempty chain would contain a label transported to that other edge.
Its contribution is therefore $1$.
Finally, if the edge is preserved but its half-edges are exchanged,
the chain is reversed. Any noncentral component would be exchanged
with another component, which is impossible because it carries a
nonempty set of individually fixed labels. Thus only the empty chain
or a single central component is allowed. For $m\ge1$ labels on the
central component, reversal acts by simultaneous inversion on
$(\mathbb G_m)^{m-1}$ and has $2^{m-1}$ fixed points. The edge series is
\beq
\label{eq:contribution_reversed_edge}
1+\sum_{m\ge1}2^{m-1}\frac{z^m}{m!}
=\frac{1+e^{2z}}2,
\eeq
which is entire.

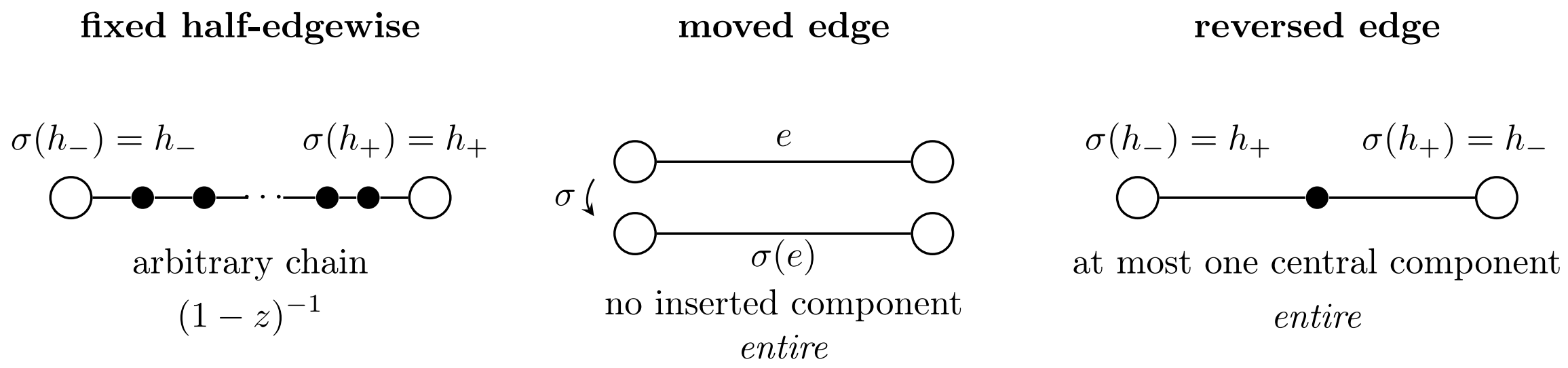
\begin{figure}[ht]
\centering

\begin{tikzpicture}[
  fixedgraph/.style={
    circle,draw,line width=.65pt,fill=white,
    minimum size=4mm,inner sep=0pt
  },
  inserted/.style={
    circle,fill=black,
    minimum size=2.2mm,inner sep=0pt
  },
  fixededge/.style={
    line width=.65pt
  },
  action/.style={
    ->,>=stealth,line width=.55pt
  },
  font=\small
]

\begin{scope}[xshift=0cm]

\node at (0,1.65) {\textbf{fixed half-edgewise}};

\node[fixedgraph] (a1) at (-1.75,0) {};
\node[fixedgraph] (b1) at ( 1.75,0) {};

\node[inserted] (u1) at (-1.05,0) {};
\node[inserted] (u2) at (-.45,0) {};
\node at (.15,0) {$\cdots$};
\node[inserted] (u3) at (.75,0) {};
\node[inserted] (u4) at (1.15,0) {};

\draw[fixededge] (a1)--(u1)--(u2);
\draw[fixededge] (u2)--(-.02,0);
\draw[fixededge] (.32,0)--(u3)--(u4)--(b1);

\node[above=7pt] at (-1.42,0)
  {$\sigma(h_-)=h_-$};
\node[above=7pt] at (1.42,0)
  {$\sigma(h_+)=h_+$};

\node at (0,-.65) {arbitrary chain};
\node at (0,-1.15) {$\displaystyle (1-z)^{-1}$};

\end{scope}

\begin{scope}[xshift=5.2cm]

\node at (0,1.65) {\textbf{moved edge}};

\node[fixedgraph] (a2) at (-1.45,.35) {};
\node[fixedgraph] (b2) at ( 1.45,.35) {};
\draw[fixededge] (a2)--(b2);
\node at (0,.58) {$e$};

\node[fixedgraph] (c2) at (-1.45,-.35) {};
\node[fixedgraph] (d2) at ( 1.45,-.35) {};
\draw[fixededge] (c2)--(d2);
\node at (0,-.58) {$\sigma(e)$};

\draw[action]
  (-1.85,.18)
  to[bend right=35]
  (-1.85,-.18);

\node[left=2pt] at (-1.82,0) {$\sigma$};

\node at (0,-1.05) {no inserted component};
\node at (0,-1.45) {\emph{entire}};

\end{scope}

\begin{scope}[xshift=10.4cm]

\node at (0,1.65) {\textbf{reversed edge}};

\node[fixedgraph] (a3) at (-1.75,0) {};
\node[fixedgraph] (b3) at ( 1.75,0) {};
\node[inserted] (u3c) at (0,0) {};

\draw[fixededge] (a3)--(u3c)--(b3);

\node[above=7pt] at (-1.35,0)
  {$\sigma(h_-)=h_+$};
\node[above=7pt] at (1.35,0)
  {$\sigma(h_+)=h_-$};

\node at (0,-.65) {at most one central component};
\node at (0,-1.15) {\emph{entire}};

\end{scope}

\end{tikzpicture}

\caption{The three possible behaviors of a graph edge under
$\sigma$. If the two half-edges are fixed individually, an arbitrary
chain of inserted components is allowed and produces the full Euler
propagator $(1-z)^{-1}$. If the edge is moved, a
$\sigma$-fixed decoration cannot contain any inserted component along
it. If the edge is preserved but reversed, all noncentral inserted
components are exchanged, so that at most one central component may
remain. Thus only the first case produces a nonentire factor.}
\label{fig:euler-fixed-edge-types}

\end{figure}

It follows that the only nonentire factors are the propagators attached
to the $r_\Gamma(\sigma)$ edges whose two half-edges are fixed
individually, proving \eqref{eq:euler-fixed-sector}.

Finally, if every half-edge of the connected graph were fixed
individually, then every edge and every incident vertex would be fixed,
and $\sigma$ would be the identity. Hence, for $\sigma\ne1$,
$r_\Gamma(\sigma)\le E(\Gamma)-1$.
\end{proof}

\begin{example}[A four-valent rational vertex with pairwise exchanged branches]
\label{ex:fixed-point-monodromy} Let us exemplify the computation of contribution 
at a vertex with permuted branches. Below we use the same notations as in point 3 of the proof of Lemma \ref{lem:euler-fixed-sectors}.

Consider the family of four-pointed rational curves
\[
(C;p_1,p_2,p_3,p_4)
=(\mathbb P^1;0,\infty,1,\lambda),
\qquad
\lambda\in B:=\mathbb C\setminus\{0,1\},
\]
parametrized by $B\simeq M_{0,4}$, and let
$\pi=(12)(34)$.
For every $\lambda\in B$, there is a unique automorphism
realizing $\pi$, namely
\(
\alpha_\lambda(x)={\lambda}/{x}.
\)
Its fixed points satisfy $x^2=\lambda$, hence
\(
Q_{\alpha_\lambda}
=\{+\sqrt{\lambda},-\sqrt{\lambda}\}.
\)
Since $H$ is trivial, the fiber series is simply
$f_\lambda(z)=e^{2z}$, independently of $\lambda$.

However, the two fixed points cannot be labeled globally
as continuous functions of $\lambda$: when $\lambda$
travels around the origin, monodromy exchanges
$+\sqrt{\lambda}$ and $-\sqrt{\lambda}$.
Thus the family of fixed points is not the disjoint union
of two copies of $B$, but rather the connected double cover
\[
\widetilde B
=\{(\lambda,x)\in B\times\mathbb C^\times:
x^2=\lambda\}
\longrightarrow B.
\]
Nevertheless,
\(
\chi(\widetilde B)=2\chi(B)=-2,
\)
since $\chi(B)=-1$.

With $n$ labeled light markings, collisions being allowed,
each label independently chooses one of the two fixed
points. The resulting family $Z_n$ is therefore a finite
étale cover of $B$ of degree $2^n$, possibly disconnected,
and
\(
\chi(Z_n)=2^n\chi(B)=-2^n.
\)
Consequently,
\[
\sum_{n\ge0}\chi(Z_n)\frac{z^n}{n!}
=\chi(B)e^{2z}=-e^{2z}.
\]
This illustrates that monodromy may permute the available
fixed points without affecting the Euler-characteristic
formula: the contribution depends only on the cardinality
of the fibers, not on their global triviality.
\qetr
\end{example}

By \eqref{eq:identitysector}, the identity term has exactly
$E(\Gamma)$ propagators. By Lemma~\ref{lem:euler-fixed-sectors}, every
nonidentity term has at most $E(\Gamma)-1$ propagators. Substituting
into \eqref{eq:graph-burnside} proves the following.

\begin{proposition}[Contribution of a stabilization graph]
\label{prop:graphcontribution}
For every stabilization graph $\Gamma$ of genus $g\ge2$,
\begin{equation}
A^\chi_{g,\Gamma}(z)
=
\frac1{|\Aut\Gamma|}
\left(
\prod_{v\in\operatorname{Vert}(\Gamma)}
V^\chi_{h_v,d_v}(z)
\right)
(1-z)^{-E(\Gamma)}
+
O\!\left((1-z)^{-E(\Gamma)+1}\right)
\label{eq:graphasym}
\end{equation}
as $z\to1$ in a punctured neighborhood of $1$.

Moreover, $A^\chi_{g,\Gamma}(z)$ is a finite sum of entire
exponential polynomials multiplied by powers of $(1-z)^{-1}$.
\end{proposition}

\subsection{Summary of the computational procedure}
\label{sec:graph-recipe}

The preceding constructions give a procedure for computing
$A_g^\chi(z)$ and determining its analytic structure. We
summarize the steps and the role of each result.

\begin{enumerate}
\item \emph{Enumerate the stabilization graphs.}
Lemma~\ref{lem:all-light-stability} excludes rational tails,
so every all-light curve is obtained from its stabilized curve
by inserting rational chains and distributing the markings.
Enumerate the finitely many graphs $\Gamma$ of genus $g$
and determine their automorphism groups.
Lemma~\ref{lem:maxedges} gives the bound
$E(\Gamma)\le3g-3$ and identifies the graphs attaining it.

\item \emph{Pass to the parameter spaces.}
For each graph $\Gamma$, the decorated expansions define
the parameter space $X_{\Gamma,n}$ with its $\Aut\Gamma$-action.
Lemma~\ref{lem:graph-stratum-euler} identifies
\(
\chi\!\left(f_n^{-1}(\mathcal Z_\Gamma)\right)
=
\chi\!\left(X_{\Gamma,n}/\Aut\Gamma\right).
\)
This reduces the contribution of the geometric stratum
to an Euler-characteristic calculation on a finite quotient.
\item \emph{Compute the identity sector.}
Lemma~\ref{lem:eulervertex} expresses each vertex series as
\[
V^\chi_{h,d}(z)
=
\sum_{k\ge0}\chi(\M_{h,d+k})
\frac{(e^z-1)^k}{k!}
\]
and proves that it is an entire exponential polynomial.
The chain calculation in \eqref{eq:eulerpropagator} supplies
one factor $(1-z)^{-1}$ for each edge. Hence
\[
F_{\Gamma,1}(z)
=
\prod_{v\in\operatorname{Vert}(\Gamma)}
V^\chi_{h_v,d_v}(z)\,
(1-z)^{-E(\Gamma)}.
\]

\item \emph{Compute the nonidentity sectors.}
For each nonidentity symmetry $\sigma\in\Aut\Gamma$, apply
the vertex and edge analysis in the proof of
Lemma~\ref{lem:euler-fixed-sectors}. The light labels remain
individually fixed throughout. The lemma shows that
\[
F_{\Gamma,\sigma}(z)
=
E_{\Gamma,\sigma}(z)
(1-z)^{-r_\Gamma(\sigma)},
\]
where $E_{\Gamma,\sigma}$ is an entire exponential polynomial
and $r_\Gamma(\sigma)$ counts the edges whose half-edges are
fixed individually. Its proof describes how to obtain the
numerator from the corresponding fixed component data.

\item \emph{Average over symmetries and sum over graphs.}
Lemma~\ref{lem:burnside} converts these fixed-locus calculations
into the Euler characteristic of the finite quotient.
Equivalently, \eqref{eq:graph-burnside} gives
\[
A_{g,\Gamma}^\chi(z)
=
\frac{1}{|\Aut\Gamma|}
\sum_{\sigma\in\Aut\Gamma}F_{\Gamma,\sigma}(z).
\]
Summing these expressions over the stabilization graphs yields
$A_g^\chi(z)$.
\end{enumerate}

For exact formulas, all symmetry sectors must be retained.
For the leading singular behavior, the procedure simplifies:
Lemma~\ref{lem:euler-fixed-sectors} bounds every nonidentity
sector by pole order $E(\Gamma)-1$ at $z=1$.
Proposition~\ref{prop:graphcontribution} therefore isolates
the highest possible pole term of each graph contribution
in the identity sector. Together with
Lemma~\ref{lem:maxedges}, this reduces the analysis of the
highest possible pole order to rational trivalent graphs.
Finally, Theorem~\ref{thm:pde} recovers the ordinary Euler
series by composition with the universal heavy/light
coordinate.

\section{Genus two: stabilization-graph formulas and Euler asymptotics}
\label{sec:g2}

We now apply the procedure developed in
Section~\ref{sec:graph-formalism} to the seven stabilization graphs
of genus two, obtaining an explicit formula for $A_2^\chi(z)$.
Theorem~\ref{thm:pde} then yields the ordinary series
$G_2^\chi(t)$, and both formulas are used to derive
the corresponding asymptotic expansions.
The genus-three calculation is presented in
Appendix~\ref{app:genus3}.

\subsection{The vertex and edge contributions required by the seven graphs}
\label{sec:g2-graph-data}
There are seven stabilization graphs of genus two, displayed in
Figure~\ref{fig:stable-graphs-g2}. To compute their exact Burnside
contributions, we first determine which local contributions their
vertices, edges, and symmetries can require. This list can be read
from the drawings: the vertex types $(h,d)$ are
$(2,0)$, $(1,1)$, $(1,2)$, $(0,3)$, and $(0,4)$; a loop reversal
exchanges its two branches, while parallel edges or equal vertices
may be permuted. 

The required building blocks are the following.
\begin{enumerate}[label=\textup{(\roman*)},leftmargin=*]
\item \emph{Edges.} An edge fixed half-edgewise contributes the
propagator $P$; a reversed edge contributes $K$; a moved edge has an
empty chain and contributes $1$.
\item \emph{Positive-genus vertices fixed with all branches fixed.}
We need $T=V^\chi_{2,0}$, $U=V^\chi_{1,1}$, and
$V=V^\chi_{1,2}$.
\item \emph{An elliptic vertex with its two branches exchanged.}
We need the fixed-locus series $W$ for exchanging the two
heavy markings while fixing every light label.
\item \emph{Fixed rational vertices.} At valence three we need the
identity, a transposition, and a three-cycle of the heavy branches.
At valence four we need the identity, a single transposition, a
double transposition, and a four-cycle. These are precisely the
branch permutations induced by the graph symmetries.
\item \emph{Moved vertices.} They carry no light labels. Their
component data must agree around each vertex cycle. In genus two
we need only the diagonal copy of $M_{1,1}$, whose Euler
characteristic is $1$, or copies of the point $M_{0,3}$.
\end{enumerate}
We now compute these building blocks in the order listed above.
The calculations rely on the collision formula of
Lemma~\ref{lem:eulervertex} and the fixed-sector analysis of
Lemma~\ref{lem:euler-fixed-sectors}.

\smallskip
\noindent\emph{Edge contributions.}
We use the edge factors computed in
\eqref{eq:eulerpropagator} and
\eqref{eq:contribution_reversed_edge}:
\begin{equation}\label{eq:g2-P-K}
P(z)=\frac1{1-z},
\qquad
K(z)=\frac{1+e^{2z}}2.
\end{equation}

\smallskip
\noindent\emph{Positive-genus identity contributions.}
By Lemma~\ref{lem:eulervertex},
the identity contribution of a vertex of genus $h$
and valence $d$ is
\[
V^\chi_{h,d}(z)
=\sum_{k\ge0}\chi(M_{h,d+k})\frac{D^k}{k!},
\qquad D=e^z-1.
\]
In particular, setting
\[
U=V^\chi_{1,1},\qquad
V=V^\chi_{1,2},\qquad
T=V^\chi_{2,0},
\]
we obtain the following explicit expressions:
\begin{align}
U&=-\frac{e^{-z}}{12}+\frac{13}{12}+\frac{11D}{12}
+\frac{D^2}{12}-\frac{D^3}{12},\label{eq:g2-U}\\
V&=\frac{e^{-2z}}{12}+\frac{11}{12}+
\frac D6-\frac{D^2}{4},\label{eq:g2-V}\\
T&=-\frac{e^{-2z}}{240}+\frac{241}{240}+\frac{239D}{120}
+\frac{81D^2}{80}-\frac{D^3}{60}
-\frac{7D^4}{48}-\frac{D^5}{40}-\frac{D^6}{240}.\label{eq:g2-T}
\end{align}

For the computations of $U$ and $V$, one can use the smooth genus-one Euler characteristics 
already invoked in Lemma \ref{lem:elliptic-collisions}, namely
\(
\chi(M_{1,1})=\chi(M_{1,2})=1,\)
\(
\chi(M_{1,3})=\chi(M_{1,4})=0,
\)
together with the Harer--Zagier formula for the tail\footnote{In genus one, with five distinct markings no
nontrivial automorphism remains.} $\chi(M_{1,k})=(-1)^k(k-1)!/12$ for $k\ge5$;
see \cite[Section~4]{BGP}\cite{HZ}. 
For the computation of $T$, the ordinary Euler characteristics of the smooth genus-two spaces are
\[
\begin{array}{c|rrrrrrr}
n&0&1&2&3&4&5&6\\\hline
\chi(M_{2,n})&1&2&2&0&-4&0&-24
\end{array}
\]
and $\chi(M_{2,n})=(-1)^{n+1}(n+1)!/240$ for $n\ge7$;
see \cite[Section~3]{BGP}. The tail is the Harer--Zagier formula,
since seven distinct markings eliminate automorphisms. Correcting its
first seven terms in the collision transform yields \eqref{eq:g2-T}.

\begin{rem}For the exponential polynomials $U$ and $V$, there is an alternative derivation. 
In fact, with $F$ as in the proof of Lemma~\ref{lem:elliptic-collisions},
we have $U=F'(D)$ and $V=F''(D)$, respectively.
\qrem
\end{rem}

\smallskip
\noindent\emph{The elliptic branch exchange.}
Consider the involution $\tau$ on $\mathcal V_{1,2,n}$
which exchanges the two heavy markings and fixes
all light labels individually. We wish to compute
\[
W(z):=\sum_{n\ge0}
\chi\bigl(\mathcal V_{1,2,n}^{\tau}\bigr)\frac{z^n}{n!}.
\]
We proceed in four steps.

\smallskip
\noindent\emph{i.~Geometric fixed-point condition.}
A configuration $(C;p,q;x_1,\ldots,x_n)$ represents
a fixed coarse point precisely when there exists
$\alpha\in\Aut(C)$ such that
\[
\alpha(p)=q,\qquad
\alpha(q)=p,\qquad
\alpha(x_i)=x_i\quad(1\le i\le n).
\]This is exactly the situation considered in
Step~3 of the proof of Lemma~\ref{lem:euler-fixed-sectors}.

Notice that an automorphism exchanging $p$ and $q$
always exists: using the group law on $C$, the involution
\(
\alpha_0(x)=p+q-x
\)
has exactly four fixed points, given by $2x=p+q$,
none of which coincides with $p$ or $q$.
More generally, since $\alpha$ exchanges the distinct
points $p,q$, it is nonidentity, and every nonidentity
automorphism of an elliptic curve has at most four
fixed points. Consequently, a fixed configuration
can contain at most four distinct light points.

\smallskip
\noindent\emph{ii.~Collision partitions.}
Partition the $n$ light labels into $k$ nonempty
collision blocks. Each block corresponds to one
distinct light point, so the corresponding stratum
is isomorphic to $M_{1,k+2}$.
Since $\tau$ fixes each light label, it also
preserves every collision block individually.

Put
\(
c_k:=\chi(M_{1,k+2}^{\tau}).
\)
The preceding fixed-point bound implies
$c_k=0$ for $k\ge5$.

Recall that $D=e^z-1$ is the exponential generating
series for nonempty sets of labels. Summing over
collision partitions gives
\beq
\label{eq:g2-W-0}
W(z)=\sum_{k=0}^4 c_k\frac{D^k}{k!}.
\eeq

\smallskip
\noindent\emph{iii.~Burnside and Bini--Gaiffi--Polito formulas.}
By Lemma~\ref{lem:burnside}, applied to the action
of $\mathfrak S_2=\langle\tau\rangle$ on $M_{1,k+2}$,
\[
\chi(M_{1,k+2}/\mathfrak S_2)
=\frac12\left(
\chi(M_{1,k+2})+c_k
\right)\qquad\Longrightarrow\qquad 
c_k=
2\chi(M_{1,k+2}/\mathfrak S_2)
-\chi(M_{1,k+2}).
\]
Using the ordinary Euler characteristics and
the quotient Euler characteristics computed in
\cite[Proposition~2.14]{BGP}, we obtain
\[
c_0=1,\qquad
c_1=c_2=c_3=c_4=2.
\]
\smallskip
\noindent\emph{iv.~Resulting series.}
Substituting these values into \eqref{eq:g2-W-0},
we obtain
\begin{equation}\label{eq:g2-W}
W(z)=1+2D+D^2+\frac{D^3}{3}+\frac{D^4}{12}.
\end{equation}
In particular, $W$ is an entire exponential polynomial,
as predicted by Lemma~\ref{lem:euler-fixed-sectors}.

For the elliptic loop graph, reversal also acts
on the chain inserted along the loop. Its contribution
is the reversed-edge factor
\(
K(z)={(1+e^{2z})}/2.
\)
Thus the complete contribution of the fixed
elliptic loop sector is $W(z)K(z)$.

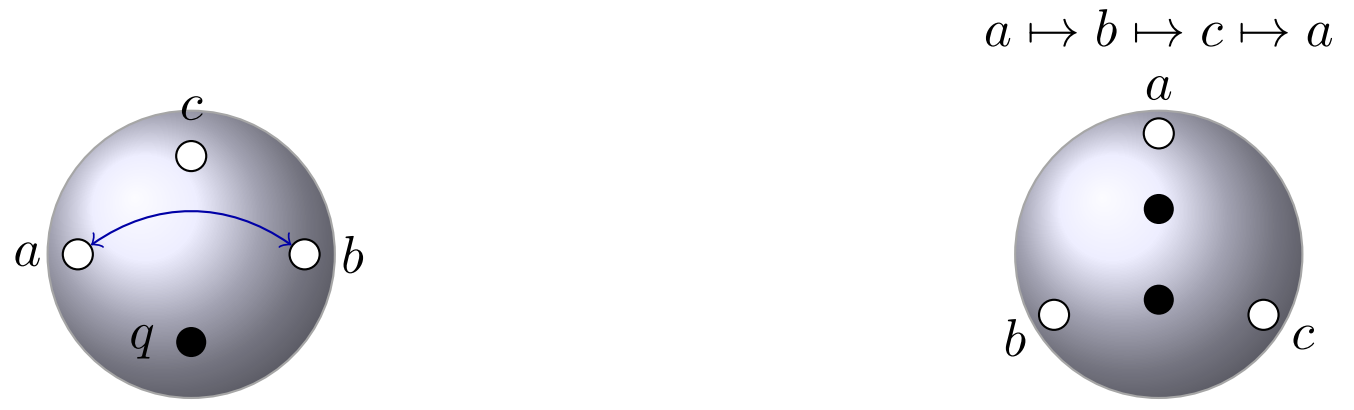
\begin{figure}[ht]
\centering
\begin{tikzpicture}[every node/.style={font=\small}]
\foreach \x in {-3.2,3.2}{
 \shade[ball color=blue!9] (\x,0) circle (.95);
 \draw[gray!70] (\x,0) circle (.95);}
\node[circle,draw,fill=white,inner sep=2pt,label=left:{$a$}] (a) at (-3.95,0) {};
\node[circle,draw,fill=white,inner sep=2pt,label=right:{$b$}] (b) at (-2.45,0) {};
\draw[<->,blue!65!black,bend left=35] (a) to (b);
\node[circle,draw,fill=white,inner sep=2pt,label=above:{$c$}] at (-3.2,.65) {};
\node[circle,fill=black,inner sep=2pt,label=left:{$q$}] at (-3.2,-.58) {};
\node[text width=5.5cm,align=center] at (-3.2,-1.8)
 {transposition $(ab)$\\one available point: $e^z$};
\foreach \angle/\lab in {90/a,210/b,330/c}{
 \node[circle,draw,fill=white,inner sep=2pt] at
 ({3.2+.8*cos(\angle)},{.8*sin(\angle)}) {};
 \node at ({3.2+1.1*cos(\angle)},{1.1*sin(\angle)}) {$\lab$};}
\node[circle,fill=black,inner sep=2pt] at (3.2,.3) {};
\node[circle,fill=black,inner sep=2pt] at (3.2,-.3) {};
\node at (3.2,1.5) {$a\mapsto b\mapsto c\mapsto a$};
\node[text width=5.5cm,align=center] at (3.2,-1.8)
 {three-cycle $(abc)$\\two available points: $e^{2z}$};
\end{tikzpicture}
\caption{Schematic fixed-point configurations on a rational component,
whose normalization is a sphere. Open circles are heavy branches;
black dots are fixed points available to light markings. For the
transposition, the other fixed point is the heavy branch $c$.
For the three-cycle, neither fixed point is heavy. Each light label
can choose any available point, with collisions allowed.}
\label{fig:rational-fixed-points}
\end{figure}

\smallskip
\noindent\emph{Rational vertex contributions.}
The rational vertex sectors have a similarly concrete interpretation.
At a fixed coarse moduli point, a permutation of the heavy branches
is realized by an automorphism of the component. Each light label must
lie at a fixed point of that automorphism, away from the heavy
branches. Figure~\ref{fig:rational-fixed-points} displays the two
nonidentity sectors of a trivalent rational vertex.

\noindent{\bf Trivalent rational vertices. } We have three possible permutations of heavy points, namely the identity, a transposition of two heavy points, a three-cycle moving the three heavy points.
\begin{itemize}
\item At a trivalent rational
vertex the identity contributes $e^{-z}$, by equation \eqref{eq:eulervertex-trivalent}. 
\item A transposition of two
branches is realized by the unique involution of $\mathbb P^1$
exchanging them and fixing the third. Of its two fixed points, one
is a heavy branch and one is available for light labels; hence its
vertex series is $e^z$. 
\item The automorphism realizing a three-cycle has two
fixed points, both outside the heavy markings.\footnote{An automorphism of $\mathbb P^1$ cyclically
permuting three distinct points has order $3$.
It is therefore conjugate in $\mathrm{PGL}_2(\mathbb C)$
to $z\mapsto\omega z$, where
$\omega=e^{2\pi i/3}$.
Consequently, it has exactly two fixed points,
neither of which belongs to the three-point orbit.}
Hence its vertex series is $e^{2z}$. 
\end{itemize}

\noindent{\bf Four-valent rational vertices.} We have four possible permutations of heavy points: the identity, 
a transposition of two heavy points, a double transposition of pairs of heavy points, and a four-cycle moving the four heavy points.
\begin{itemize}
\item At a four-valent rational vertex the identity contributes $-e^{-2z}$, by \eqref{eq:eulervertex-rational}.
\item For a single transposition of heavy branches, the fixed locus in
$M_{0,4}$ is one point; the two fixed points of the realizing
involution are the two individually fixed heavy branches. Thus no
light labels are available and the sector equals $1$.
\item For a double transposition, the contributing series is $-e^{2z}$,
as computed in Example \ref{ex:fixed-point-monodromy}. \item Finally, a four-cycle has
one fixed heavy configuration\footnote{Suppose that an automorphism $f\in\operatorname{Aut}(\mathbb P^1)$
cyclically permutes four distinct points $a,b,c,d$. Since $f^4$
fixes at least three points, it is the identity; moreover,
$f^2\neq\mathrm{id}$, so $f$ has order exactly $4$.
Consequently, $f$ is conjugate to $z\mapsto iz$ and hence has
exactly two distinct fixed points, which coincide with those of
the involution $f^2$. In particular, $f^2$ exchanges the pairs
$(a,c)$ and $(b,d)$. Such an automorphism exists if and only if
the quadruple is harmonic, namely
$[a,c;b,d]=-1$, where
$[a,c;b,d]:=\frac{(a-b)(c-d)}{(a-d)(c-b)}$.}, realized by an automorphism of order
four. Its two fixed points are outside the heavy branches, giving
the series $e^{2z}$.
\end{itemize}

For reference, the rational vertex series just computed are
\[
\begin{array}{c|ccc}
 d=3 & \mathrm{id} & (12) & (123)\\\hline
 \text{series} & e^{-z} & e^z & e^{2z}
\end{array}
\qquad
\begin{array}{c|cccc}
 d=4 & \mathrm{id} & (12) & (12)(34) & (1234)\\\hline
 \text{series} & -e^{-2z} & 1 & -e^{2z} & e^{2z}
\end{array}
\]
\smallskip
\noindent\emph{Moved vertices.}
When two genus-one vertices are exchanged, their light-label sets
are empty and their heavy-pointed curves must be isomorphic. Their
fixed locus is the diagonal copy of $M_{1,1}$, of Euler characteristic
$1$. Exchanged trivalent rational vertices likewise contribute $1$,
since $M_{0,3}$ is a point. 

\smallskip
This completes the calculation of all local building blocks
required by the seven genus-two stabilization graphs.
We now combine them, using Burnside's formula for the
automorphism groups of the graphs, to obtain their exact
contributions.

\subsection{The seven exact graph contributions}
\label{sec:g2-seven-graphs}
Denote the graphs by $\Gamma_2$ (a single genus-two vertex),
$\Gamma_{1\ell}$ (a genus-one vertex with a loop),
$\Gamma_{11}$ (two genus-one vertices joined by an edge),
$\Gamma_{\ell\ell}$ (a rational vertex with two loops),
$\Gamma_{\ell1}$ (a rational vertex with a loop joined to a
genus-one vertex), $\Gamma_{\mathrm{db}}$ (two rational looped
vertices joined by a bridge), and $\Theta$ (two rational vertices
joined by three edges). This notation names the seven drawings in
Figure~\ref{fig:stable-graphs-g2} without changing their half-edge
conventions.

Table~\ref{tab:g2-graph-groups} records the graph automorphism groups.
These groups preserve the vertex genera and include permutations of
half-edges, in particular the reversal of a loop. They are distinct
from the automorphism groups of the individual curves parameterized
by the vertex spaces.

\begin{table}[ht]
\centering\small
\begin{tabular}{@{}lclp{7cm}@{}}
\toprule
Graph & Group & Order & Generators / actions\\
\midrule
$\Gamma_2$ & $\{1\}$ & $1$ & No edges or nontrivial graph symmetries.\\
$\Gamma_{1\ell}$ & $C_2$ & $2$ & Reverse the loop and exchange its branches.\\
$\Gamma_{11}$ & $C_2$ & $2$ & Exchange the elliptic vertices and reverse the bridge.\\
$\Gamma_{\ell\ell}$ & $(C_2\times C_2)\rtimes C_2$ & $8$ & Reverse either loop; exchange the loops.\\
$\Gamma_{\ell1}$ & $C_2$ & $2$ & Reverse the rational loop.\\
$\Gamma_{\mathrm{db}}$ & $(C_2\times C_2)\rtimes C_2$ & $8$ & Reverse either loop; exchange the vertices and their loops.\\
$\Theta$ & $\mathfrak S_3\times C_2$ & $12$ & Permute the three edges; exchange the vertices.\\
\bottomrule
\end{tabular}
\caption{Automorphism groups of the seven genus-two stabilization graphs.
Here $C_2$ is the cyclic group of order two. In each semidirect
product the last $C_2$ exchanges the two factors; this is the
dihedral group of order eight.}
\label{tab:g2-graph-groups}
\end{table}

\begin{proposition}[Exact genus-two stabilization-graph decomposition]
\label{prop:g2-seven-graphs}
With $U,V,T,W,P,K$ defined above, the exact contributions are
\begin{align}
A^\chi_{2,\Gamma_2}&=T,\label{eq:g2-graph-genus2}\\
A^\chi_{2,\Gamma_{1\ell}}&=\tfrac12(VP+WK),
\label{eq:g2-graph-elliptic-loop}\\
A^\chi_{2,\Gamma_{11}}&=\tfrac12(U^2P+K),
\label{eq:g2-graph-two-elliptic}\\
A^\chi_{2,\Gamma_{\ell\ell}}&=
\tfrac18(-e^{-2z}P^2+2KP-e^{2z}K^2),
\label{eq:g2-graph-bouquet}\\
A^\chi_{2,\Gamma_{\ell1}}&=
\tfrac12UP(e^{-z}P+e^zK),
\label{eq:g2-graph-lollipop}\\
A^\chi_{2,\Gamma_{\mathrm{db}}}&=
\tfrac18(e^{-2z}P^3+2KP^2+e^{2z}K^2P+4K),
\label{eq:g2-graph-dumbbell}\\
A^\chi_{2,\Theta}&=
\tfrac1{12}(e^{-2z}P^3+3e^{2z}P+2e^{4z}+K^3+3K+2).
\label{eq:g2-graph-theta}
\end{align}
Their sum is $A_2^\chi$.
\end{proposition}

\begin{proof}
We apply the exact identity \eqref{eq:graph-burnside}, retaining all
sectors rather than just the maximal propagator degree. For each symmetry we multiply the vertex contributions and the edge
contributions listed in Section~\ref{sec:g2-graph-data}, then average
over the group in Table~\ref{tab:g2-graph-groups}. Three rules make this
assembly explicit: moved vertices carry no light labels; moved edges
have empty chains; and an edge preserved as a set contributes $P$ or
$K$ according to whether its half-edges are fixed or exchanged.

\emph{The three graphs involving only positive-genus vertices.}
The graph $\Gamma_2$ has no graph automorphisms and contributes $T$.
For $\Gamma_{1\ell}$ there are two symmetries. The identity leaves
the elliptic vertex unrestricted ($V$) and the loop fixed half-edgewise
($P$), giving $VP$. Loop reversal imposes the heavy-branch exchange
($W$) and reverses the inserted chain ($K$), giving $WK$. Thus its
Burnside average is $(VP+WK)/2$.

For $\Gamma_{11}$ the identity gives one factor $U$ at each vertex
and one factor $P$ at the bridge, hence $U^2P$.
The vertex exchange prohibits all light labels on the graph vertices;
its fixed parameter locus is the diagonal copy of $M_{1,1}$, of Euler
characteristic one. The bridge is reversed and contributes $K$.
Averaging the identity and exchange sectors gives $(U^2P+K)/2$.
This proves \eqref{eq:g2-graph-genus2}--\eqref{eq:g2-graph-two-elliptic}.

\emph{The bouquet $\Gamma_{\ell\ell}$.}
The eight symmetries split into two sets: four preserve each loop
as a set, and four exchange the loops. We first consider those
preserving the loops. The identity
gives $-e^{-2z}P^2$. The two single-loop reversals each give $KP$,
using the single-transposition vertex sector $1$, one reversed
edge $K$, and one half-edgewise fixed edge $P$. Reversing both loops
gives $-e^{2z}K^2$. Now consider the four elements exchanging the loops. Two
are double transpositions of the four branches and two are four-cycles.
The loops are moved and carry no inserted components. Their vertex
contributions $-2e^{2z}$ and $+2e^{2z}$ cancel. Averaging gives
\eqref{eq:g2-graph-bouquet}.

\emph{The mixed graph $\Gamma_{\ell1}$.}
The rational and elliptic vertices cannot be exchanged because
their genera differ. The only nonidentity symmetry reverses the
rational loop. The identity gives $e^{-z}UP^2$: the two vertices
contribute $e^{-z}$ and $U$, and the loop and bridge each contribute
$P$. Under loop reversal the rational vertex contribution becomes
$e^z$ and the loop contribution becomes $K$; the elliptic vertex
and the bridge still give $U$ and $P$. This sector is $e^zUKP$.
Their average is
\eqref{eq:g2-graph-lollipop}.

\emph{The dumbbell $\Gamma_{\mathrm{db}}$.}
As for the bouquet, four of the eight symmetries preserve the
vertices and four exchange them. With the vertices preserved, the
identity gives $e^{-2z}P^3$: both rational vertices have their identity
sector, and the two loops and bridge contribute $P$.
The two single-loop reversals give $KP^2$ each: the factors $e^z$
and $e^{-z}$ from the two rational vertices cancel.
Reversing both loops gives $e^{2z}K^2P$.
The four remaining elements exchange the vertices and the loops.
No graph vertex or moved loop may carry a light label; the vertex
parameter spaces are points. Only the reversed bridge contributes
labels, yielding $K$ for each of these four elements.
This proves \eqref{eq:g2-graph-dumbbell}.

\emph{The theta graph.}
Here $\Aut\Theta=\mathfrak S_3\times\mathbb Z/2\mathbb Z$.
We first keep both vertices fixed and classify the six edge
permutations. The identity gives $e^{-2z}P^3$. Each of the three
transpositions induces a branch transposition at both vertices,
giving $e^z e^z$; two edges are moved, while the remaining edge is
fixed half-edgewise and gives $P$. Its sector is therefore $e^{2z}P$.
Each of the two three-cycles gives $e^{2z}$ at both vertices and
moves all edges, so its sector is $e^{4z}$.

Next exchange the vertices, optionally composing with any of the
six edge permutations. No
light labels occur on the graph vertices. Exchange without an edge
permutation reverses all three edges and gives $K^3$.
For each of the three edge transpositions, two edges are moved
and the third is preserved with its ends exchanged. This gives
one factor $K$, hence the sector $K$ for each of these three symmetries. The two exchanges accompanied by a
three-cycle move every edge and give $1$ each. Averaging these twelve
terms gives \eqref{eq:g2-graph-theta}.

Every all-light genus-two curve has exactly one of these seven stable
graphs, so additivity gives the asserted sum.
\end{proof}

\begin{figure}[ht]
\centering
\begin{tikzpicture}[every node/.style={font=\small}]
\cv{bou}{(-3.5,0)}{0}
\draw[graph edge] (bou) edge[loop left] node[left] {$e_1$} (bou);
\draw[graph edge] (bou) edge[loop right] node[right] {$e_2$} (bou);
\draw[<->,blue!65!black] (-4.15,1.1)--(-2.85,1.1);
\node[above] at (-3.5,1.1) {exchange loops};
\node[text width=5.5cm,align=center] at (-3.5,-1.2)
 {bouquet: four heavy branches\\loop reversals also act on branches};
\cv{tha}{(2,0)}{0} \cv{thb}{(4.4,0)}{0}
\draw[graph edge,bend left=35] (tha) to node[above] {$e_1$} (thb);
\draw[graph edge] (tha)--node[above] {$e_2$}(thb);
\draw[graph edge,bend right=35] (tha) to node[below] {$e_3$}(thb);
\draw[<->,blue!65!black] (2,1.35)--(4.4,1.35);
\node[above] at (3.2,1.35) {exchange vertices};
\node[text width=5.5cm,align=center] at (3.2,-1.5)
 {theta: permute $e_1,e_2,e_3$\\and optionally exchange vertices};
\end{tikzpicture}
\caption{Generators of the bouquet and theta symmetries. In the
bouquet, reversing a loop exchanges its two incident half-edges.
For theta, edge permutations form $\mathfrak S_3$ and commute with
vertex exchange, giving $\Aut\Theta=\mathfrak S_3\times C_2$.
The blue arrows indicate symmetries, not edges.}
\label{fig:g2-symmetries}
\end{figure}
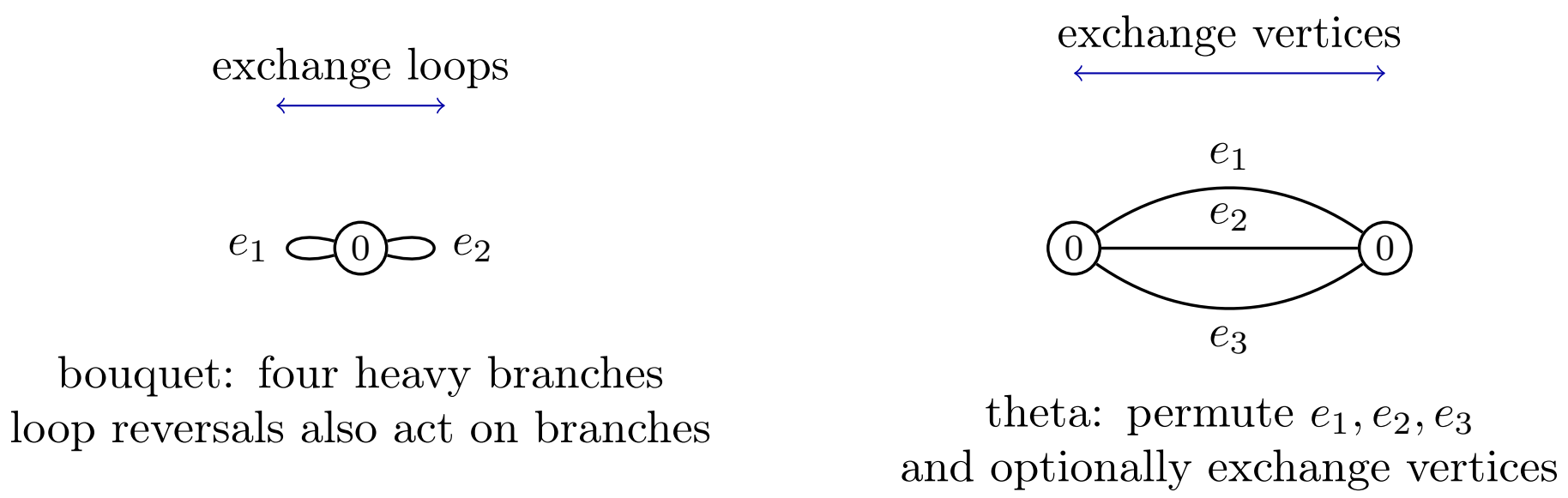

Tables~\ref{tab:bouquet-sectors} and
\ref{tab:theta-sectors} summarize the sectors of the two
most symmetric rational graphs;
Figure~\ref{fig:g2-symmetries} illustrates their basic symmetries.

\begin{table}[ht]
\centering\small
\begin{tabular}{@{}lcl@{}}
\toprule
Action on the bouquet & Number & Fixed-sector series\\
\midrule
Identity & $1$ & $-e^{-2z}P^2$\\
Reverse one loop & $2$ & $KP$\\
Reverse both loops & $1$ & $-e^{2z}K^2$\\
Exchange loops; branch double transposition & $2$ & $-e^{2z}$\\
Exchange loops; branch four-cycle & $2$ & $e^{2z}$\\
\bottomrule
\end{tabular}
\par\vspace{.6em}
\caption{The eight bouquet automorphisms. The last two rows cancel
in the Burnside average. Branch permutations act on the four heavy
points of the rational component of the stabilized curve.}
\label{tab:bouquet-sectors}
\end{table}

\begin{table}[ht]
\centering\small
\begin{tabular}{@{}llcl@{}}
\toprule
Graph vertices & Edge permutation & Number & Fixed-sector series\\
\midrule
Fixed & Identity & $1$ & $e^{-2z}P^3$\\
Fixed & Transposition & $3$ & $e^{2z}P$\\
Fixed & Three-cycle & $2$ & $e^{4z}$\\
Exchanged & Identity & $1$ & $K^3$\\
Exchanged & Transposition & $3$ & $K$\\
Exchanged & Three-cycle & $2$ & $1$\\
\bottomrule
\end{tabular}
\par\vspace{.6em}
\caption{The twelve theta automorphisms. With vertices exchanged,
no light labels can lie on graph vertices. An edge fixed as a set is
then reversed; moved edges have empty chains. Summing the last column
with its multiplicities and dividing by twelve recovers
\eqref{eq:g2-graph-theta}.}
\label{tab:theta-sectors}
\end{table}

\begin{remark}[The trivalent leading coefficient]
Only $\Gamma_{\mathrm{db}}$ and $\Theta$ have three edges. Their
identity sectors give
\[
\left(\frac18+\frac1{12}\right)e^{-2z}P^3
=\frac5{24}\frac{e^{-2z}}{(1-z)^3}.
\]
Thus the exact calculation already displays the trivalent mechanism
underlying the fixed-genus asymptotics of
Section~\ref{sec:fixed-genus}, with $\gamma_2=5/24$.
\qrem\end{remark}
\subsection{Closed Euler generating functions}

Summing all the contributions listed in Proposition~\ref{prop:g2-seven-graphs} gives an elementary
closed expression for the all-light series $A_2^\chi(z)$. The ordinary weight series $G_2^\chi(t)$ will subsequently follow by
composition with $\Psi$. We write it in the
expanded normalization used by Bini--Gaiffi--Polito in order to compare
the two constructions.

Set
\[
 D=e^E-1
\]
and define
\begin{align}
 N(D,E)={}&-2D^8(E-1)^2(7+3E)
 -24D^7(E-1)^2(17E-7) \notag\\
 &-3D^6(E-1)^2(201E+259)
 +30D^5(E-1)^2(61E-221) \notag\\
 &+15D^4(631E^3-2640E^2+3395E-1386) \notag\\
 &+60D^3(341E^3-1322E^2+1633E-652) \notag\\
 &+180D^2(138E^3-519E^2+635E-254) \notag\\
 &+360D(45E^3-167E^2+206E-84) \notag\\
 &+60(73E^3-270E^2+336E-144).
\label{eq:N2}
\end{align}

\begin{theorem}[Closed genus-two formulas]
\label{thm:g2closed}
Let
\[
 E=\Psi(t),\qquad D=e^E-1.
\]
Then
\begin{equation}
 G_2^{\chi}(t)
 =
 \frac{N(D,E)}
      {1440(1+D)^2(E-1)^3}
 +\frac{D^5}{4}.
\label{eq:K2closed}
\end{equation}
Consequently,
\begin{equation}
 A_2^{\chi}(z)
 =
 G_2^{\chi}(\Phi(z))
 =
 \frac{N(e^z-1,z)}
      {1440e^{2z}(z-1)^3}
 +\frac{(e^z-1)^5}{4}.
\label{eq:A2closed}
\end{equation}
\end{theorem}

\begin{proof}
Substitute \eqref{eq:g2-U}--\eqref{eq:g2-W} and \eqref{eq:g2-P-K}
into the seven formulas of Proposition~\ref{prop:g2-seven-graphs}.
Putting their sum over the common denominator
$1440e^{2z}(z-1)^3$ and collecting powers of $D=e^z-1$ gives
\eqref{eq:A2closed}. 
Theorem~\ref{thm:pde} gives $G_2^\chi=A_2^\chi\circ\Psi$.
Setting $E=\Psi(t)$ and $D=e^E-1$ therefore yields
\eqref{eq:K2closed}.
\end{proof}

\begin{remark}[The two Bini--Gaiffi--Polito displays and the separate sign correction]
\label{rem:g2-source-discrepancies}
The rational term in \eqref{eq:K2closed} is the concluding formula
of Section~7 in the arXiv version of Bini--Gaiffi--Polito
\cite{BGP} (arXiv:math/9806048v1, p.~27).
That version contains two different displays of the genus-two formula:
Theorem~1.2 on p.~2 has the first numerator term
$-2D^8(E-1)(7+3E)$, whereas the concluding display on p.~27 has
$-2D^8(E-1)^2(7+3E)$. Only the second display has the correct
factor $(E-1)^2$; the introductory display is missing the square.

This typographical discrepancy is separate from the sign error in
the graph-type-$5$ contribution identified by Bini--Harer
\cite[pp.~511--512]{BH}. Their correction replaces $-D^5/8$
by $+D^5/8$ in that contribution and therefore adds $D^5/4$
to the final generating function. This additional correction is
included in \eqref{eq:K2closed}.

The direct stabilization-graph computation recovers the
concluding display, including its factor $(E-1)^2$,
together with the Bini--Harer correction $D^5/4$.
\qrem\end{remark}

\subsection{Minimal-weight asymptotics}

We first consider the all-light Euler characteristics
\(
\chi\bigl(\Mb_{2,(0^+)^n}\bigr).
\)
The closed formula \eqref{eq:A2closed} shows that $A_2^\chi$ has a unique
finite singularity, namely a pole of order three at $z=1$.

Write
\begin{equation}
A_2^\chi(1-\delta)
=
\frac{\ell_{-3}}{\delta^3}
+\frac{\ell_{-2}}{\delta^2}
+\frac{\ell_{-1}}{\delta}
+O(1).
\label{eq:A2Laurent}
\end{equation}
A direct expansion of \eqref{eq:A2closed} gives
\[
\ell_{-3}
=
\frac{5}{24e^2},
\]
\[
\ell_{-2}
=
\frac{
2e^4+4e^3+9e^2+4e+6
}{
24e^2
},
\]
and
\[
\ell_{-1}
=
\frac{
4e^8+16e^7+52e^6+88e^5
+133e^4+112e^3+106e^2+40e+37
}{
288e^2
}.
\]

It is convenient to collect the corresponding factorial-scale
coefficients in the quadratic polynomial
\beq
\label{eq:P2poly}
P_2(n)=\alpha_2n^2+\alpha_1n+\alpha_0,
\eeq
where
\begin{equation}
\alpha_2:=\frac{\ell_{-3}}2,
\qquad
\alpha_1:=\frac32\ell_{-3}+\ell_{-2},
\qquad
\alpha_0:=\ell_{-3}+\ell_{-2}+\ell_{-1}.
\label{eq:alpha2def}
\end{equation}
Explicitly,
\[
\alpha_2
=
\frac{5}{48e^2},
\]
\[
\alpha_1
=
\frac{
4e^4+8e^3+18e^2+8e+27
}{
48e^2
},
\]
and
\[
\alpha_0
=
\frac{
4e^8+16e^7+52e^6+88e^5
+157e^4+160e^3+214e^2+88e+169
}{
288e^2
}.
\]Numerically,
\[
\alpha_2=0.0140974253371472,
\qquad
\alpha_1=1.58124098333656,
\qquad
\alpha_0=36.3166969595549.
\]

\begin{theorem}[All-light genus-two asymptotics]
\label{thm:g2light}
For every fixed $c>6$, as $n\to\infty$,
\begin{equation}
\chi\bigl(\Mb_{2,(0^+)^n}\bigr)
=
n!\,P_2(n)+O_c(c^n),
\label{eq:g2lightasym}
\end{equation}
where $P_2$ is defined in \eqref{eq:P2poly}-\eqref{eq:alpha2def}.
In particular,
\[
\chi\bigl(\Mb_{2,(0^+)^n}\bigr)
\sim
\frac{5}{48e^2}\,n!\,n^2.
\]
Thus the complete asymptotic expansion at factorial scale terminates
after three terms; the remaining contribution is exponentially smaller.
\end{theorem}

\begin{proof}
Subtract from $A_2^\chi$ the principal part of its Laurent expansion at
$z=1$:
\[
R_2(z)
:=
A_2^\chi(z)
-
\frac{\ell_{-3}}{(1-z)^3}
-
\frac{\ell_{-2}}{(1-z)^2}
-
\frac{\ell_{-1}}{1-z}.
\]
By \eqref{eq:A2closed}, the apparent singularity of $R_2$ at $z=1$ is
removable, and $R_2$ is entire.

The singular part contributes exactly
\[
[z^n](1-z)^{-3}
=
\binom{n+2}{2},
\qquad
[z^n](1-z)^{-2}
=
n+1,
\qquad
[z^n](1-z)^{-1}
=
1.
\]
Hence
\begin{align*}
n![z^n]\bigl(A_2^\chi(z)-R_2(z)\bigr)
&=
n!
\left[
\ell_{-3}\binom{n+2}{2}
+\ell_{-2}(n+1)
+\ell_{-1}
\right]\\
&=
n!\left(
\alpha_2n^2+\alpha_1n+\alpha_0
\right).
\end{align*}

It remains to estimate the entire remainder. From
\eqref{eq:A2closed}, after the factor $e^{-2z}$ is taken into account,
the largest exponential occurring in
$N(e^z-1,z)$ is $e^{6z}$; the correction $(e^z-1)^5/4$ has smaller
exponential type. Dividing by $(z-1)^3$ after cancellation of the
principal part does not change the exponential type. Thus $R_2$ is
entire of exponential type at most $6$.

Standard Cauchy estimates therefore give, for every $c>6$,
\[
n![z^n]R_2(z)=O_c(c^n),
\]
which proves \eqref{eq:g2lightasym}.
\end{proof}

\begin{example}
For $n=50$ one has
\[
\chi\bigl(\Mb_{2,(0^+)^{50}}\bigr)
=
4581040958455114586283032342264655965929154747770239568893313089874.
\]
On the other hand, the three factorial-scale terms of
\eqref{eq:g2lightasym} give
\[
50!\,P_2(50)
=
4.581040958455114586283032342269517\ldots\times10^{66}.
\]
Thus the approximation agrees with the exact value through the first
$30$ significant digits; the first differing digit is the next one.
The relative error is approximately
\(
1.06\times10^{-30}.
\)
This illustrates quantitatively the fact that, after the finite
factorial-scale expansion, the remaining contribution is only
exponentially large.
\qetr\end{example}

\begin{remark}
To our knowledge, the closed formula \eqref{eq:A2closed} for the
all-light generating function, and hence Theorem~\ref{thm:g2light}, do
not appear previously in the literature. Here they follow directly from the exact stabilization-graph decomposition.
The comparison with the corrected Bini--Gaiffi--Polito expression
provides an independent consistency check.
\qrem\end{remark}

\subsection{Ordinary-weight asymptotics}

We now turn to the complete asymptotic expansion of
$\chi(\Mb_{2,n})$.

As in genus one, set
\[
\rho=e-2,\qquad
\zeta=1-\frac{t}{\rho},\qquad
a=\sqrt{\frac{2\rho}{e}}.
\]
Recall the analytic germ $\mathfrak d$ defined by
\eqref{eq:def-d}. By the local inversion formula
\eqref{eq:delta-d},
\[
1-\Psi(t)
=
\mathfrak d(a\sqrt{\zeta}).
\]

Recall the Laurent expansion at $z=1$
\begin{equation}
A_2^\chi(1-\delta)
=
\sum_{j\ge-3}\ell_j\delta^j,
\label{eq:A2LaurentRecall}
\end{equation}
whose principal part was computed above. Define the coefficients $B_j$
by the convergent Laurent--Puiseux expansion
\begin{equation}
A_2^\chi\bigl(1-\mathfrak d(as)\bigr)
=
\sum_{j\ge-3}B_js^j.
\label{eq:Bj2}
\end{equation}
Since the coefficients of $\mathfrak d$ are recursively computable,
\eqref{eq:Bj2} gives an effective algorithm for all $B_j$.

Recall the coefficients $C_r(\alpha)$ defined by
\eqref{eq:Cr-gamma}.
Finally, define
\begin{equation}
\kappa_0:=\frac{B_{-3}}{\Gamma(3/2)}
=\frac{2B_{-3}}{\sqrt{\pi}},
\qquad
\kappa_1:=B_{-2},
\label{eq:kappa01}
\end{equation}
and, for every $m\ge0$,
\begin{equation}
\kappa_{m+2}
:=
\sum_{k=-2}^{m-1}
\frac{
B_{2k+1}\,
C_{m-k-1}\!\left(k+\frac12\right)
}{
\Gamma\!\left(-k-\frac12\right)
}.
\label{eq:kappageneral}
\end{equation}

\begin{theorem}[Full genus-two expansion]
\label{thm:g2full}
The radius of convergence of $G_2^\chi$ is
\(
\rho=e-2,
\)
and $\rho$ is its unique dominant singularity. In a dented neighborhood of $\rho$, one has the convergent
Laurent--Puiseux expansion
\beq
G_2^\chi(t)
=
\sum_{j\ge-3}B_j
\left(1-\frac{t}{\rho}\right)^{j/2},
\label{eq:G2LaurentPuiseaux}
\eeq
where the coefficients $B_j$ are defined by \eqref{eq:Bj2}.
As $n\to\infty$,
\begin{equation}
\chi(\Mb_{2,n})
\sim
n!\rho^{-n}
\left(
\kappa_0n^{1/2}
+\kappa_1
+\sum_{m\ge0}\kappa_{m+2}n^{-m-1/2}
\right),
\label{eq:g2fullasym}
\end{equation}
where the coefficients $\kappa_j$ are defined by
\eqref{eq:kappa01}--\eqref{eq:kappageneral}.

In particular, apart from the constant term $\kappa_1$, no negative
integer powers of $n$ occur.
\end{theorem}

\begin{proof}By the analysis of $\Psi$ in Section~\ref{sec:sing1},
$\Psi$ is analytic for $|t|<\rho$ and admits continuation
to a dented neighborhood of its unique dominant singularity
$\rho=e-2$. Since $A_2^\chi$ has its only finite singularity
at $z=1$, with a nonzero cubic pole, the composition
$G_2^\chi=A_2^\chi\circ\Psi$ has the same unique dominant
singularity at $t=\rho$.

Substituting
$1-\Psi(t)=\mathfrak d(a\sqrt{\zeta})$
into \eqref{eq:Bj2} gives the stated local expansion \eqref{eq:G2LaurentPuiseaux}.

For $\alpha\notin\mathbb Z_{\ge0}$, coefficient transfer gives
\[
[t^n]\left(1-\frac{t}{\rho}\right)^\alpha
=
\rho^{-n}
\frac{\Gamma(n-\alpha)}
{\Gamma(-\alpha)\Gamma(n+1)}.
\]
The nonnegative even powers of $\zeta$ are analytic and therefore do
not contribute to the algebraic asymptotic expansion. Among the
negative even powers only $\zeta^{-1}$ occurs, and its contribution is
exactly
\(
B_{-2}\rho^{-n}.
\)
All the remaining singular terms have half-integer exponent.

Applying \eqref{eq:Cr-gamma} to each odd power
$\zeta^{k+1/2}$ and collecting the coefficient of
$n^{-m-1/2}$ gives precisely \eqref{eq:kappageneral}.
The leading term $\zeta^{-3/2}$ gives
$\kappa_0=2B_{-3}/\sqrt{\pi}$, while $\zeta^{-1}$ gives
$\kappa_1=B_{-2}$. This proves \eqref{eq:g2fullasym}.
\end{proof}

The first coefficients $B_j$ are already determined by the principal
part of \eqref{eq:A2LaurentRecall}:
\begin{align}
B_{-3}
&=
a^{-3}\ell_{-3},
\label{eq:Bm3}\\
B_{-2}
&=
a^{-2}(\ell_{-2}-\ell_{-3}),
\label{eq:Bm2}\\
B_{-1}
&=
a^{-1}
\left(
\ell_{-1}
-\frac23\ell_{-2}
+\frac5{24}\ell_{-3}
\right).
\label{eq:Bm1}
\end{align}
Consequently,
\[
\kappa_0
=
\frac{5}
{24\sqrt{2\pi e}\,(e-2)^{3/2}},
\]
whereas
\[
\kappa_1
=
\frac{
2e^4+4e^3+9e^2+4e+1
}{
48e(e-2)
}.
\]

The next coefficients are
\begin{align*}
\kappa_2
&=
\frac1{\sqrt{\pi}}
\left(
\frac34 B_{-3}+B_{-1}
\right),\\
\kappa_3
&=
\frac1{\sqrt{\pi}}
\left(
-\frac7{64}B_{-3}
-\frac18B_{-1}
-\frac12B_1
\right),\\
\kappa_4
&=
\frac1{\sqrt{\pi}}
\left(
\frac9{512}B_{-3}
+\frac1{128}B_{-1}
-\frac3{16}B_1
+\frac34B_3
\right).
\end{align*}

There is also a simple arithmetic structure behind these coefficients.
Since $\ell_j\in\mathbb Q(e)$ and the Taylor coefficients of
$\mathfrak d$ are rational, \eqref{eq:Bj2} implies
\[
B_j\in a^j\mathbb Q(e).
\]
Since $a^2=2(e-2)/e\in\mathbb Q(e)$, it follows in particular that
\[
B_{2m}\in\mathbb Q(e),
\qquad
B_{2m+1}\in a\,\mathbb Q(e),
\]
and hence
\[
\kappa_j=\kappa_0R_j(e),
\qquad
R_j(e)\in\mathbb Q(e),
\qquad j\ge2.
\]
For example,
\[
\frac{\kappa_2}{\kappa_0}
=
\frac38+
\frac{a^2}{2\ell_{-3}}
\left(
\ell_{-1}
-\frac23\ell_{-2}
+\frac5{24}\ell_{-3}
\right).
\]

Numerically,
\[
\kappa_0=0.0828093295004206,
\qquad
\kappa_1=2.85866535466901,
\]
and
\[
\kappa_2=26.2080174505291,
\qquad
\kappa_3=-51.1775027401697,
\qquad
\kappa_4=32.0204982386165.
\]

\begin{example}
For $n=50$ one has
\[
\chi(\Mb_{2,50})
=
3265182350471382337037628271422579740566802244399315434767713752283351612.
\]
Truncating \eqref{eq:g2fullasym} after the term
$\kappa_4n^{-5/2}$ gives
\[
\chi(\Mb_{2,50})
\approx
3.26515057624\times10^{72},
\]
with relative error approximately
\(
9.73\times10^{-6}.
\)
For comparison, truncation after
$\kappa_2n^{-1/2}$ and $\kappa_3n^{-3/2}$ gives relative errors
approximately
\(
2.04\times10^{-2}\)
and \(
2.68\times10^{-4},
\)
respectively.
\qetr\end{example}

\section{Euler-characteristic asymptotics in fixed genus}
\label{sec:fixed-genus}

We now use the stabilization-graph formalism of Section~\ref{sec:graph-formalism}
to extract the leading asymptotic behavior for every fixed $g\ge2$.
The genus-one and genus-two calculations have exhibited the two
mechanisms: unbounded rational chains produce propagator poles, and
composition with $\Psi$ converts these poles into Puiseux singularities.
The genus-zero comparison was given in Section~\ref{sec:g0comparison}.

The argument below extends the trivalent mechanism already
observed in genus two: maximal-edge graphs determine the
leading pole, while nonidentity sectors have lower pole order.

\subsection{The leading pole}
\label{sec:leadingpole}

Let $\mathfrak T_g$ denote the finite set of connected rational
trivalent graphs of first Betti number $g$, allowing loops and multiple
edges, and define
\begin{equation}
\gamma_g
:=
\sum_{\Gamma\in\mathfrak T_g}
\frac1{|\Aut\Gamma|}.
\label{eq:gammagdef}
\end{equation}

\begin{theorem}[All-light asymptotics in fixed genus]
\label{thm:alllightfixedg}
For every fixed $g\ge2$, the all-light Euler EGF has a pole of order
exactly $3g-3$ at $z=1$, and
\begin{equation}
A_g^\chi(z)
=
\frac{\gamma_g\,e^{-(2g-2)z}}
{(1-z)^{3g-3}}
+
O\!\left((1-z)^{-3g+4}\right).
\label{eq:alllightpole}
\end{equation}
In particular,
\begin{equation}
\chi\bigl(\Mb_{g,(0^+)^n}\bigr)
\sim
\frac{\gamma_g e^{-2g+2}}
{\Gamma(3g-3)}
\,n!\,n^{3g-4}.
\label{eq:alllightfixedgasym}
\end{equation}
\end{theorem}

\begin{proof}
For fixed $g$, there are only finitely many stabilization graphs. By
Lemma~\ref{lem:maxedges}, every graph has at most $3g-3$ edges, and
equality occurs precisely for rational trivalent graphs.

For such a graph,
\[
h_v=0,\qquad d_v=3
\]
at every vertex. By \eqref{eq:eulervertex-trivalent},
\[
V^\chi_{0,3}(z)=e^{-z}.
\]
Since a rational trivalent graph of genus $g$ has $2g-2$ vertices and
$3g-3$ edges, its identity sector in \eqref{eq:identitysector}
contributes
\[
\frac{e^{-(2g-2)z}}
{|\Aut\Gamma|(1-z)^{3g-3}}.
\]
All nonidentity sectors have pole order at most $3g-4$, by
Proposition~\ref{prop:graphcontribution}. Moreover, every nontrivalent
graph has at most $3g-4$ edges, and its vertex factors are entire by
Lemma~\ref{lem:eulervertex}. Hence such graphs also contribute poles of
order at most $3g-4$. Summing the leading contributions over the
rational trivalent graphs gives \eqref{eq:alllightpole}.

Finally,
\[
[z^n](1-z)^{-3g+3}
=
\frac{\Gamma(n+3g-3)}
{\Gamma(3g-3)\Gamma(n+1)}
\sim
\frac{n^{3g-4}}{\Gamma(3g-3)}.
\]
Since $e^{-(2g-2)z}$ is analytic at $z=1$, standard singularity
analysis yields
\[
[z^n]A_g^\chi(z)
\sim
\frac{\gamma_g e^{-2g+2}}
{\Gamma(3g-3)}
\,n^{3g-4}.
\]
Multiplying by $n!$ gives \eqref{eq:alllightfixedgasym}.
\end{proof}

\begin{proposition}[Complete all-light asymptotics in fixed genus]
\label{prop:alllightfullfixedg}
For every fixed $g\geq2$, there exist a unique polynomial
$P_g\in\mathbb R[n]$ and a constant $b>0$ such that
\[
 \chi\bigl(\Mb_{g,(0^+)^n}\bigr)=n!P_g(n)+O(b^n)
 \qquad(n\to\infty).
\]
The polynomial $P_g$ has degree $3g-4$ and leading coefficient
\[
 \frac{\gamma_g e^{-2g+2}}{\Gamma(3g-3)}.
\]
\end{proposition}

\begin{proof}
By Proposition~\ref{prop:graphcontribution} and the finiteness of the
set of stabilization graphs, $A_g^\chi$ is a finite sum of terms
$p(z)e^{\lambda z}(1-z)^{-r}$, where $p$ is a polynomial and
$0\leq r\leq3g-3$. Write its principal part at $z=1$ as
\[
 \sum_{r=1}^{3g-3}\frac{\ell_r}{(1-z)^r},
\]
and define
\[
 P_g(n)=\sum_{r=1}^{3g-3}\ell_r\binom{n+r-1}{r-1}.
\]
The coefficients $\ell_r$ are real. By
Theorem~\ref{thm:alllightfixedg},
\[
 \ell_{3g-3}=\gamma_g e^{-2g+2}>0.
\]
Consequently, $P_g$ has degree $3g-4$ and the stated leading
coefficient.

The function
\[
 R_g(z)=A_g^\chi(z)-\sum_{r=1}^{3g-3}\frac{\ell_r}{(1-z)^r}
\]
is entire. Its finite exponential-rational representation implies
that it is of exponential type. Indeed, outside a fixed disk centered
at the origin, each term is bounded by a polynomial in $|z|$ times
$e^{|\lambda||z|}$; on the remaining compact set, $R_g$ is bounded.
Thus, for some $C,\tau>0$,
\[
 |R_g(z)|\leq C e^{\tau|z|},\qquad z\in\mathbb C.
\]
Cauchy's coefficient estimate on the circle of radius $n/\tau$,
followed by Stirling's formula, gives
\[
 n!\,\bigl|[z^n]R_g(z)\bigr|
 \leq Cn!\left(\frac{e\tau}{n}\right)^n
 =O\bigl(\sqrt n\,\tau^n\bigr)=O(b^n)
\]
for any fixed $b>\tau$. Since
$[z^n](1-z)^{-r}=\binom{n+r-1}{r-1}$, this proves the expansion.

Finally, if two polynomials satisfied such an estimate, their
difference would satisfy $n!\,|P(n)|=O(b^n)$ for some $b>0$.
A nonzero polynomial cannot satisfy this bound, so $P_g$ is unique.
\end{proof}

\subsection{Ordinary-weight asymptotics}
\label{sec:ordinaryfixedg}

We now return to ordinary markings. Recall
\[
G_g^\chi(t)=A_g^\chi(\Psi(t)),
\qquad
\rho=e-2,
\qquad
a=\sqrt{\frac{2\rho}{e}},
\]
and, with
\[
\zeta=1-\frac{t}{\rho},
\qquad
1-\Psi(t)
=
\mathfrak d(a\sqrt{\zeta})
=
a\sqrt{\zeta}+O(\zeta).
\]

Substitution into \eqref{eq:alllightpole} gives
\begin{equation}
G_g^\chi(t)
=
\gamma_g e^{-2g+2}
a^{-3g+3}
\left(1-\frac{t}{\rho}\right)^{-(3g-3)/2}
+
O\!\left(
\left(1-\frac{t}{\rho}\right)^{-(3g-4)/2}
\right).
\label{eq:ordinaryleadinglocal}
\end{equation}

\begin{theorem}[Ordinary Euler asymptotics in fixed genus]
\label{thm:ordinaryfixedg}
For every fixed $g\ge2$, as $n\to\infty$,
\begin{equation}
\chi(\Mb_{g,n})
\sim
\frac{\gamma_g}
{\Gamma\!\left(\frac{3g-3}{2}\right)
e^{(g-1)/2}(2\rho)^{(3g-3)/2}}
\,n!\rho^{-n}n^{(3g-5)/2}.
\label{eq:ordinaryfixedgasym}
\end{equation}
\end{theorem}

\begin{proof}
The function $\Psi$ has its unique dominant singularity at
$t=\rho$, and admits analytic continuation to a dented neighborhood
of that point. Equation \eqref{eq:ordinaryleadinglocal} therefore
falls directly under the standard transfer theorem \cite{FO,FS}. Since
\[
[t^n]
\left(1-\frac{t}{\rho}\right)^{-(3g-3)/2}
\sim
\frac{\rho^{-n}n^{(3g-5)/2}}
{\Gamma((3g-3)/2)},
\]
the leading coefficient is
$\gamma_g e^{-2g+2}a^{-3g+3}/\Gamma((3g-3)/2)$.
Using $a^2=2\rho/e$ gives \eqref{eq:ordinaryfixedgasym}.
\end{proof}

\begin{remark}[Complete expansions]
\label{rem:fullfixedg}
The preceding argument also gives a complete asymptotic expansion.
Indeed, $A_g^\chi$ is a finite sum of entire exponential polynomials
times powers of $(1-z)^{-1}$, and therefore has a convergent Laurent
expansion at $z=1$. Substituting
\[
1-\Psi(t)=\mathfrak d(a\sqrt\zeta)
\]
produces a convergent Laurent--Puiseux expansion in powers of
$\zeta^{1/2}$. Termwise coefficient transfer then yields the complete
fixed-genus expansion of $\chi(\Mb_{g,n})$.

As already seen in genera one and two, integer and half-integer powers
play different roles: negative integer powers of $\zeta$ give finite
polynomial contributions in $n$, whereas nonintegral powers produce
infinite descending asymptotic series. No uniformity in $g$ is asserted
here.
\qrem\end{remark}

\begin{remark}[The fixed-genus exponent]
Combining \eqref{eq:g0leading} with the genus-one result and
Theorem~\ref{thm:ordinaryfixedg}, one obtains the common exponent law
\(
n^{(3g-5)/2}
\)
in genera $g=0$ and $g\ge2$. Genus one is the logarithmic marginal
case: its leading singularity is logarithmic and gives
\[
\chi(\Mb_{1,n})
\sim
\frac14\,n!\rho^{-n}n^{-1}.
\]
Thus the same dominant scale $\rho=e-2$ governs every fixed genus,
while the nature of the singularity changes with $g$.
\qrem\end{remark}

Kannan, Serpente, and Yun asked in \cite[Question~5.3]{KSY} for the
asymptotic behavior of
\[
\frac{
\chi\bigl(\Mb_{g,(0^+)^n}\bigr)
}{
\chi(\Mb_{g,n})
},
\qquad g\ge2,
\]
and, in particular, whether this ratio tends to zero as
$n\to\infty$. The preceding results give a precise answer.

\begin{corollary}
\label{cor:KSYquestion}
For every fixed $g\ge2$,
\begin{equation}
\frac{
\chi\bigl(\Mb_{g,(0^+)^n}\bigr)
}{
\chi(\Mb_{g,n})
}
\sim
\left(\frac{2(e-2)}{e}\right)^{(3g-3)/2}
\frac{\Gamma\!\left(\frac{3g-3}{2}\right)}{\Gamma(3g-3)}
\,(e-2)^n n^{(3g-3)/2}.
\label{eq:KSYratio}
\end{equation}
In particular,
\[
\lim_{n\to\infty}
\frac{
\chi\bigl(\Mb_{g,(0^+)^n}\bigr)
}{
\chi(\Mb_{g,n})
}
=0.
\]
\end{corollary}

\begin{proof}
Taking the quotient of \eqref{eq:alllightfixedgasym} and
\eqref{eq:ordinaryfixedgasym}, the common factor
$\gamma_g e^{-2g+2}$ cancels, and we obtain
\[
\frac{
\chi\bigl(\Mb_{g,(0^+)^n}\bigr)
}{
\chi(\Mb_{g,n})
}
\sim
a^{3g-3}
\frac{
\Gamma\!\left(\frac{3g-3}{2}\right)
}{
\Gamma(3g-3)
}
\rho^n n^{(3g-3)/2},
\]
with $\rho=e-2$. Since $0<e-2<1$, the limit follows.
\end{proof}

\section{The trivalent constants and cubic field theory}
\label{sec:gammag}

We now study the constants
\[
\gamma_g
=
\sum_{\Gamma\in\mathfrak T_g}
\frac1{|\Aut\Gamma|},
\qquad g\ge2,
\]
introduced in \eqref{eq:gammagdef}, where $\mathfrak T_g$ denotes the set of
connected rational trivalent graphs of first Betti number $g$.

Although they arise here from the geometry of stabilization graphs, these
constants are classical combinatorial quantities. Up to a shift of
indices, they coincide with the total symmetry weights of connected
vacuum diagrams in zero-dimensional cubic, or $\phi^3$, field theory;
see \cite{CLP,BenderMilton}. The same constants also occur in the
asymptotic enumeration of connected graphs of fixed excess; see
\cite{FSS}. Our purpose in this section is to give a self-contained
derivation of their generating function, recurrence, and large-genus
behavior, making their relation with trivalent graphs explicit.

The argument is entirely combinatorial. In particular, no familiarity
with Feynman rules will be assumed: the trivalent graphs will emerge
directly from Gaussian pairings. For general background on diagrammatic expansions
and their combinatorics, see
\cite[Chapter~1]{MarcolliFM} and
\cite[Chapter~7]{ItzyksonDrouffeII}.

\subsection{Gaussian pairings and trivalent graphs}
\label{sec:gaussiantrivalent}

Consider the formal Gaussian integral
\begin{equation}
Z(s)
:=
\frac1{\sqrt{2\pi}}
\int_{\mathbb R}
\exp\!\left(
-\frac{x^2}{2}
+\frac{\sqrt{s}}6x^3
\right)\,dx .
\label{eq:formalZ}
\end{equation}
The integral is understood formally: we expand the cubic interaction in
powers of $\sqrt{s}$ and integrate coefficientwise against the Gaussian
measure. In other words, if $X$ is a standard real Gaussian random
variable, then
\[
Z(s)
=
\mathbb E\!\left[
\exp\!\left(\frac{\sqrt{s}}6X^3\right)
\right]
\]
as a formal power series.

Expanding the exponential gives
\[
Z(s)
=
\sum_{k\ge0}
\frac{s^{k/2}}{6^k k!}\,
\mathbb E(X^{3k})
=
\sum_{k\ge0}
\frac{s^{k/2}}{6^k k!}
\frac1{\sqrt{2\pi}}
\int_{\mathbb R}x^{3k}e^{-x^2/2}\,dx .
\]
Since the odd Gaussian moments vanish, only $k=2m$ contribute. Using
\(
\mathbb E(X^{6m})=(6m-1)!!,
\)
we obtain
\begin{equation}
Z(s)
=
\sum_{m\ge0}
\frac{(6m-1)!!}{(2m)!\,6^{2m}}\,s^m.
\label{eq:Zseries}
\end{equation}
Here, as usual, $(-1)!!=1$.

The identity $\mathbb E(X^{2N})=(2N-1)!!$ has a simple combinatorial
meaning: $(2N-1)!!$ is the number of pairings of $2N$ objects. This is
the elementary Gaussian form of Wick's pairing rule \cite{Wick}.

We now apply this interpretation to \eqref{eq:Zseries}. Fix $m\ge1$ and
consider $2m$ labeled vertices, each equipped with three labeled
half-edges. There are altogether $6m$ half-edges, and pairing them
produces a trivalent graph with $2m$ vertices; loops and multiple edges
are allowed. Conversely, every trivalent graph equipped with labels on
its vertices and on the three half-edges incident to each vertex arises
from such a pairing.

Let $\mathcal P_m$ denote the set of pairings of these $6m$ half-edges.
Then
\(
|\mathcal P_m|=(6m-1)!!.
\)
The group
\(
\mathfrak S_3^{\,2m}\rtimes\mathfrak S_{2m}
\)
acts on $\mathcal P_m$ by permuting the three half-edges at each vertex
and by permuting the vertices, and has order
\(
(3!)^{2m}(2m)!=6^{2m}(2m)!.
\)
Its orbits are precisely the isomorphism classes of trivalent graphs
with $2m$ vertices. Moreover, the stabilizer of a pairing representing
$\Gamma$ is naturally identified with $\Aut\Gamma$. The
orbit--stabilizer formula therefore gives
\begin{equation}
\frac{(6m-1)!!}{(2m)!\,6^{2m}}
=
\sum_{\substack{\Gamma\ {\rm trivalent}\\
V(\Gamma)=2m}}
\frac1{|\Aut\Gamma|},
\label{eq:alltrivalent}
\end{equation}
where disconnected graphs are included in the sum. Figure~\ref{fig:connected-exponential} illustrates the
factorial arising from repeated isomorphic components
in the case $\Theta\sqcup\Theta\sqcup\Gamma_{\mathrm{db}}$.

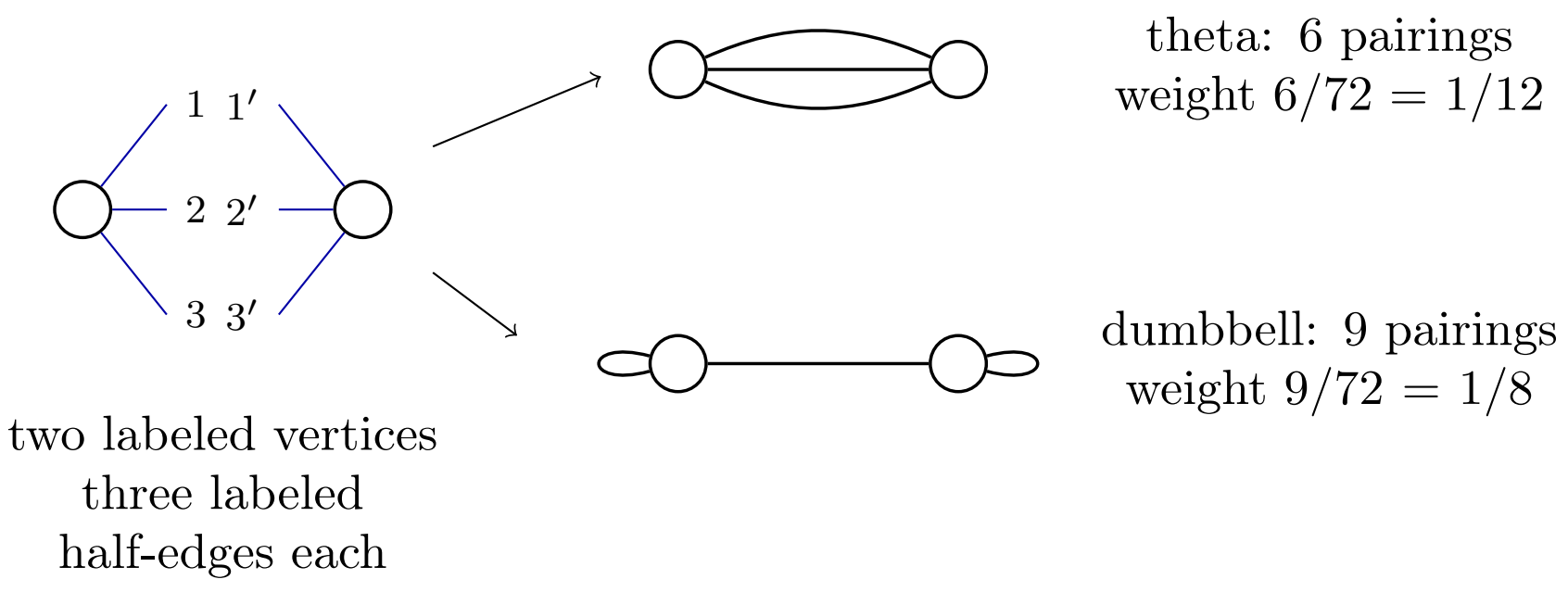
\begin{figure}[ht]
\centering
\begin{tikzpicture}[every node/.style={font=\small}]
\cv{paira}{(-4,0)}{} \cv{pairb}{(-2,0)}{}
\foreach \y/\lab in {.75/1,0/2,-.75/3}{
\draw[blue!65!black] (paira)--(-3.4,\y);
\node[right,font=\scriptsize] at (-3.4,\y) {$\lab$};
\draw[blue!65!black] (pairb)--(-2.6,\y);
\node[left,font=\scriptsize] at (-2.6,\y) {$\lab'$};}
\node[text width=3.5cm,align=center] at (-3,-2.05)
 {two labeled vertices\\three labeled half-edges each};
\draw[->] (-1.5,.45)--(-.3,.95);
\draw[->] (-1.5,-.45)--(-.9,-.9);
\cv{pairth1}{(.25,1)}{} \cv{pairth2}{(2.25,1)}{}
\tripleedge{pairth1}{pairth2}
\node[text width=3.3cm,align=center] at (4.9,1)
 {theta: $6$ pairings\\weight $6/72=1/12$};
\cv{pairdb1}{(.25,-1.1)}{} \cv{pairdb2}{(2.25,-1.1)}{}
\draw[graph edge] (pairdb1)--(pairdb2);
\draw[graph edge] (pairdb1) edge[loop left] (pairdb1);
\draw[graph edge] (pairdb2) edge[loop right] (pairdb2);
\node[text width=3.3cm,align=center] at (4.9,-1.1)
 {dumbbell: $9$ pairings\\weight $9/72=1/8$};
\end{tikzpicture}
\caption{Pairings of six half-edges and the two genus-two trivalent
graphs. All pairings here produce connected graphs. Relabeling vertices
and their half-edges has group order $72$; its orbit stabilizers have
orders $12$ and $8$. Thus $\gamma_2=1/12+1/8=5/24$.}
\label{fig:pairings-genus-two}
\end{figure}

By equations \eqref{eq:Zseries} and \eqref{eq:alltrivalent}, we have
\begin{equation}
Z(s)
=
\sum_{\Gamma\ {\rm trivalent}}
\frac{s^{V(\Gamma)/2}}{|\Aut\Gamma|}.
\label{eq:Zgraphs}
\end{equation}
Thus $Z(s)$ is the generating function of all trivalent graphs, weighted
by the reciprocal of their automorphism groups.

\begin{remark}[Zero-dimensional cubic field theory]
The expression \eqref{eq:formalZ} is the partition function of a
zero-dimensional cubic field theory. The integral itself is not
convergent for nonzero real $\sqrt{s}$ and is used here only as a
compact notation for its formal perturbative expansion. The preceding
argument is precisely its Wick expansion, written in elementary
combinatorial terms. Zero-dimensional field theories as generating
functions for weighted Feynman diagrams are discussed, for example, in
\cite{CLP,Borinsky}.
\qrem\end{remark}

\subsection{Connected graphs and the constants \texorpdfstring{$\gamma_g$}{gamma_g}}
\label{sec:connectedgammag}

Every trivalent graph is uniquely a disjoint union of connected trivalent graphs. The standard exponential formula, in its symmetry-weighted form, therefore implies that the generating function of connected trivalent graphs is the logarithm of the generating function of all trivalent graphs; see, for instance, \cite[Chapter~II]{FS}\cite[Chapter 1]{MarcolliFM}. Hence
\begin{equation}
\log Z(s)
=
\sum_{\substack{\Gamma\ {\rm connected}\\
\Gamma\ {\rm trivalent}}}
\frac{s^{V(\Gamma)/2}}{|\Aut\Gamma|}.
\label{eq:logZconnected}
\end{equation}

In the present weighted setting this identity can also be seen directly.
Indeed, if a connected graph $\Gamma_i$ occurs with multiplicity $m_i$
in a disconnected graph, the corresponding automorphisms contribute
the wreath-product factor
\(
\Aut(\Gamma_i)^{m_i}\rtimes\mathfrak S_{m_i},
\)
and hence the weight
\[
\frac1{m_i!\,|\Aut\Gamma_i|^{m_i}}.
\]
Multiplying over the isomorphism classes of connected components and
summing over their multiplicities gives exactly the exponential
relation above.

\begin{figure}[ht]
\centering
\begin{tikzpicture}[every node/.style={font=\small}]
\cv{compa}{(-4.8,0)}{} \cv{compb}{(-3.4,0)}{}
\tripleedge{compa}{compb}
\cv{compc}{(-1.6,0)}{} \cv{compd}{(-.2,0)}{}
\tripleedge{compc}{compd}
\cv{compe}{(2,0)}{} \cv{compf}{(3.4,0)}{}
\draw[graph edge] (compe)--(compf);
\draw[graph edge] (compe) edge[loop left] (compe);
\draw[graph edge] (compf) edge[loop right] (compf);
\draw[<->,blue!65!black] (-4.1,1)--(-.9,1);
\node[above] at (-2.5,1) {exchange identical components};
\node at (-4.1,-.75) {$\Theta$};
\node at (-.9,-.75) {$\Theta$};
\node at (2.7,-.75) {$\Gamma_{\mathrm{db}}$};
\end{tikzpicture}
\caption{A disconnected graph with two identical components.
Internal automorphisms contribute $12^2\cdot8$; exchanging the two
theta components contributes the additional factor $2!$.
In general, multiplicity $r$ of a connected graph contributes
$1/(r!\,|\Aut\Gamma|^r)$, giving the exponential formula.}
\label{fig:connected-exponential}
\end{figure}
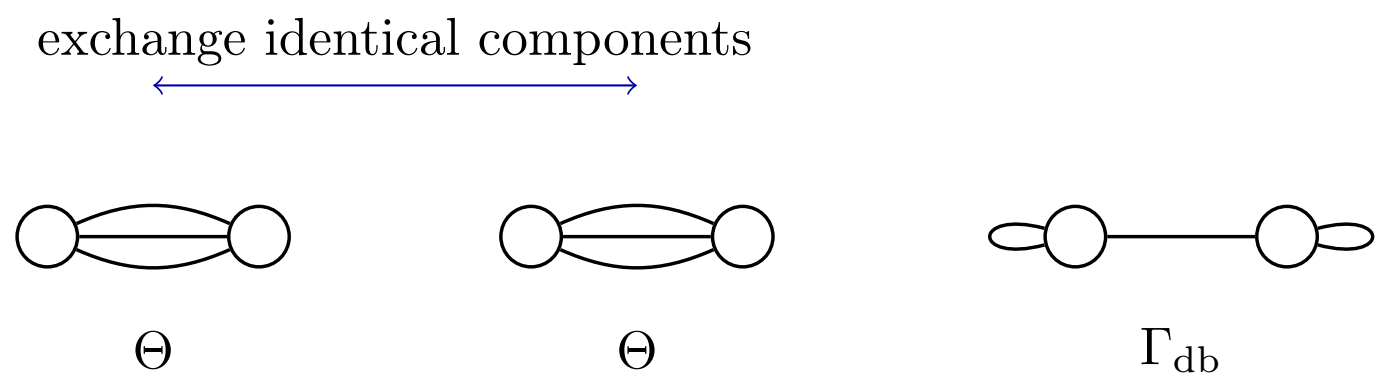

\begin{example}Figure~\ref{fig:connected-exponential} illustrates why repeated
isomorphic components introduce a factorial. If a graph is
$\Theta\sqcup\Theta\sqcup\Gamma_{\mathrm{db}}$, the two theta
components can be exchanged. Its symmetry weight is therefore
\[
\frac1{2!\,|\Aut\Theta|^2\,|\Aut\Gamma_{\mathrm{db}}|}
=\frac1{2!\,12^2\,8}.
\]
The factorial is exactly the one produced by the exponential series.
Thus taking the logarithm of $Z$ retains one connected component.
\qetr
\end{example}

For a connected trivalent graph,
\[
3V(\Gamma)=2E(\Gamma),
\]
and therefore
\[
b_1(\Gamma)
=
E(\Gamma)-V(\Gamma)+1
=
\frac{V(\Gamma)}2+1.
\]
Thus a connected trivalent graph of genus $g$ has exactly
$V(\Gamma)=2g-2$ vertices. Combining this observation with
\eqref{eq:logZconnected} yields the following identity.

\begin{theorem}[Cubic partition function]
\label{thm:cubicpartition}
As an identity of formal power series,
\begin{equation}
\log Z(s)
=
\sum_{g\ge2}\gamma_gs^{g-1}.
\label{eq:logZgammag}
\end{equation}
Moreover,
\begin{equation}
Z(s)
=
\sum_{m\ge0}
\frac{(6m-1)!!}{(2m)!\,6^{2m}}s^m
=
{}_2F_0\!\left(
\frac16,\frac56;;
\frac32s
\right).
\label{eq:Zhypergeom}
\end{equation}
\end{theorem}

\begin{proof}
Only the hypergeometric expression remains to be checked. Using the
Pochhammer symbol,
\[
\frac{(6m-1)!!}{(2m)!\,6^{2m}}
=
\frac{(1/6)_m(5/6)_m}{m!}
\left(\frac32\right)^m,
\]
which gives \eqref{eq:Zhypergeom}.
\end{proof}

The first values are
\begin{equation}
\gamma_2=\frac5{24},
\qquad
\gamma_3=\frac5{16},
\qquad
\gamma_4=\frac{1105}{1152},
\qquad
\gamma_5=\frac{565}{128},
\qquad
\gamma_6=\frac{82825}{3072}.
\label{eq:firstgammag}
\end{equation}
These values agree with the classical symmetry-weighted enumeration of
connected vacuum graphs in cubic field theory; in the notation of
\cite{BenderMilton},
\(
\gamma_g=C^{[3]}_{2g-2}.
\)
The corresponding rational sequence is recorded in the OEIS through
its numerators and denominators as
\href{https://oeis.org/A226260}{A226260} and
\href{https://oeis.org/A226261}{A226261}, respectively.
The same constants also occur as the leading kernel constants in the
enumeration of connected graphs of fixed excess \cite{FSS}.

\subsection{A recurrence for \(\gamma_g\)}
\label{sec:gammagrecurrence}

The hypergeometric representation gives a particularly simple
recurrence. Related differential and recursive descriptions of
zero-dimensional diagram counts already appear in \cite{CLP}; we
include the short derivation adapted to our normalization.

Let
\[
\vartheta:=s\frac{d}{ds}.
\]
The series $Z$ satisfies the hypergeometric differential equation
\begin{equation}
\vartheta Z
=
\frac32s
\left(\vartheta+\frac16\right)
\left(\vartheta+\frac56\right)Z.
\label{eq:Zode}
\end{equation}
Set
\[
F(s):=\log Z(s)
=
\sum_{g\ge2}\gamma_gs^{g-1},
\qquad
Y:=\vartheta F.
\]
Since
\[
\frac{\vartheta^2Z}{Z}
=
\vartheta Y+Y^2,
\]
equation \eqref{eq:Zode} becomes
\begin{equation}
Y
=
\frac32s
\left(
\vartheta Y+Y^2+Y+\frac5{36}
\right).
\label{eq:Yriccati}
\end{equation}

Comparing coefficients gives the following recurrence.

\begin{proposition}[Recurrence for the trivalent constants]
\label{prop:gammagrecurrence}
One has $\gamma_2=5/24$, and for every $g\ge3$,
\begin{equation}
\gamma_g
=
\frac32(g-2)\gamma_{g-1}
+
\frac{3}{2(g-1)}
\sum_{h=2}^{g-2}
(h-1)(g-1-h)\gamma_h\gamma_{g-h}.
\label{eq:gammagrecurrence}
\end{equation}
\end{proposition}

\begin{proof}
Write
\(
Y(s)
=
\sum_{g\ge2}(g-1)\gamma_gs^{g-1}.
\)
The coefficient of $s$ in \eqref{eq:Yriccati} gives
\(
\gamma_2=\frac32\frac5{36}=\frac5{24}.
\)
For $g\ge3$, comparison of the coefficient of $s^{g-1}$ gives
\[
(g-1)\gamma_g
=
\frac32
\left[
(g-1)(g-2)\gamma_{g-1}
+
\sum_{h=2}^{g-2}
(h-1)(g-1-h)\gamma_h\gamma_{g-h}
\right],
\]
which is equivalent to \eqref{eq:gammagrecurrence}.
\end{proof}

\subsection{Large-genus behavior}
\label{sec:gammaglargegenus}

The large-order behavior of zero-dimensional Feynman-diagram
enumeration is classical; see \cite{CLP}, and \cite{Borinsky} for
complete asymptotic expansions and a modern systematic treatment. For
completeness, we give a short direct derivation of the leading term in
our normalization.

\begin{theorem}[Large-genus asymptotics]
\label{thm:gammaglargegenus}
As $g\to\infty$,
\begin{equation}
\gamma_g
\sim
\frac1{2\pi}
\left(\frac32\right)^{g-1}
(g-2)!.
\label{eq:gammaglargegenus}
\end{equation}
\end{theorem}

\begin{proof}
Write
\(
Z(s)=1+\sum_{m\ge1}a_ms^m.
\)
By \eqref{eq:Zhypergeom},
\[
a_m
=
\frac{(1/6)_m(5/6)_m}{m!}
\left(\frac32\right)^m.
\]
Using Euler's reflection formula
\(
\Gamma\!\left(\frac16\right)
\Gamma\!\left(\frac56\right)
=
2\pi,
\)
we obtain
\begin{equation}
a_m
=
\frac1{2\pi}
\left(\frac32\right)^m
\frac{
\Gamma(m+1/6)\Gamma(m+5/6)
}{
\Gamma(m+1)
}.
\label{eq:amgamma}
\end{equation}
The standard gamma-ratio expansion therefore gives
\begin{equation}
a_m
=
\frac1{2\pi}
\left(\frac32\right)^m
\Gamma(m)
\left(1+O(m^{-1})\right).
\label{eq:amasym}
\end{equation}

It remains to compare the coefficients of $Z$ and $\log Z$. Put
\[
F(s)=\log Z(s)=\sum_{m\ge1}f_ms^m,
\qquad
f_m=\gamma_{m+1}.
\]
From $Z'=F'Z$,
\[
ma_m
=
\sum_{k=1}^m
kf_ka_{m-k},
\qquad
a_0=1,
\]
and hence
\begin{equation}
f_m
=
a_m
-
\frac1m
\sum_{k=1}^{m-1}
kf_ka_{m-k}.
\label{eq:fmam}
\end{equation}

Since connected graphs form a subset of all graphs and the symmetry weights are positive, \(0\leq f_k\leq a_k\).
Moreover,
\eqref{eq:amasym} implies that there exist constants $C_1,C_2>0$
such that, for every $m\ge1$,
\[
C_1\left(\frac32\right)^m\Gamma(m)
\le
a_m
\le
C_2\left(\frac32\right)^m\Gamma(m).
\]
Hence, for $1\le k\le m-1$,
\[
\frac{a_ka_{m-k}}{a_m}
\le
C\,
\frac{\Gamma(k)\Gamma(m-k)}{\Gamma(m)}
\]
for some constant $C>0$ independent of $k$ and $m$. Therefore
\begin{equation}
\frac1{a_m}
\sum_{k=1}^{m-1}k f_k a_{m-k}
\le
\frac1{a_m}
\sum_{k=1}^{m-1}k a_k a_{m-k}
\le
C\sum_{k=1}^{m-1}
k\,\frac{\Gamma(k)\Gamma(m-k)}{\Gamma(m)}=
C\sum_{k=1}^{m-1}
\binom{m-1}{k}^{-1}.
\end{equation}
For $m\ge4$, the last sum satisfies
\[
\sum_{k=1}^{m-1}\binom{m-1}{k}^{-1}
=
1+O(m^{-1}).
\]
Indeed, the term $k=m-1$ equals $1$, the terms $k=1$ and
$k=m-2$ are both $(m-1)^{-1}$, while for
$2\le k\le m-3$ one has
\(
\binom{m-1}{k}
\ge
\binom{m-1}{2},
\)
so that the sum of the remaining terms is $O(m^{-1})$. Hence
\[
\sum_{k=1}^{m-1}k f_k a_{m-k}
=
O(a_m).
\]
Equation \eqref{eq:fmam} consequently gives
\[
f_m
=
a_m\left(1+O(m^{-1})\right).
\]
Combining this with \eqref{eq:amasym} and setting $m=g-1$, we obtain
\[
\gamma_g=f_{g-1}
\sim
\frac1{2\pi}
\left(\frac32\right)^{g-1}
\Gamma(g-1)
=
\frac1{2\pi}
\left(\frac32\right)^{g-1}
(g-2)!,
\]
which proves \eqref{eq:gammaglargegenus}.
\end{proof}

\begin{remark}
The series $Z(s)$ has zero radius of convergence, since
$a_{m+1}/a_m\sim 3m/2$. This factorial divergence is natural from the
viewpoint of perturbative zero-dimensional field theory. For the
purposes of the present paper, however, only the coefficientwise
combinatorial interpretation of $Z$ is needed. Much finer all-order
asymptotic information is available; see \cite{Borinsky}.
\qrem\end{remark}

\appendix
\section{An explicit genus-three calculation}\label{app:genus3}
This appendix applies the same fixed-graph calculation as in genus two
and records an exact formula for both genus-three Euler series.
There are $42$ stabilization graphs, with respectively
$1,2,5,9,12,8,5$ graphs having $0,1,\ldots,6$ edges, see Figure \ref{fig:stable-graphs-g3}.
The formula is organized by propagator degree; the graph-by-graph
calculation is recorded below.

\subsection{Local inputs and the finite enumeration}
For a partition $\mu$ of $d$, let $\pi_\mu$ be a permutation of cycle
type $\mu$ of the $d$ heavy markings, fixing every light label.
With $x=e^z$, write
$B_{h,\mu}(x):=\sum_{n\ge0}\chi(\mathcal V_{h,d,n}^{\pi_\mu})z^n/n!$;
for the empty partition, $B_{3,\varnothing}=V^\chi_{3,0}$.
The collision decomposition gives
\begin{equation}\label{eq:g3-local-trace}
B_{h,\mu}(x)=\sum_{k\ge0}
\operatorname{Tr}\!\left(\pi_\mu\mid
\chi_{\mathfrak S_{d+k}}(M_{h,d+k})\right)
\frac{(x-1)^k}{k!}.
\end{equation}
Here the equivariant Euler characteristic is the alternating sum of
compactly supported rational cohomology representations, and
$\pi_\mu$ has cycle type $\mu\cup1^k$ on the $d+k$ distinct markings.
The equality of this trace and the fixed-locus Euler characteristic
is the finite-order Lefschetz identity used in
Lemma~\ref{lem:burnside}.
For $\mu\ne1^d$, the sum stops at $k=2$ in genus zero and at
$k=2h+2$ in positive genus: a nonidentity automorphism has at most
that many fixed points. For $\mu=1^d$, use the vertex formula of
Lemma~\ref{lem:eulervertex}.

The equivariant smooth-space inputs are the tables in
\cite{BergstromData}, specialized to ordinary Euler characteristics.
Only $M_{1,N}$ and $M_{2,N}$ with $N\le8$, and the ordinary
characteristics of $M_{3,N}$ with $N\le8$, are required.
For $N=0,\ldots,8$ the latter values are, in order,
$3,6,8,6,4,-10,30,-660,6540$.
For $N\ge9$, the ordinary characteristic equals the
Harer--Zagier value $(-1)^N(N+3)!/6048$.
In particular, writing $Z=B_{3,\varnothing}$, we obtain
\[
Z=\frac{-3x^{12}+112x^9-84x^8+1440x^7+3192x^6+1200x^5+189x^4+2}{2016x^4}.
\]
All data used here are within the unconditional ranges described
in \cite{BergstromData}.

For the remaining positive-genus identity sectors, write
$U_d=B_{1,1^d}$ ($1\le d\le4$) and $T_d=B_{2,1^d}$ ($d=1,2$).
They can be recovered from the genus-two inputs by
$U_d=\partial_x^{d-1}U$ and $T_d=\partial_x^dT$, where the
$U$ and $T$ of Section~\ref{sec:g2-graph-data} are viewed as functions
of $x=e^z$. For elliptic nonidentity sectors write $E_\mu=B_{1,\mu}$,
and set $J=B_{2,(2)}$. The required values are
\begin{center}
\small
\begin{tabular}{@{}ll@{\qquad}ll@{}}
$E_{(2)}$ & $(x^4+6x^2+8x-3)/12$ & $E_{(2,1)}$ & $(x^3+3x+2)/3$\\
$E_{(2,1,1)}$ & $x^2+1$ & $E_{(2,2)}$ & $(1-x^4)/2$\\
$E_{(3)}$ & $(-x^3+3x-2)/3$ & $E_{(4)}$ & $-x^2-1$\\
$J$ & $(x^6-15x^2+24x-10)/30$ & &
\end{tabular}
\end{center}
The rational identity sectors are
$B_{0,1^d}=(-1)^{d-3}(d-3)!x^{-(d-2)}$.
Table~\ref{tab:g3-rational-sectors} records every rational nonidentity
sector used by a genus-three graph, including return permutations at
moved vertices.
\begin{table}[ht]
\centering\small
\begin{tabular}{@{}cl@{\qquad}cl@{\qquad}cl@{}}
\toprule
$\mu$ & $B_{0,\mu}$ & $\mu$ & $B_{0,\mu}$ & $\mu$ & $B_{0,\mu}$\\\midrule
$(2,1)$ & $x$ & $(2,2,1,1)$ & $-2$ & $(3,3)$ & $- 3 x^{2}$\\
$(2,1,1)$ & $1$ & $(2,2,2)$ & $4 x^{2}$ & $(4)$ & $x^{2}$\\
$(2,1,1,1)$ & $0$ & $(3)$ & $x^{2}$ & $(4,1)$ & $2 x$\\
$(2,1,1,1,1)$ & $0$ & $(3,1)$ & $2 x$ & $(4,1,1)$ & $2$\\
$(2,2)$ & $- x^{2}$ & $(3,1,1)$ & $2$ & $(4,2)$ & $0$\\
$(2,2,1)$ & $- 2 x$ & $(3,2)$ & $0$ & $(6)$ & $x^{2}$\\
\bottomrule\end{tabular}
\caption{Rational branch-permutation sectors required in genus three.}
\label{tab:g3-rational-sectors}
\end{table}

For completeness, the enumeration uses $1\le v\le4$ vertices,
nonnegative edge multiplicities $m_{ij}$, and genus weights $h_i$ satisfying
\[
\sum_i h_i+\sum_{i\le j}m_{ij}-v+1=3,\qquad
2h_i-2+2m_{ii}+\sum_{j\ne i}m_{ij}>0.
\]
The underlying graph must be connected, and arrays related by a
simultaneous vertex relabeling are identified. These conditions give
precisely the $42$ rows of Table~\ref{tab:g3-graphs}.
Let $\Pi_\Gamma$ be the vertex-permutation group preserving the weights
and edge multiplicities. Then
\[
\begin{aligned}
\Aut\Gamma&\cong
\left(\prod_i(C_2\wr\mathfrak S_{m_{ii}})
\times\prod_{i<j}\mathfrak S_{m_{ij}}\right)\rtimes\Pi_\Gamma,\\
|\Aut\Gamma|&=|\Pi_\Gamma|
\prod_i2^{m_{ii}}m_{ii}!\prod_{i<j}m_{ij}!.
\end{aligned}
\]
Here $\wr$ denotes the \emph{wreath product}:
$C_2\wr\mathfrak S_m=C_2^{\,m}\rtimes\mathfrak S_m$, where
$\mathfrak S_m$ permutes the $m$ factors of $C_2^{\,m}$.
For $m$ loops at one vertex, the $m$ factors record their independent
reversals, and $\mathfrak S_m$ records their permutations.
Factors with zero multiplicity are trivial.

For each graph, enumerate the vertex permutations preserving the genus
weights and edge multiplicities, the permutations of parallel edges,
and the permutations and independent reversals of loops.
For a fixed vertex use $B_{h,\mu}(x)$. For a vertex cycle of length
$\ell>1$, no light labels are allowed; its factor is
$B_{h,\mu}(1)$, with $\mu$ the return permutation induced by
$\sigma^\ell$. An edge fixed half-edgewise contributes $P=(1-z)^{-1}$,
a reversed edge contributes $K=(1+x^2)/2$, and a moved edge contributes
$1$. Multiplication followed by the Burnside average gives the
graph contribution.

\subsection{Closed formulas}
Write $N_r(x)=\sum_{k=0}^{12}c_{k,r}x^k$, with the integer
coefficients in Table~\ref{tab:g3-numerators}; blank entries are zero.
\begin{table}[ht]
\centering\small\setlength{\tabcolsep}{5pt}
\begin{tabular}{@{}rrrrrrrr@{}}
\toprule
$k\backslash r$ & $0$ & $1$ & $2$ & $3$ & $4$ & $5$ & $6$\\\midrule
12 & 31 & -6 & 120 & 8 &  &  & \\
11 &  & -12 & 640 & 48 &  &  & \\
10 & -294 & 163 & 2380 & 204 &  &  & \\
9 & -2800 & 776 & 5440 & 544 &  &  & \\
8 & 12390 & 2094 & 11454 & 1890 & 4 &  & \\
7 & 57024 & 3372 & 17040 & 4140 & 16 &  & \\
6 & 120050 & 4998 & 28065 & 8085 & 52 &  & \\
5 & 63984 & 3876 & 29208 & 8820 & 88 &  & \\
4 & 46935 & 1549 & 24435 & 9441 & 167 & 10 & \\
3 &  & -64 & 4416 & 3028 & 120 & 20 & \\
2 &  & 9 & -505 & 15 & -11 & 45 & \\
1 &  & 4 & -520 & 1308 & -72 & 20 & \\
0 & 40 & -19 & 1307 & -2539 & 140 & -55 & 5\\
\midrule
$d_r$ & 40320 & 1440 & 17280 & 10368 & 288 & 96 & 16\\
\bottomrule\end{tabular}
\caption{Coefficients of the seven genus-three numerator polynomials
and their normalizing denominators.}
\label{tab:g3-numerators}
\end{table}
\begin{theorem}[Exact genus-three Euler series]\label{thm:g3-exact}
With $d_r$ as in Table~\ref{tab:g3-numerators} and $E=\Psi(t)$,
\begin{align}
A_3^\chi(z)&=e^{-4z}\sum_{r=0}^6
\frac{N_r(e^z)}{d_r(1-z)^r},\label{eq:g3-A-exact}\\
G_3^\chi(t)&=e^{-4E}\sum_{r=0}^6
\frac{N_r(e^E)}{d_r(1-E)^r}.\label{eq:g3-G-exact}
\end{align}
\end{theorem}
\begin{proof}
The $42$ graph contributions in Table~\ref{tab:g3-graphs} account
for all stabilization graphs of genus three. Applying
\eqref{eq:graph-burnside} with the local inputs above and collecting
the powers of $P$ gives \eqref{eq:g3-A-exact}.
The heavy/light composition identity gives \eqref{eq:g3-G-exact}.
\end{proof}

The five rational trivalent graphs have automorphism orders
$8,16,16,24,48$. Consequently $\gamma_3=5/16$, agreeing both with
$N_6/d_6$ and with the recurrence of Proposition~\ref{prop:gammagrecurrence}.

\subsection{Asymptotics to arbitrary order}
Set $q_r(z)=e^{-4z}N_r(e^z)/d_r$, and define
\begin{equation}\label{eq:g3-L-def}
L_r:=\sum_{s=r}^6\frac{(-1)^{s-r}}{(s-r)!}
q_s^{(s-r)}(1),\qquad 1\le r\le6.
\end{equation}
Thus the principal part at $z=1$ is
$\sum_{r=1}^6L_r(1-z)^{-r}$.
For example, $L_6=5/(16e^4)$ and
$L_5=5(2e^4+4e^3+9e^2+4e+13)/(96e^4)$; all six coefficients are
explicitly determined by \eqref{eq:g3-L-def}.
\begin{theorem}[Genus-three all-light asymptotics]\label{thm:g3-light-asym}
As $n\to\infty$,
\begin{equation}\label{eq:g3-light-full}
\chi(\Mb_{3,(0^+)^n})
=n!\sum_{r=1}^6L_r\binom{n+r-1}{r-1}+O(8^n).
\end{equation}
In particular,
$\chi(\Mb_{3,(0^+)^n})\sim e^{-4}n!n^5/384$.
\end{theorem}
\begin{proof}
Subtract the principal part from \eqref{eq:g3-A-exact}.
Each summand is an exponential divided by a power of $1-z$, with
its pole removed; its remaining coefficients are bounded by
$O(8^n/n!)$, since all exponential exponents lie between $-4$ and $8$.
This follows directly from the integral form of the exponential
Taylor remainder at $z=1$. The principal part has the coefficients
in \eqref{eq:g3-light-full}, and $L_6/5!=e^{-4}/384$.
\end{proof}

For the ordinary series set $\rho=e-2$, $a=\sqrt{2\rho/e}$ and define
$B_j^{(3)}$ by
$A_3^\chi(1-\mathfrak d(as))=\sum_{j\ge-6}B_j^{(3)}s^j$,
where $\mathfrak d$ is defined by \eqref{eq:def-d}.
This convergent expansion determines all the coefficients explicitly
from \eqref{eq:g3-A-exact}. Its first three principal coefficients are
\[
\begin{aligned}
B_{-6}^{(3)}&=a^{-6}L_6,&
B_{-5}^{(3)}&=a^{-5}(L_5-2L_6),\\
B_{-4}^{(3)}&=a^{-4}\left(L_4-\tfrac53L_5+\tfrac{17}{12}L_6\right).
\end{aligned}
\]
Set $p_2=B_{-6}^{(3)}/2$,
$p_1=3B_{-6}^{(3)}/2+B_{-4}^{(3)}$ and
$p_0=B_{-6}^{(3)}+B_{-4}^{(3)}+B_{-2}^{(3)}$.
Using the coefficients $C_r(\alpha)$ of \eqref{eq:Cr-gamma}, define
\begin{equation}\label{eq:g3-beta}
\beta_m=\sum_{k=0}^m
\frac{B_{2k-5}^{(3)}C_{m-k}(k-5/2)}{\Gamma(5/2-k)}
\quad(m\ge0).
\end{equation}
\begin{theorem}[Full genus-three ordinary-weight expansion]
\label{thm:g3-ordinary-asym}
The ordinary Euler series has radius of convergence $\rho=e-2$,
with $\rho$ its unique dominant singularity. As $n\to\infty$,
\begin{equation}\label{eq:g3-ordinary-full}
\chi(\Mb_{3,n})\sim n!\rho^{-n}
\left(p_2n^2+p_1n+p_0+
\sum_{m\ge0}\beta_m n^{3/2-m}\right),
\end{equation}
where the terms are understood in decreasing order of powers of $n$.
In particular, $\chi(\Mb_{3,n})\sim5n!\rho^{-n}n^2/(256e\rho^3)$.
There are no negative integer powers of $n$ in the full expansion.
\end{theorem}
\begin{proof}
The local inversion $1-\Psi(t)=\mathfrak d(a\sqrt{1-t/\rho})$
turns the Laurent expansion into
$\sum_{j\ge-6}B_j^{(3)}(1-t/\rho)^{j/2}$.
The terms with $j=-6,-4,-2$ transfer exactly to
$B_{-6}^{(3)}\binom{n+2}{2}+B_{-4}^{(3)}(n+1)+B_{-2}^{(3)}$,
which gives the polynomial above. The odd $j$ yield
\eqref{eq:g3-beta} by the gamma-ratio expansion; the nonnegative
even $j$ are analytic and do not contribute algebraic terms.
The continuation and unique dominant singularity follow as in
Theorem~\ref{thm:g2full}.
\end{proof}

\subsection{Checks and the graph-by-graph calculation}
The same enumeration reconstructs the existing genus-two closed
formula exactly. In genus three, the ordinary coefficients through
$n=8$ agree with the independently tabulated Euler characteristics
of $\Mb_{3,n}$ in \cite{BergstromData}:
\begin{table}[ht]
\centering\small
\begin{tabular}{@{}rrr@{\qquad}rrr@{}}\toprule
$n$ & Light & Ordinary & $n$ & Light & Ordinary\\\midrule
0 & 32 & 32 & 5 & 125128 & 156094\\
1 & 102 & 102 & 6 & 1073977 & 1526677\\
2 & 454 & 454 & 7 & 10282332 & 17004830\\
3 & 2510 & 2612 & 8 & 108636992 & 212619148\\
4 & 16495 & 18515 &  &  & \\
\bottomrule\end{tabular}
\caption{Initial coefficients of the exact genus-three Euler series.}
\label{tab:g3-coefficients}
\end{table}

In Table~\ref{tab:g3-graphs}, vertices are numbered $1,\ldots,v$;
$\mathbf h$ lists their genus weights, and $ij$ denotes an edge between
$i$ and $j$, with $ii$ a loop. A superscript denotes multiplicity.
The final column is the Burnside numerator
$S_\Gamma=\sum_{\sigma\in\Aut\Gamma}F_{\Gamma,\sigma}$, so the graph
contribution is $S_\Gamma/|\Aut\Gamma|$.
Here $x=e^z$, $P=(1-z)^{-1}$, $K=(1+x^2)/2$, and
$U_d,T_d,E_\mu,J,Z$ are the local series defined above.
The table enumerates $416$ graph automorphisms in total.

\begingroup
\small\setlength{\tabcolsep}{3pt}
\begin{longtable}{@{}rlp{2.5cm}c>{\raggedright\arraybackslash}p{9.35cm}@{}}
\caption{The $42$ genus-three graphs and their exact Burnside numerators.}
\label{tab:g3-graphs}\\
\toprule
No. & $\mathbf h$ & Edges & $|\Aut\Gamma|$ & $S_\Gamma$\\\midrule
\endfirsthead
\toprule No. & $\mathbf h$ & Edges & $|\Aut\Gamma|$ & $S_\Gamma$\\\midrule
\endhead
\bottomrule\endfoot
1 & $(3)$ & $\varnothing$ & $1$ & $Z$\\[2pt]
2 & $(2)$ & $11$ & $2$ & $J K + P T_{2}$\\[2pt]
3 & $(1,2)$ & $12$ & $1$ & $P T_{1} U_{1}$\\[2pt]
4 & $(1)$ & $11^{2}$ & $8$ & $2 E_{(2,1,1)} K P + E_{(2,2)} K^{2} + 2 E_{(2,2)} + 2 E_{(4)} + P^{2} U_{4}$\\[2pt]
5 & $(0,2)$ & $11,\allowbreak 12$ & $2$ & $K P T_{1} x + \frac{P^{2} T_{1}}{x}$\\[2pt]
6 & $(1,1)$ & $11,\allowbreak 12$ & $2$ & $E_{(2,1)} K P U_{1} + P^{2} U_{1} U_{3}$\\[2pt]
7 & $(1,1)$ & $12^{2}$ & $4$ & $E_{(2)}^{2} + K^{2} + P^{2} U_{2}^{2} + 1$\\[2pt]
8 & $(1,1,1)$ & $12,\allowbreak 13$ & $2$ & $E_{(2)} + P^{2} U_{1}^{2} U_{2}$\\[2pt]
9 & $(0)$ & $11^{3}$ & $48$ & $4 K^{3} x^{2} - 6 K^{2} P + 24 K x^{2} - \frac{6 P^{3}}{x^{4}} - 16 x^{2}$\\[2pt]
10 & $(0,1)$ & $11^{2},\allowbreak 12$ & $8$ & $- 2 K^{2} P U_{1} x + \frac{2 P^{3} U_{1}}{x^{3}}$\\[2pt]
11 & $(0,1)$ & $11,\allowbreak 12^{2}$ & $4$ & $- E_{(2)} K x^{2} + E_{(2)} P + K P^{2} U_{2} - \frac{P^{3} U_{2}}{x^{2}}$\\[2pt]
12 & $(0,1)$ & $11,\allowbreak 12,\allowbreak 22$ & $4$ & $E_{(2,1)} K^{2} P x + \frac{E_{(2,1)} K P^{2}}{x} + K P^{2} U_{3} x + \frac{P^{3} U_{3}}{x}$\\[2pt]
13 & $(0,1)$ & $12^{3}$ & $6$ & $3 E_{(2,1)} P x + 2 E_{(3)} x^{2} + \frac{P^{3} U_{3}}{x}$\\[2pt]
14 & $(0,1,1)$ & $11,\allowbreak 12,\allowbreak 13$ & $4$ & $K P^{2} U_{1}^{2} - K x^{2} - \frac{P^{3} U_{1}^{2}}{x^{2}} + P$\\[2pt]
15 & $(0,1,1)$ & $11,\allowbreak 12,\allowbreak 23$ & $2$ & $K P^{2} U_{1} U_{2} x + \frac{P^{3} U_{1} U_{2}}{x}$\\[2pt]
16 & $(0,1,1)$ & $12^{2},\allowbreak 13$ & $2$ & $E_{(2)} P U_{1} x + \frac{P^{3} U_{1} U_{2}}{x}$\\[2pt]
17 & $(0,1,1,1)$ & $12,\allowbreak 13,\allowbreak 14$ & $6$ & $\frac{P^{3} U_{1}^{3}}{x} + 3 P U_{1} x + 2 x^{2}$\\[2pt]
18 & $(0,0)$ & $11^{2},\allowbreak 12,\allowbreak 22$ & $16$ & $- 2 K^{3} P x^{2} - 2 K^{2} P^{2} + \frac{2 K P^{3}}{x^{2}} + \frac{2 P^{4}}{x^{4}}$\\[2pt]
19 & $(0,0)$ & $11,\allowbreak 12^{3}$ & $12$ & $- 6 K P x^{2} + \frac{2 P^{4}}{x^{4}} + 4 P x^{2}$\\[2pt]
20 & $(0,0)$ & $11,\allowbreak 12^{2},\allowbreak 22$ & $16$ & $K^{2} P^{2} + K^{2} x^{4} - \frac{2 K P^{3}}{x^{2}} - 2 K P x^{2} + \frac{P^{4}}{x^{4}} + P^{2}$\\[2pt]
21 & $(0,0)$ & $12^{4}$ & $48$ & $- K^{4} - 6 K^{2} + 16 K + \frac{P^{4}}{x^{4}} + 6 P^{2} + 32 P x^{2} + 9 x^{4} - 9$\\[2pt]
22 & $(0,0,1)$ & $11,\allowbreak 12^{2},\allowbreak 23$ & $4$ & $\frac{K P^{3} U_{1}}{x} - K P U_{1} x^{3} - \frac{P^{4} U_{1}}{x^{3}} + P^{2} U_{1} x$\\[2pt]
23 & $(0,0,1)$ & $11,\allowbreak 12,\allowbreak 13,\allowbreak 22$ & $4$ & $K^{2} P^{2} U_{1} x - \frac{P^{4} U_{1}}{x^{3}}$\\[2pt]
24 & $(0,0,1)$ & $11,\allowbreak 12,\allowbreak 23^{2}$ & $4$ & $E_{(2)} K P x^{2} + E_{(2)} P^{2} + K P^{3} U_{2} + \frac{P^{4} U_{2}}{x^{2}}$\\[2pt]
25 & $(0,0,1)$ & $12^{3},\allowbreak 13$ & $6$ & $- \frac{P^{4} U_{1}}{x^{3}} + 3 P^{2} U_{1} x + 4 P U_{1} x^{3}$\\[2pt]
26 & $(0,0,1)$ & $12^{2},\allowbreak 13,\allowbreak 23$ & $4$ & $E_{(2)} K^{2} + E_{(2)} + \frac{P^{4} U_{2}}{x^{2}} + P^{2} U_{2} x^{2}$\\[2pt]
27 & $(0,1,0)$ & $11,\allowbreak 12,\allowbreak 23,\allowbreak 33$ & $8$ & $4 E_{(2)} + K^{2} P^{2} U_{2} x^{2} + 2 K P^{3} U_{2} + \frac{P^{4} U_{2}}{x^{2}}$\\[2pt]
28 & $(0,0,1,1)$ & $11,\allowbreak 12,\allowbreak 23,\allowbreak 24$ & $4$ & $K P^{3} U_{1}^{2} + K P x^{2} + \frac{P^{4} U_{1}^{2}}{x^{2}} + P^{2}$\\[2pt]
29 & $(0,0,1,1)$ & $12^{2},\allowbreak 13,\allowbreak 24$ & $4$ & $K^{2} + \frac{P^{4} U_{1}^{2}}{x^{2}} + P^{2} U_{1}^{2} x^{2} + 1$\\[2pt]
30 & $(0,0,0)$ & $11,\allowbreak 12^{2},\allowbreak 23,\allowbreak 33$ & $8$ & $K^{2} P^{3} - K^{2} P x^{4} - \frac{P^{5}}{x^{4}} + P^{3}$\\[2pt]
31 & $(0,0,0)$ & $11,\allowbreak 12,\allowbreak 13,\allowbreak 22,\allowbreak 33$ & $16$ & $K^{3} P^{2} x^{2} + K^{2} P^{3} - \frac{K P^{4}}{x^{2}} - 4 K x^{2} - \frac{P^{5}}{x^{4}} + 4 P$\\[2pt]
32 & $(0,0,0)$ & $11,\allowbreak 12,\allowbreak 13,\allowbreak 23^{2}$ & $8$ & $- K^{3} x^{2} + K^{2} P + \frac{K P^{4}}{x^{2}} + K P^{2} x^{2} - K x^{2} - \frac{P^{5}}{x^{4}} - P^{3} + P$\\[2pt]
33 & $(0,0,0)$ & $11,\allowbreak 12,\allowbreak 23^{3}$ & $12$ & $- \frac{K P^{4}}{x^{2}} + 3 K P^{2} x^{2} + 4 K P x^{4} - \frac{P^{5}}{x^{4}} + 3 P^{3} + 4 P^{2} x^{2}$\\[2pt]
34 & $(0,0,0)$ & $12^{2},\allowbreak 13^{2},\allowbreak 23$ & $8$ & $- \frac{P^{5}}{x^{4}} + 2 P^{3} - P x^{4}$\\[2pt]
35 & $(0,0,0,1)$ & $11,\allowbreak 12,\allowbreak 23^{2},\allowbreak 34$ & $4$ & $\frac{K P^{4} U_{1}}{x} + K P^{2} U_{1} x^{3} + \frac{P^{5} U_{1}}{x^{3}} + P^{3} U_{1} x$\\[2pt]
36 & $(0,0,0,1)$ & $11,\allowbreak 12,\allowbreak 23,\allowbreak 24,\allowbreak 33$ & $8$ & $K^{2} P^{3} U_{1} x + \frac{2 K P^{4} U_{1}}{x} + \frac{P^{5} U_{1}}{x^{3}} + 4 P U_{1} x$\\[2pt]
37 & $(0,0,0,1)$ & $12^{2},\allowbreak 13,\allowbreak 23,\allowbreak 34$ & $4$ & $K^{2} P U_{1} x + \frac{P^{5} U_{1}}{x^{3}} + P^{3} U_{1} x + P U_{1} x$\\[2pt]
38 & $(0,0,0,0)$ & $11,\allowbreak 12,\allowbreak 23^{2},\allowbreak 34,\allowbreak 44$ & $16$ & $K^{2} P^{4} + K^{2} P^{2} x^{4} + 4 K^{2} + \frac{2 K P^{5}}{x^{2}} + 2 K P^{3} x^{2} + \frac{P^{6}}{x^{4}} + P^{4} + 4$\\[2pt]
39 & $(0,0,0,0)$ & $11,\allowbreak 12,\allowbreak 23,\allowbreak 24,\allowbreak 33,\allowbreak 44$ & $48$ & $K^{3} P^{3} x^{2} + 3 K^{2} P^{4} + \frac{3 K P^{5}}{x^{2}} + 12 K P x^{2} + \frac{P^{6}}{x^{4}} + 12 P^{2} + 16 x^{2}$\\[2pt]
40 & $(0,0,0,0)$ & $11,\allowbreak 12,\allowbreak 23,\allowbreak 24,\allowbreak 34^{2}$ & $8$ & $K^{3} P x^{2} + K^{2} P^{2} + \frac{K P^{5}}{x^{2}} + K P^{3} x^{2} + K P x^{2} + \frac{P^{6}}{x^{4}} + P^{4} + P^{2}$\\[2pt]
41 & $(0,0,0,0)$ & $12^{2},\allowbreak 13,\allowbreak 24,\allowbreak 34^{2}$ & $16$ & $K^{4} + 6 K^{2} + \frac{P^{6}}{x^{4}} + 2 P^{4} + P^{2} x^{4} + 5$\\[2pt]
42 & $(0,0,0,0)$ & $12,\allowbreak 13,\allowbreak 14,\allowbreak 23,\allowbreak 24,\allowbreak 34$ & $24$ & $3 K^{2} + 6 K P x^{2} + \frac{P^{6}}{x^{4}} + 8 x^{2} + 6$\\[2pt]
\end{longtable}\endgroup

\end{document}